\documentclass[11pt]{article}
\usepackage{float}
\usepackage{amsmath,amssymb,amsthm,graphicx,mathrsfs,enumerate}
\newcounter{probs}
\usepackage[margin=0.9in]{geometry}
\usepackage{graphicx}
\allowdisplaybreaks
\usepackage[ruled,longend]{algorithm2e}
\usepackage{enumitem}
\usepackage{dsfont}
\usepackage{graphicx}
\usepackage{subcaption}
\usepackage{mathtools}

\usepackage{hyperref}

\newlist{steps}{enumerate}{1}
\setlist[steps, 1]{label = Step \arabic*:}

\usepackage{xcolor}
\usepackage[normalem]{ulem}

\newcommand{\cturc}[1]{{\color{black} #1}}

\begin{document}

\topmargin -.5in
\oddsidemargin 0pt
\textheight 8.8in
\textwidth 6.5in
\title{High-Order Nystr\"om/Convolution-Quadrature Solution of Time-Domain Scattering from Closed and Open Lipschitz Boundaries with Dirichlet and Neumann Boundary Conditions}
\author{Erli Wind-andersen, Peter G. Petropoulos, Catalin Turc}

\setcounter{MaxMatrixCols}{14}
\newtheorem{theorem}{Theorem}[section]
\newtheorem{lemma}[theorem]{Lemma}
\newtheorem{proposition}[theorem]{Proposition}
\newtheorem{corollary}[theorem]{Corollary}
\newtheorem{remark}[theorem]{Remark}

\setcounter{footnote}{0}
\date{}
\newcommand{\triple}[1]{{\left\vert\kern-0.25ex\left\vert\kern-0.25ex\left\vert #1 
    \right\vert\kern-0.25ex\right\vert\kern-0.25ex\right\vert}}

\bibliographystyle{plain}
\maketitle
\begin{abstract}
 We investigate high-order Convolution Quadratures methods for the solution of the wave equation in unbounded domains in two and three dimensions that rely on Nystr\"om discretizations for the Boundary Integral Equation formulations of the ensemble of associated Laplace domain modified Helmholtz problems. Both Dirichlet and Neumann boundary conditions, imposed on open-arc/open surfaces as well  as Lipschitz closed scatterers, are considered. Two classes of CQ discretizations are employed, one based on linear multistep methods and the other based on Runge-Kutta methods,  in conjunction with Nystr\"om discretizations based on Alpert and QBX quadratures of Boundary Integral Equation (BIE) formulations of the Laplace domain Helmholtz problems with complex wavenumbers. A variety of accuracy tests are presented that showcase the high-order in time convergence (up to and including fifth order) that the Nystr\"om CQ discretizations are capable of delivering and we compare to numerical results in the literature pertaining to time-domain multiple scattering problems solved with other methods.
   
\end{abstract}

\section{Introduction}
\label{intro}
Numerical solutions of the wave equation based on retarded potentials boundary integral equations~\cite{et1986formulation,yilmaz2004time,banjai2009rapid,barnett2020high,anderson2020high} enjoy certain advantages over volumetric solvers such as FDTD and FEM, especially when the wave equation is posed in unbounded domains. In particular, Convolution Quadratures (CQ) have emerged in the last two decades as a powerful machinery for the numerical solution of time domain boundary integral equations. Indeed, the CQ provide a means to extend boundary integral equation (BIE) based frequency domain solvers to time domain BIE solvers via Laplace transforms. We present in this paper a CQ based methodology to extend high-order Nystr\"om BIE frequency domain solvers to high-order solvers of the time domain equation in the exterior of two dimensional obstacles that exhibit geometric singularities (e.g. corners, arc tips). 

Whenever applicable, BIE solvers are attractive choices for the solution of wave equations in the frequency domain (e.g. Helmholtz equation), owing to the dimensional reduction they provide and the explicit enforcement of radiation conditions when dealing with unbounded regions. In addition, high-order Nystr\"om BIE solvers for Helmholtz equations in geometrically complex two dimensional geometries (including non-smooth geometries) have been proposed in the literature~\cite{dominguez2008dirac,dominguez2012fully,dominguez2014nystrom,dominguez2015fully,helsing2008corner,bremer2012nystrom,ColtonKress,hao2014high,klockner2013quadrature} and can be more efficient than BEM solvers. Most of these BIE solvers were designed to handle the case of Helmholtz equations with \emph{real} wavenumbers, which is the most relevant in physical applications. CQ methods for the solution of the time-dependent wave equation, on the other hand, require the solution of an ensemble of \emph{modified} Helmholtz boundary value problems, that is Helmholtz equations with \emph{complex} wavenumbers, and the size of the ensemble grows with the order in time that CQ achieve. We focus in this paper on the second order in time BDF2 CQ~\cite{banjai2010multistep} and the third and fifth order in time Runge Kutta CQ~\cite{banjai2011runge}, whose algorithmic details are presented in those references. BIE solvers can be incorporated into the CQ machinery to solve wave equations provided their underpinning quadratures can handle integral kernels that feature complex wavenumbers. As it turns out, the wavenumbers associated with the CQ modified Helmholtz equations may have very large (positive) imaginary parts which become larger as the order of the CQ time integrator increases. 

While some of the high-order Helmholtz BIE quadratures proposed in the literature are oblivious to the nature of the wavenumbers, others struggle to resolve the additional challenges that complex wavenumbers with large imaginary parts pose. For instance, and as illustrated in this paper, the popular Kussmaul Martensen singularity splitting high-order quadratures~\cite{kusmaul,martensen}, while very effective in the case of real wavenumbers, cannot be extended in a straightforward manner to handle the complex wavenumbers countered in the CQ context. A method that does not suffer from this issue is the delta-bem method presented in~\cite{dominguez2008dirac,dominguez2012fully,dominguez2014nystrom,dominguez2015fully} which produces third order quadratures schemes for the discretization of the four boundary integral operators needed for BIE solutions of Helmholtz boundary value problems in two dimensions, irrespective of whether the wavenumbers are real or complex. Owing to this desirable feature, the delta-bem quadratures have been incorporated into high-order CQ solvers of the wave equation~\cite{dominguez2012fully,dominguez2014nystrom,dominguez2015fully}, and the code has been freely available for almost a decade~\cite{sayas2013retarded} serving as an important computational platform in the CQ community.

 We focus in this paper on two different types of high-order Nystr\"om quadratures that can also handle seamlessly discretizations of Helmholtz BIOs with complex wavenumbers: (1) one that relies on Alpert quadratures~\cite{alpert1999hybrid} for weakly singular kernels, and (2) the other which is the increasingly popular Quadratures by Expansion (QBX) introduced in~\cite{klockner2013quadrature}. Unlike the delta-bem methods, both aforementioned quadratures (i) do not rely on global interpolation and therefore are amenable to  panel discretizations and (ii) can be extended to three dimension configurations (albeit the Alpert quadrature can only handle axisymmetric scatterers~\cite{hao2014high}). Most of the delta-bem and Alpert quadrature literature is concerned with smooth boundaries, and therefore we focus mostly in this paper on two dimensional Lipschitz boundaries. High-order Nystr\"om discretizations of the Helmholtz BIOs associated with Lipschitz curves have been proposed~\cite{bremer2012nystrom,helsing2008corner,dominguez2015well}, and they resolve the singularities of the BIO functional densities  through the use of graded meshes that allow for the definition of weighted densities that are more regular. We follow this program and we present in this paper a simple extension of Alpert quadratures that produces high-order Nystr\"om discretizations of Helmholtz BIOs corresponding to Lipschitz boundaries via sigmoid graded meshes~\cite{ColtonKress}. In addition, we present a new version of QBX methods based on Chebyshev meshes  and F\'ejer quadratures that leads to high-order discretizations of Helmholtz BIOs in non smooth (including open) boundaries. An attractive feature of QBX methods based on Chebyshev meshes is their natural treatment of the corner and/or open end singularities of the Dirichlet and Neumann traces of Helmholtz solutions on non-smooth boundaries.
 
The high-order discretizations of Helmholtz BIOs with complex wavenumbers lead to high-order BIE solvers for Helmholtz problems, which can then be incorporated in the (blackbox) CQ machinery to deliver high-order in time solutions of the wave equation. While the choice of the BIE formulations can affect the accuracy of CQ schemes~\cite{betcke2017overresolving}, we focus on the simplest possible formulations, that is the single layer and respectively the double layer formulations in the case of Dirichlet and respectively Neumann boundary conditions. Besides the fact that both of these formulations are well posed for the CQ Laplace domain problems, we prefer them because they can be equally applied to both closed and open boundaries. In the case of Neumann boundary conditions the formulation we chose entails the Helmholtz hypersingular operator whose discretization is challenging for Nystr\"om discretizations. We investigate an alternative formulation of the hypersingular operator which is the Current and Charge (CC) formulation~\cite{taskinen2006current,helsing2020extended}. The CC formulation uses two functional densities instead of one in order to bypass the need to effect numerical differentiation, and in addition it involves Cauchy Principal Value BIOs. Nevertheless, we show that the QBX Nystr\"om discretization can still deliver high-order discretizations of the Current and Charge (CC) formulations on closed and open Lipschitz boundaries.

Our BIE based CQ solvers can be extended to solve time domain scattering problems involving axisymmetric scatterers in three dimensions. Indeed, we follow the methodology in~\cite{young2012high,epstein2019high} which leverages the Fourier modal representation of the rotationally invariant Green functions~\cite{cohl1999compact,conway2010exact} to extend two dimensional BIE Nystr\"om solvers to three dimensional axi symmetric solvers.  Our three dimensional solvers can handle both closed and open axisymmetric scatterers. We mention that a different class of Nystr\"om discretization BIE based CQ time domain solvers was recently introduced in~\cite{ganesh2023high}. The Nystr\"om solvers in~\cite{ganesh2023high} rely on the addition theorem for Helmholtz Green's functions and are applicable to scatterers which are (globally) diffeomorphic to spheres. As such, the axisymmetric scatterers we considered in this paper, e.g. torii and annuli, fall beyond the scope of the CQ solvers in~\cite{ganesh2023high}. The extension to three dimensional surfaces of the CQ methods based on QBX BIE discretizations will be presented in a future publication.

The paper is organized as follows: in Section~\ref{briefrev} we summarize the application of the BDF2 and RK CQ to the wave equation in two dimensions; in Section~\ref{Helmholtz} we present Nystr\"om discretizations of the ensemble of Laplace domain modified Helmholtz problems related to BDF2 and RK CQ methods for the wave equation; finally, in Section~\ref{Numerical_results} we present numerical results that illustrate the high order in time (up and including fifth order) achieved by CQ solutions of the wave equation in the exterior of geometrically complex two and three dimensional obstacles.

\subsection{Brief Review of Convolution Quadratures for TDBIE}\label{briefrev}
\cturc{We are interested in solving the wave equation in an infinite domain $\Omega^{+}:=\mathbb{R}^d\setminus \overline{\Omega}$, where $\Omega$ is a bounded domain in $\mathbb{R}^d$ with $d=2,3$:
\begin{equation}\label{eq:wave}
	\begin{cases}
	c^{-2}\frac{\partial^2 u}{\partial t^2} &= \Delta u \quad {\rm in}\  \Omega^{+} \times (0,\infty)\\
	u(x,0) &= \frac{\partial u(x,0)}{\partial t } = 0 \quad {\rm in}\ \Omega^{+}\\
	u(x,t) &= g(x,t) \quad {\rm on}\ \Gamma\times (0,\infty), 
	\end{cases}
\end{equation}
and we denoted $\Gamma:=\partial\Omega $. Using the fundamental solution of the wave equation
\[
k_{2D}(x,t):=\frac{H(t-c^{-1}|x|)}{2\pi\sqrt{t^2-c^{-2}|x|^2}}\qquad k_{3D}(x,t):=\frac{\delta(t-c^{-1}|x|)}{4\pi|x|}
\]
we seek for a solution $u(x,t)$ of the wave equation in the form of a Time Domain Single-Layer (TDSL) potential:
\begin{equation*}
    u(x,t)= (S*\lambda)(x,t) =\int_{\Gamma}\int_0^t k(x-y,t-\tau)\varphi(y,\tau)d s(y)d\tau
\end{equation*}
in terms of the unknown density function $\varphi$ defined on $\Gamma\times(0,\infty)$ where the kernel $k(x,t)$ is either the two dimensional or three dimensional fundamental solution of the wave equation. Imposing the Dirichlet boundary conditions leads to the following Time \cturc{Domain} Boundary Integral Equation (TDBIE) 
\begin{equation}\label{eq:TDBIE}
\int_{\Gamma}\int_0^t k(x-y,t-\tau)\varphi(y,\tau)d s(y)d\tau=g(x,t)\quad {\rm on}\ \Gamma\times (0,\infty).
\end{equation}
We assume in what follows that the wave speed is $c=1$ and we solve the TDBIE~\eqref{eq:TDBIE} using the Convolution Quadrature (CQ) methodology which we briefly review in the remainder of this Section.  concerned exclusively with the algorithmic aspects of CQ and is based on~\cite{banjai2009rapid,betcke2017overresolving,sayas2013retarded}.} Also, it

\subsubsection{Convolution Quadratures for TDBIE}\label{cq}
CQ are hybrid time stepping methods that combine Laplace transforms and A-stable ODE solvers to produce approximations of convolution integrals
\begin{equation}
	(f*g)(t) \coloneqq \int_0^t f(\tau)g(t-\tau)d\tau,
	\label{cq1}
\end{equation}
for values of $t$ on an equidistant mesh  $0=t_0\leq\dots\leq t_{N_t}=N_t\Delta t:=T$. In a nutshell, multistep CQ methods apply A-stable multistep solvers (e.g. the second order in time multistep method BDF2 which is considered in this paper) to  deliver a discrete convolution
$$(f*g)(t_n)\approx\sum_{m=0}^n \omega_m(\Delta t)g_{n-m} \quad \text{ where   }g_n = g(t_n)$$
that approximates the continuous convolution $(f*g)(t_n), 0\leq n\leq N_t$. The convolution weights $\omega_m(\Delta t)$ in the formula above are defined as
\begin{equation}
	F\left( \frac{P(\zeta)}{\Delta t}\right )= \sum_{m=0}^\infty \omega_m(\Delta t )\zeta^{m} \quad \text{ where   } \omega_m(\Delta t) = \frac{1}{m!}\frac{d^m}{d\zeta^m} \left ( F\left( \frac{P(\zeta)}{\Delta t} \right) \right)\Big|_{\zeta=0},
	\label{cq6}
\end{equation}
where where $P(\zeta) = \frac{1}{2} (\zeta^2-4\zeta+3)$ is the characteristic polynomial of BDF2. 
Runge-Kutta based CQ methods can achieve orders of convergence in time higher than three when applied to the wave equation~\cite{banjai2010multistep,banjai2011runge}. Our exposition of these methods follows closely the one in~\cite{betcke2017overresolving}. Applying an $m$ stage Runge-Kutta (RK) scheme with the associated Butcher Tableau consisting of the ensemble of matrices and vectors
\begin{equation*}
    A= [a_{ij}]_{1\leq i,j\leq m }\quad b = [b_j]_{1\leq j \leq m} \quad c=[c_j]_{1\leq j\leq m},
\end{equation*}
to the first order system reformulation of the wave equation, one arrives via the $\zeta$-transforms at the following ensemble of Laplace domain modified Helmholtz equations
\begin{equation}\label{eq:CQh2}
    \begin{cases}
        \Delta_x W_j(x;\zeta) - \left ( \frac{\gamma_j(\zeta)}{\Delta t} \right)^2 W_j(x;\zeta) =0 \quad &x\in \Omega^+ \\
        W_j(x;\zeta)= \sum_{\ell=1}^m\mathds P^{-1}_{j\ell}(\zeta)G_\ell(x;\zeta)\quad x\in\Gamma
    \end{cases}
\end{equation}
where $G_\ell(x;\zeta)$ is the following $\zeta$-transform
\begin{equation}\label{eq:Gl}
 G_\ell(x;\zeta):=\sum_{n \geq 0}g(x;t_n + c_\ell\Delta t)\zeta^n\quad x\in \Gamma.
\end{equation}
We remark that the Laplace domain Helmholtz equations featured in equations~\eqref{eq:CQh2} are all of the form
\begin{equation}\label{eq:modH1}
	\Delta_x w(x;s) -s^2 w(x;s)= 0 \quad x \in  \Omega^+ 
\end{equation}
and thus we can seek the radiative fields $w(x;s)$ in the form of single layer potentials
\[
w(x;s):=\int_{\Gamma} K(x-y,s)\varphi(y) ds(y),\quad x\in \mathbb{R}^d\setminus\Gamma
\]
where the kernels $K(x,s)$ denotes the Laplace transform of the wave equation Green functions $k_{2D}(x,t)$ and respectively $k_{3D}(x,t)$ with respect to the variable $t$, that is
$$ K_{2D}(x,s)= \frac{i}{4}H_0^{(1)}(is|x|)\qquad K_{3D}(x,s)=\frac{e^{-s|x|}}{4\pi|x|}.$$
The enforcement of the boundary conditions in equations~\eqref{eq:CQh2} leads to the solution of BIE featuring the single layer BIO $V(s)$ associated with a modified Helmholtz equation
\begin{equation}
	\left( V(s)\varphi   \right )(x) = \int_{\Gamma} K(x-y,s)\varphi(y) ds(y),\quad x\in\Gamma.
	\label{cq19}
\end{equation}
Once the densities $\varphi$ associated with the single layer BIE formulations of the ensemble of Helmholtz equations~\eqref{eq:CQh2} are obtained (and hence the solutions of the Helmholtz equations $W_j(\cdot;\zeta)$ themselves), the $\zeta$-transform of the solution of the wave equation
\[
U_{d}(x;\zeta):=\sum_{n\geq 0}u_d(x;t_n)\zeta^n
\]
is in turn immediately retrieved from the formula
\begin{equation}
    U_d(x;\zeta) = \zeta\sum_{j=1}^m\mathds P_{mj}(\zeta)W_j(x;\zeta).
\end{equation}
In equations~\eqref{eq:CQh2} the wavenumbers $\gamma_j(\zeta)$ are related to the diagonalization of the matrix valued function $\Delta(\zeta)$ defined by 
\begin{equation*}
    \Delta(\zeta) = \left( A+\frac{\zeta}{1-\zeta}\mathds 1  b^T     \right )^{-1}
\end{equation*}
where $\mathds 1 = [1,\ldots,1]^T$. Specifically, $\mathds P(\zeta)$ is the matrix consisting of the eigenvectors of $\Delta(\zeta)$ and $\mathds D(\zeta)=diag(\gamma_1(\zeta),\ldots,\gamma_m(\zeta))$ is the diagonal matrix with the corresponding eigenvalues of $\Delta(\zeta)$ so that $\Delta(\zeta)=\mathds P(\zeta)\mathds D(\zeta)\mathds P(\zeta)^{-1}$. Finally, reverting from the Laplace domain to the physical domain is performed via FFTs~\cite{banjai2010multistep,banjai2011runge}.

We will use in our numerical solutions two-stage and three-stage Implicit RK Radau IIa CQ methods. The two-stage RK Radau IIa method gives rise to third order in time CQ and will be denoted by the acronym RK3. The three-stage RK Radau IIa method gives rise to fifth order in time CQ and will be denoted by the acronym RK5. We note that the higher order in time convergence that can be achieved by RK CQ  methods entail commensurately more solutions of frequency domain modified Helmholtz equations.

\section{Nystr\"om discretizations for the solution of modified Helmholz BIE}\label{Helmholtz}

As we have already remarked, a key component of CQ methods for the solution of TDBIE is solving ensembles of modified Helmholtz equations~\eqref{eq:CQh2} in unbounded domains.  We will use Nystr\"om discretizations for the solution of the BIE formulations of those modified Helmholtz problems that involve complex wavenumbers. Whereas in the case of BIE associated with Helmholtz equations with real wavenumbers several Nystr\"om discretization strategies exist in the literature (see for example~\cite{hao2014high} for an excellent review of some of those strategies), the case of Helmholtz equations with complex wavenumbers has received less attention. We discuss in what follows the challenges that are encountered in extending Nystr\"om discretizations to the case of BIOs with complex wavenumbers and we present several methods that can deal seamlessly with ensemble evaluation of Helmholtz BIOs for a wide range of complex wavenumbers, and for Lipschitz and smooth boundaries (both closed and open) in two and three dimensions. 

\subsection{Dirichlet Boundary Conditions, two dimensional case} 
We note that the CQ solution of TDBIE~\eqref{eq:TDBIE} presented in Section~\ref{cq} requires in turn the solution of Helmholtz problems of the type~\eqref{eq:modH1} via associated single layer operators~\eqref{cq19}.  Therefore, we discuss in what follows the single layer potential approach for the solution of Helmholtz problems in unbounded domains with Dirichlet boundary conditions in two dimensions. Specifically, we seek to solve the following two dimensional Helmholtz problem with wavenumber $k\in\mathbb{C}$
\begin{equation}
    \begin{cases}
        \Delta u + k^2u = 0 \quad &x\in \Omega^+ \\
        u(x) =g(x) \quad & x \in \Gamma \\
        u(x) = \sum_{m=-\infty}^\infty u_m H_m^{(1)}(k|x|)e^{im\theta},\quad &|x|>r_0,\quad x=|x|(\cos{\theta},\sin{\theta}).
        \label{DIRICH}
    \end{cases}
\end{equation}
We note that the last condition at infinity imposed in equations~\eqref{DIRICH} is equivalent to the Sommerfeld radiation condition $ \lim_{|x|\to \infty} \sqrt |x| (\frac{\partial}{\partial |x|} -ik)u(|x|,\theta)=0$  in the case when $k\in\mathbb{R}^{+}$.

For wavenumbers $k\in\mathbb{C}$ such that $\Im{k}>0$ (which turn out to be all of the wavenumbers involved in the Helmholtz problems~\eqref{eq:modH1} and~\eqref{eq:CQh2}) we seek the solution $u$ of the Helmholtz problem~\eqref{DIRICH} in the form of the single layer potential
\begin{equation*}
  u(x)=SL(k)[\varphi](x):=\int_\Gamma G_k(x-y)\varphi(y)ds(y),x\in\mathbb{R}^2\setminus\Gamma, \quad G_k(z)=\frac{i}{4}H_0^{(1)}(k|z|),\ z\neq 0,
\end{equation*}
and thus the boundary density $\varphi$ solves the BIE of the first kind
\begin{equation}\label{eq:SL}
  V(k)[\varphi](x):=\int_\Gamma G_k(x-y)\varphi(y)ds(y)=g(x),\quad x\in\Gamma.
\end{equation}
We note that the BIE~\eqref{eq:SL} is uniquely solvable for all wavenumbers $k$ such that $\Im{k}>0$, and all types of boundaries $\Gamma$, be the boundaries closed (Lipschitz) or open. Given its versatility, the single layer BIE formulation~\eqref{eq:SL} seems to be the preferred method of solution in the CQ literature~\cite{banjai2009rapid,banjai2010multistep,sayas2013retarded}. Furthermore, coercivity estimates on the operators $V(k):H^{-1/2}(\Gamma)\to H^{1/2}(\Gamma)$ are available when $\Im{k}>0$~\cite{et1986formulation} and they greatly facilitate the error analysis~\cite{banjai2009rapid,sayas2013retarded} of CQ solution of TDBIE~\eqref{eq:TDBIE}. We note here that in the case when $\Gamma$ is closed, a variety of other BIE formulations are available for the solution of the Helmholtz problems~\eqref{DIRICH}, such as the double layer or the combined field approach, and their influence on the accuracy of the CQ schemes in connection with the location of the complex poles of the BIE formulation was discussed in~\cite{betcke2017overresolving}.

The main challenge faced by Nystr\"om discretizations of the single layer formulations of the type~\eqref{eq:SL} is the resolution of the logarithmic singularity of the integral kernel. To this end, several singular quadratures have been proposed in the literature. For instance, the Kussmaul-Martensen quadrature~\cite{kusmaul,martensen} relies on global trigonometric interpolation of the functional density $\varphi$, and a singular splitting of the kernels in the form
\begin{equation}\label{eq:splitH}
\frac{i}{4}H_0^{(1)}(k|z|)=-\frac{1}{2\pi}J_0(k|z|)\log{|z|}+C_k(z)
\end{equation}
where $C_k(z)$ is a smooth function. Unfortunately, the Kussmaul-Martensen quadratures cannot be applied directly in the case when the wavenumber $k$ is complex with positive imaginary part because the Bessel function $J_0(k|z|)$ in the logarithmic splitting~\eqref{eq:splitH} grows exponentially as $|z| \to \infty$. This setback can be remedied in various ways, either using cutoff functions to mollify the exponential growth of the Bessel functions~\cite{turc1} or carefully modifying the splitting formula~\eqref{eq:splitH}~\cite{lu2014efficient}. However, neither of these approaches are particularly suited for use with CQ to solve time-domain problems. For this reason, we use alternative quadrature schemes that can be applied seamlessly to the resolution of general logarithmic singularities.  

\subsubsection{Alpert Quadratures}\label{alpquad}
In contrast to the Kusmaul-Martensen quadratures that rely on the logarithmic splitting techniques described above, Alpert quadratures~\cite{alpert1999hybrid} can handle the same kernel singularities in a manner that is agnostic to the argument of the Hankel functions. We note that the use of Alpert quadratures in the CQ context was already advocated in~\cite{Labarca2019convolution} where it was shown that these quadratures can handle in a seamless manner the discretization of the BIOs involved in the transmission analogues of the ensemble of Laplace domain Helmholtz problems~\eqref{eq:modH1} and respectively~\eqref{eq:CQh2}. We describe briefly in what follows the details of Alpert quadratures, following the exposition in~\cite{hao2014high}. In essence, these quadratures are designed to handle logarithmic singularities in the following manner
\begin{align*}
\int_0^T \rho(x_i, x')\mu(x')dx' & \approx h\sum_{p=0}^{N-2a} \rho(x_i,x_i+ah+ph)\mu(x_i+ah +ph) \\
    &+ h\sum_{p=1}^m w_p\rho(x_i,x_i+\chi_p h)\mu(x_i+\chi_p h ) \\
    & + h\sum_{p=1}^m w_p \rho(x_i,x_i+T-\chi_p h)\mu(x_i+T-\chi_p h ),
    \end{align*}
where $\{x_i:=ih\}$ is an equispaced mesh on the interval $[0,T]$, the kernel $\rho(x_i,x')$ has a logarithmic singularity as $x'\to x_i$, $\mu$ is assumed to be a regular enough density function, and the weights $w_p$ and the nodes $\chi_p$ are selected so that the ensuing quadratures achieve a prescribed rate of convergence provided that the integer quantities $a$ and $m$ are chosen appropriately~\cite{alpert1999hybrid}. The endpoint correction nodes $\chi_p$ are typically not integers, and as such the Alpert quadrature requires evaluation of the density $\mu$ outside of the equispaced mesh $\{x_i\}$. This is achieved by resorting to Lagrange interpolation of order $m+3$ with $m+4$ equispaced nodes around $x_i$ and $x_i+T$
\[
    \sigma(x) = \sum_{q=0}^{m+3}L_q^{(x_i)}(x) \mu(x_i+qh), \quad\quad \text{where} \quad\quad L_q ^{(x_i)}(x) = \prod_{r=0} \frac{x - (x_i+rh)}{(x_i+qh) -(x_i+rh)}.
    \]
Consequently, the Alpert quadrature takes on the explicit form
    \begin{align*}
    \int_0^T \rho(x_i,x')\mu(x')dx' &\approx h\sum_{p=0}^{N-2a}\rho(x_i,x_i+ah+ph)\mu(x_i+ah+ph) \\
    &+ h\sum_{q=0}^{m+3}\left ( \sum_{p=1}^m w_p\rho(x_i,x_i+\chi_ph) L_q^{(x_i)}(x_i+\chi_p) \right ) \mu(x_i+qh) \\
    & +h\sum_{q=0}^{m+3} \left ( \sum_{p=1}^m w_p\rho(x_i,x_i+T-\chi_p)L_q^{(x_i+T)}(x_i+T-\chi_p h)  \right ) \mu(x_i+T-qh),
    \end{align*}
which can be readily re-expressed in matrix-vector product form~\cite{hao2014high}. The Alpert quadratures described above can be applied directly to obtain high-order Nystr\"om discretizations of the parametrized versions of the BIE~\eqref{eq:SL} (because the kernels of the latter exhibit a logarithmic singularity) for any value of the wavenumber $k$ such that $\Im{k}\geq 0$. Clearly, the high-order of Alpert quadratures Nystr\"om discretizations of the BIE~\eqref{eq:SL} also depends on the regularity of the density function $\varphi$. Extensive comparisons of Alpert quadrature Nystr\"om discretizations with other high-order Nystr\"om discretizations were performed in~\cite{hao2014high} in the case when the boundary $\Gamma$ is a smooth closed curve. We present in what follows a simple extension of Alpert quadrature Nystr\"om discretizations of the single layer BIE~\eqref{eq:SL} to the case when the boundaries $\Gamma$ are either (a) Lipschitz and closed or (b) open arcs. To this end we rely on sigmoid transforms~\cite{KressCorner} and weighted formulations~\cite{dominguez2015well} as we explain next. We mention that this procedure is readily extended to the case of the double layer and combined field BIE formulations of the Helmholtz equation~\eqref{DIRICH}.  

We assume that the curve $\Gamma$ is piecewise smooth has corners at $x_1,x_2,\ldots,x_P$; in the case when $\Gamma$ is an open arc, we denote by $x_1$ and $x_P$ its end points. We further assume that the boundary curve $\Gamma$ has a $2\pi$ periodic parametrization so that each of the curved segments $[x_j,x_{j+1}]$ is parametrized by $\gamma(t)=(\gamma_1(w_j(t)), \gamma_2(w_j(t)))$ with $t\in[T_j,T_{j+1}]$ (so that $x_j=\gamma(T_j)$,  and $\gamma_j:\mathbb{R}\to\mathbb{R}$ are $2\pi$  periodic with $(\gamma_1'(t))^2+(\gamma_2'(t))^2>0$ for all $t$), where $0=T_0<T_1<T_2<\ldots <T_P<T_{P+1}=2\pi$ in the case when $\Gamma$ is closed, and $0=T_1<T_2<\ldots <T_P=2\pi$ when $\Gamma$ is an open arc. The maps $w_j:[T_j,T_{j+1}]\to[T_j,T_{j+1}],\ 0\leq j\leq P$, are the sigmoid transforms introduced by Kress~\cite{KressCorner}
\begin{eqnarray}\label{eq:cov_w}
w_j(s)&=&\frac{T_{j+1}[v_j(s)]^\sigma+T_j[1-v_j(s)]^\sigma}{[v_j(s)]^\sigma+[1-v_j(s)]^\sigma},\ T_j\leq s\leq T_{j+1},\ 0\leq j\leq P\\
v_j(s)&=&\left(\frac{1}{\sigma}-\frac{1}{2}\right)\left(\frac{T_{j}+T_{j+1}-2s}{T_{j+1}-T_j}\right)^3+\frac{1}{\sigma}\ \frac{2s-T_j-T_{j+1}}{T_{j+1}-T_j}+\frac{1}{2}\nonumber
\end{eqnarray} 
where $\sigma>2$.  The functions $w_j$ are smooth and increasing bijective functions on each of the intervals $[T_j,T_{j+1}]$ for $0\leq j\leq P$, with $w_j^{(k)}(T_j)=w_j^{(k)}(T_{j+1})=0$
for $1\leq k\leq \sigma-1$ and all $0\leq j \leq P$. With the aid of the graded meshes just introduced, we define \emph{weighted} densities $\varphi^w$ in the following manner 
\[
\varphi^w(t):=|\gamma'(t)|\varphi(\gamma(t)),\quad 0\leq t\leq 2\pi
\]
as well as \emph{weighted} parametrized single layer formulations
\begin{equation}\label{eq:SLw}
  \frac{i}{4}\int_0^{2\pi}H_0^{(1)}(k|\gamma(t)-\gamma(\tau)|)\varphi^w(\tau)d\tau=g(\gamma(t)),\quad t\in[0,2\pi].
\end{equation}
On account of the vanishing property of the sigmoid transforms at the breakpoints $T_j$, the weighted density $\varphi^w$ are $2\pi$ periodic functions whose regularity can be control by the value of $p$ in the graded mesh used. As such, we apply the Alpert quadrature described above to produce high-order Nystr\"om discretizations of the weighted single layer formulations~\eqref{eq:SLw}. Finally, in order to avoid the evaluation of the density function at corner points, we employ the equispaced mesh
$t_j = (j-1/2)\left(\frac{2\pi}{N}\right), 1\leq j\leq N$ in the Alpert Nystr\"om method. We present in Table~\ref{compA1} the errors in the near field achieved by the Alpert discretization of the weighted formulation~\eqref{eq:SLw} for real wavenumbers with parameters $a=2$ and $m=3$ for various Lipschitz and open arc $\Gamma$ configurations. Specifically, we considered the cases when $\Gamma$ is (i) the teardrop geometry $x(t)=\left(2\sin{\frac{t}{2}},-\beta\sin{t}\right),\ \beta=\tan{\frac{\alpha\pi}{2}},\ \alpha=1/2,\ 0\leq t\leq 2\pi$, (ii) the boomerang geometry $x(t)=\left(-\frac{2}{3}\sin{\frac{3t}{2}},-\sin{t}\right),\ \ 0\leq t\leq 2\pi$~\cite{ColtonKress} (both domains (i)-(ii) have diameters equal to 2), (iii) the strip connecting $(-1,0)$ to $(1,0)$, and (iv) the V-shaped strip connecting $(-1,1)$ to $(0,0)$ and then to $(1,1)$. We report errors in the near field evaluated at $512$ equispaced points placed on circles located at distance $1$ from the scatterers. We focus on the case when the wavenumber is $k=8$ under normal plane wave incidence, because in the case of real wavenumbers extremely accurate reference solutions can be produced via alternative spectral Nystr\"om solvers (e.g. those described in~\cite{ColtonKress} in the case of closed curves and respectively in~\cite{bruno2012second} in the case of open arcs). We report the estimated orders of convergence achieved by the Alpert methods under the acronym ``e.o.c.''.  We mention that similar convergence rates were observed in the case when the wavenumber $k$ is complex in the formulation~\eqref{eq:SLw}. 

\begin{table}
  \begin{center}
\begin{tabular}{|c|c|c|c|c|c|c|c|c|}
\hline
$N$ & \multicolumn{2}{c|} {Teardrop} &  \multicolumn{2}{c|} {Boomerang} & \multicolumn{2}{c|}{Strip} & \multicolumn{2}{c|}{V-shaped Strip}\\
\hline
 & $\varepsilon_\infty$ & e.o.c.  & $\varepsilon_\infty$ & e.o.c. & $\varepsilon_\infty$ &  e.o.c. & $\varepsilon_\infty$ &  e.o.c.\\
\hline
64 & 4.9 $\times$ $10^{-4}$ & & 1.9 $\times$ $10^{-3}$& & 4.1 $\times$ $10^{-6}$ & & 8.5 $\times$ $10^{-5}$ & \\
128 & 1.2 $\times$ $10^{-5}$ & 4.63 & 6.1 $\times$ $10^{-5}$ & 4.95 & 1.9 $\times$ $10^{-7}$ & 4.37 & 3.2 $\times$ $10^{-6}$ & 4.71\\
256 & 7.3 $\times$ $10^{-7}$ & 4.75 & 2.3 $\times$ $10^{-6}$& 4.72 & 1.2 $\times$ $10^{-8}$& 4.02& 1.3 $\times$ $10^{-7}$ & 4.57\\
512 & 2.8 $\times$ $10^{-8}$ & 4.66 & 9.0 $\times$ $10^{-8}$ & 4.67 & 7.9 $\times$ $10^{-10}$ & 3.93& 6.7 $\times$ $10^{-9}$ & 4.35\\
\hline
\end{tabular}
\caption{Errors in the near field and estimated orders of convergence corresponding to the Alpert discretization of the weighted formulation~\eqref{eq:SLw} for real wavenumber $k=8$ and plane wave normal incidence. We used the following Alpert quadrature parameters: $\sigma=4$ in the sigmoid transform, and respectively $a=2, m=3$. \label{compA1}}\end{center}
\end{table}

\subsubsection{QBX Quadratures}\label{qbxquad}
The quadrature by expansion (QBX) method is another powerful method to discretize the boundary layer potentials and operators that arise in Helmholtz problems~\cite{klockner2013quadrature}. Just like the Alpert quadratures, QBX can be straightforwardly applied to discretizations of layer potentials whose associated Green functions feature either real or complex wavenumbers. We review in what follows the details of the QBX method in the case of evaluations of single layer potentials and BIOs. The QBX method is based on the continuity of the single layer potential across the boundary $\Gamma$
\[
V(k)[\varphi](x)=\lim_{x^+\to x\in\Gamma} SL(k)[\varphi](x^+),
\]
and it relies on smooth extensions of the single layer potentials $SL(k)[\varphi](x^+),\ x^+\in\Omega^+$ to $\overline{\Omega^+}$ via certain series expansions. Specifically, the method uses expansion centers, that is for $x\in\Gamma$ we define $x^\pm=x\pm\varepsilon n(x)$, where the vector $n(x)$ is the exterior unit normal at $x\in\Gamma$ and the parameter $\varepsilon=\varepsilon(x)>0$ is small enough,  and relies on the addition theorem
\begin{equation}\label{eq:add_thm}
  H_0^{(1)}(k|x-x'|)=\sum_{\ell=-\infty}^{\infty}H_\ell^{(1)}(k|x'-x^+|)e^{i\ell\theta'}J_\ell(k|x-x^+|)e^{-i\ell\theta},\ x'\in\Gamma
\end{equation}
where $\theta$ and $\theta'$ are the angular coordinates of $x$ and respectively $x'$ in the polar coordinate system centered at $c=x^+$ to produce the following series representation for the single layer potential
\begin{equation}\label{eq:qbx}
  V(k)[\varphi](x)=\sum_{\ell=-\infty}^{\infty}\alpha_\ell[\varphi] J_\ell(k|x-x^+|)e^{-i\ell\theta},\quad \alpha_\ell[\varphi]:=\frac{i}{4}\int_\Gamma H_\ell^{(1)}(k|x'-x^+|)e^{i\ell\theta'}\varphi(x')ds(x').
\end{equation}
In practice a truncation parameter $p$ is chosen in the series representation~\eqref{eq:qbx} leading to the following QBX approximation
\begin{equation}\label{eq:qbxa}
  V(k)[\varphi](x)\approx\sum_{\ell=-p}^{p}\alpha_\ell[\varphi] J_\ell(k|x-x^+|)e^{-i\ell\theta}
\end{equation}
where the coefficients $\alpha_\ell$ have to be evaluated numerically for all $-p\leq \ell\leq p$. We note that the integrands in the definition~\eqref{eq:qbx} of the coefficients $\alpha_\ell$ do not contain kernel singularities as $x^+\notin\Gamma$. Owing to the fact that the evaluation of the expansion coefficients $\alpha_\ell$ requires integration on the whole boundary $\Gamma$, this version of QBX is referred to as global; local versions have also been designed~\cite{siegel2018local}, but we do not consider them herein. 

We will use in what follows panel Chebyshev meshes and Clenshaw-Curtis quadratures to evaluate numerically the coefficients $\alpha_\ell$ in equation~\eqref{eq:qbx}. Indeed, we assume that the piecewise smooth curve $\Gamma$ can be represented as the union of (disjoint) smooth panels $\Gamma_m $ in the form $\Gamma=\bigcup_{m=1}^M \Gamma_m$ so that the corner points $x_1,x_2,\ldots,x_P$ of $\Gamma$ are end points of these panels. We thus have
\begin{equation*}
  \alpha_\ell[\varphi]=\sum_{m}\alpha_{\ell,m}[\varphi],\qquad \alpha_{\ell,m}[\varphi]:=\frac{i}{4}\int_{\Gamma_m} H_\ell^{(1)}(k|x'-x^+|)e^{i\ell\theta'}\varphi(x')ds(x').
\end{equation*}
We apply Clenshaw-Curtis quadratures for the evaluation of each of the integrals in the definition of $\alpha_{\ell,m}$, that is
\begin{equation}\label{eq:Clensh-Curtis}
  \alpha_{\ell,m}[\varphi]\approx\frac{i}{4}\sum_{j=1}^{N_m}\omega_jH_\ell^{(1)}(k|\gamma_m(t_j)-x^+|)e^{i\ell\theta'_j}\widetilde{\varphi}(\gamma_m(t_j))
\end{equation}
where we assume that $\Gamma_m$ is parametrized in the form $\Gamma_m=\{\gamma_m(t): t\in[-1,1]\}$ and $\gamma_m:[-1,1]\to\mathbb{R}^2$ is smooth, $\widetilde{\varphi}(\gamma_m(t_j))=\varphi(\gamma_m(t_j))|\gamma'_m(t_j)|$, the quadrature points $t_j$ are the Chebyshev zero points
\begin{equation*}
  t_j:=\cos(\vartheta_j),\quad \vartheta_j:=\frac{(2j-1)\pi}{2N_m},\quad j=1,\ldots,N_m
\end{equation*}
and the Fej\'er quadrature weights $\omega_j$ are given by
\begin{equation*}
  \omega_j:=\frac{2}{N_m}\left(1-2\sum_{q=1}^{[N_m/2]}\frac{1}{4q^2-1}\cos(2q\vartheta_j)\right),\quad j=1,\ldots,N_m.
\end{equation*}
The steps of the QBX discretization algorithm of the single layer BIE~\eqref{eq:SL} are summarized in what follows: (1) choose a grid on $\Gamma$ such that $x_{j,m}=\gamma_m(t_j),\ 1\leq m\leq M,\ 1\leq j\leq n_m$ and $t_j:=\cos(\vartheta_j),\quad \vartheta_j:=\frac{(2j-1)\pi}{2n_m},\quad j=1,\ldots,n_m$; (2) for each $x_j$ on $\Gamma$ define the expansion centers $x_j^\pm:=x_j\pm\varepsilon(x_j)n(x_j)$; (3) use the QBX approximation~\eqref{eq:qbxa} for a given truncation parameter $p$; and (4) evaluate the ensuing coefficients $\alpha_{\ell,m}[\varphi]$ according to the Clenshaw-Curtis quadrature~\eqref{eq:Clensh-Curtis} using oversampling, that is choose $N_m=\beta n_m,\ \beta\geq 1$ in equation~\eqref{eq:Clensh-Curtis} and use Chebyshev interpolation to access the values of the density $\varphi^w$ on the fine Chebyshev grid with $N_m$ nodes on each $\Gamma_m$. Finally, we choose $\varepsilon(x_j)=\min(|x_j-x_{j-1}|,|x_j-x_{j+1}|)$ in the definition of the centers $x_j^\pm$.

It is also possible to consider weighted versions of the single layer BIO $V(k)$ by defining weighted functional densities that are more regular in the neighborhood of the geometric singularities of the curve $\Gamma$. Specifically we define the weighted densities
\begin{equation}\label{eq:weighted_phi}
  \varphi_m^w(t):=\varphi(\gamma_m(t))|\gamma'_m(t)|\sqrt{1-t^2},\quad -1\leq t\leq 1,\quad 1\leq m\leq M
\end{equation}
whose square root weights are chosen to resolve the worst case single layer density singularities that arise in the neighborhood of open endpoints of $\Gamma$ (in the case when $\Gamma$ is an open arc the solutions of the single layer BIE~\eqref{eq:SL} behave asymptotically as $\varphi\sim\mathcal{O}(1/\sqrt{d})$ in the neighborhood of the open ends of the arc, where $d$ is the distance to the open end of the arc). Incidentally, the incorporation of weighted densities~\eqref{eq:weighted_phi} in the definition of single layer BIO leads to a simpler Fej\'er quadrature rule for the evaluation of the expansion coefficients $\alpha_{\ell,m}$ in the form
\begin{equation}\label{eq:Clensh-Curtis-1}
  \alpha_{\ell,m}[\varphi^w]\approx\frac{i\pi}{4N_m}\sum_{j=1}^{N_m}H_\ell^{(1)}(k|\gamma_m(t_j)-x^+|)e^{i\ell\theta'_j} \varphi^w_m(t_j).
\end{equation}
In addition, the condition numbers of the QBX Nystr\"om matrices corresponding to the weighted expansion coefficients $\alpha_{\ell,m}[\varphi^w]$ according to equations~\eqref{eq:Clensh-Curtis-1} are about one order of magnitude smaller than those corresponding to the expansion coefficients $\alpha_{\ell,m}[\varphi]$ according to equations~\eqref{eq:Clensh-Curtis}. Also, should iterative solvers such as GMRES be used for the solution of the QBX Nystr\"om linear systems, the iterations counts are more favorable in the case of the single layer formulation that incorporates the weighted unknown $\varphi^w$ defined in equation~\eqref{eq:weighted_phi}. For these reasons, the use of weighted single layer formulations is preferred.

We report in Table~\ref{compQBX1} numerical results concerning QBX quadrature Nystr\"om discretizations of the single layer formulation of the Helmholtz equation~\eqref{DIRICH} using the weighted unknown~\eqref{eq:weighted_phi}, a number of Chebyshev panels that coincides with the number of corners in the case of closed scatterers and respectively the number of corners plus one in the case of open arcs (that is, one global Chebyshev mesh in the case of the teardrop, boomerang, and strip, and two Chebyshev meshes in the case of the V-shaped scatterer), as well as the Clenshaw-Curtis quadratures~\eqref{eq:Clensh-Curtis-1} to evaluate the QBX expansion coefficients $ \alpha_{\ell,m}[\varphi^w]$ in the expansion~\eqref{eq:qbx}. We focused on the case of real wavenumbers, so that highly accurate reference solutions can be produced by alternative high-order discretizations~\cite{ColtonKress,bruno2012second}; qualitatively similar results are obtained in the cases when the wavenumber is complex. In these numerical experiments we fixed the QBX expansion order parameter $p=8$ and the upsampling factor $\beta=6$, and we present results as the discretization Chebyshev mesh is increasingly refined. We note that in the case of open arcs the convergence of the Nystr\"om QBX discretizations with respect to the size of the mesh appears to be second order, which can be accounted for by the fact that the endpoints of the open arc are the nearest singularities of the smooth extension onto $\Gamma$ of single layer potentials evaluated at nearby exterior points. 
\begin{table}
  \begin{center}
\begin{tabular}{|c|c|c|c|c|c|c|c|c|}
\hline
$n_1$ & \multicolumn{2}{c|} {Teardrop} &  \multicolumn{2}{c|} {Boomerang} & \multicolumn{2}{c|}{Strip} & \multicolumn{2}{c|}{V-shaped Strip}\\
\hline
 & $\varepsilon_\infty$ & e.o.c.  & $\varepsilon_\infty$ & e.o.c. & $\varepsilon_\infty$ &  e.o.c. & $\varepsilon_\infty$ &  e.o.c. \\
\hline
32 & 1.3 $\times$ $10^{-2}$ & & 4.4 $\times$ $10^{-2}$& & 2.6 $\times$ $10^{-4}$ & & 4.7 $\times$ $10^{-4}$ &\\
64 & 6.5 $\times$ $10^{-6}$ & 10.9 & 1.0 $\times$ $10^{-5}$ & 12.01 & 5.9 $\times$ $10^{-5}$ & 2.18 & 1.1 $\times$ $10^{-4}$ & 1.99\\
128 & 1.0 $\times$ $10^{-6}$ & 2.70 & 4.9 $\times$ $10^{-8}$& 7.73 & 1.3 $\times$ $10^{-5}$& 2.10 & 2.9 $\times$ $10^{-5}$ & 1.99\\
256 & 1.6 $\times$ $10^{-7}$ & 2.64 & 2.2 $\times$ $10^{-10}$ & 7.80 & 3.3 $\times$ $10^{-6}$ & 2.07& 7.4 $\times$ $10^{-6}$ & 1.99 \\
\hline
\end{tabular}
\caption{Errors in the near field and estimated orders of convergence obtained from the QBX Nystr\"om quadrature discretization of the single layer formulation for real wavenumber $k=8$ and plane wave normal incidence. We used the weighted unknown~\eqref{eq:weighted_phi} and the Clenshaw-Curtis quadratures~\eqref{eq:Clensh-Curtis-1} in the expansion~\eqref{eq:qbx}, as well as the QBX parameters $p=8$ and $\beta=6$.\label{compQBX1}}
\end{center}
\end{table}

With Table~\ref{compQBXAlp} we investigate the role played by nearby singularities in QBX field expansions. Specifically, we present near field accuracy errors achieved by Nystr\"om QBX discretizations of the single layer formulation of the Helmholtz equation with real wavenumber $k=8$ in the case of teardrop scatterers with very acute apertures  whose parametrizations are given by $x(t)=\left(2\sin{\frac{t}{2}},-\beta\sin{t}\right),\ \beta=\tan{\frac{\alpha\pi}{2}},\ \alpha=1/32,\ 0\leq t\leq 2\pi$. We used in our numerical experiments both one panel (a global Chebyshev mesh---second column) as well as 8 Chebyshev panels dyadically refined near the corner (third column), we fixed the same expansion parameter $p=6$ and the up sampling factor $\beta=4$ for each panel, and we refined the size $n_m,\ 1\leq m\leq M$ of the Chebyshev meshes.  We observe that the QBX discretizations give rise to second order convergent Nystr\"om solvers (with respect to the size $n_m$ of Chebyshev meshes), and we note that Alpert solvers (first column in Table~\ref{compQBXAlp}) do not exhibit a deterioration of the order of convergence upon Chebyshev mesh refinement for very acute corners. Nevertheless, the accuracy levels of the Alpert and QBX Nystr\"om discretizations appear to be quite similar for the same mesh size. On the other hand, following~\cite{klockner2013quadrature}, we present in the fourth column of Table~\ref{compQBXAlp} the convergence properties of QBX Nystr\"om discretizations from a different perspective, that is we considered 64 Chebyshev panels dyadically refined near the corner so that the smallest size of the panels is of the order $10^{-9}$, the size of the Chebyshev meshes on each panel was fixed at $n_m=8$, and we progressively increased the values of the parameters $(p,\beta)$. In this setting the QBX Nystr\"om discretizations exhibit high order convergence with respect to the parameters $(p,\beta)$. We applied in the experiments presented in Table~\ref{compQBXM} the same dyadic panel refinement QBX strategy  in the case of scattering off strips and high order convergence with respect to the parameters $(p,\beta)$ is observed again. However, since the single layer formulation is a BIE of the first kind, the condition numbers and numbers of GMRES iterations required for convergence of the QBX Nystr\"om matrices tend to grow with the number $M$ of panels used for the same global size of the discretization. One possibility to mitigate this issue is to resort to recursively compressed inverse preconditioning (RCIP)~\cite{helsing2008corner}, which will be pursued elsewhere. For the scatterers considered in this paper and the levels of accuracy of CQ methods for the solution of the time domain scattering (e.g. $10^{-8}$), the use of geometrically large panels in connection to QBX appears to be adequate as either the direct or iterative solvers of the linear systems involving QBX Nystr\"om matrices are well behaved.

\begin{table}
  \begin{center}	
\begin{tabular}{|c|c|c|c|c|c|c|c|c|c|c|}
\hline
\multicolumn{3}{|c|} {Alpert} &  \multicolumn{3}{c|} {QBX $M=1$} & \multicolumn{3}{c|}{ QBX $M=8$} & \multicolumn{2}{c|}{ QBX $M=64, n_m=8 $} \\
\hline
 $N$ & $\varepsilon_\infty$ & e.o.c.  & $n_1$ & $\varepsilon_\infty$ & e.o.c. & $n_m$ & $\varepsilon_\infty$ & e.o.c.& $(p,\beta)$ & $\varepsilon_\infty$ \\
\hline
64 & 9.8 $\times$ $10^{-3}$ & & 64 &2.0 $\times$ $10^{-4}$& & 8 & 4.3 $\times$ $10^{-4}$ & & (1,1) & 5.8 $\times$ $10^{-2}$ \\
128 & 3.5 $\times$ $10^{-4}$ & 4.81 & 128 & 4.9 $\times$ $10^{-5}$ & 2.02 & 16 & 1.1 $\times$ $10^{-4}$ & 2.04& (2,2) & 7.7 $\times$ $10^{-3}$ \\
256 & 1.0 $\times$ $10^{-5}$ & 5.01& 256 & 1.1 $\times$ $10^{-5}$ & 2.02 &  32 & 2.6 $\times$ $10^{-5}$ & 2.00& (4,2) & 3.9 $\times$ $10^{-4}$  \\
512 & 5.7 $\times$ $10^{-7}$ & 4.25  & 512 & 2.9 $\times$ $10^{-6}$ & 2.01 & 64 & 6.4 $\times$ $10^{-6}$ & 2.01& (6,4)  & 2.5 $\times$ $10^{-5}$ \\
\hline
\end{tabular}
\caption{Errors in the near field and estimated orders of convergence obtained from the Alpert and QBX Nystr\"om quadrature discretizations of the single layer formulation for real wavenumber $k=8$ and plane wave normal incidence impinging on a teardrop with a very acute aperture. We used the parameters $\sigma=4,\ a=2,\ m=3$ for the Alpert quadrature, the single layer formulation with the weighted unknown~\eqref{eq:weighted_phi}, $M=1$ and $M=8$ panels with QBX parameters $p=6$ and $\beta=4$ and increasingly refined Chebyshev meshes, as well as $M=64$ dyadically refined panels around the corner with Chebyshev meshes of a fixed size $n_m=8$ and various values of the $(p,\beta)$ QBX parameters.\label{compQBXAlp}}
\end{center}
\end{table}

\begin{table}
  \begin{center}
\begin{tabular}{|c|c|c|c|}
\hline
\multicolumn{2}{|c|}{ QBX $M=16, n_m=16$} & \multicolumn{2}{c|}{ QBX $M=64, n_m=8 $} \\
\hline
  $(p,\beta)$ & $\varepsilon_\infty$ & $(p,\beta)$ & $\varepsilon_\infty$ \\
\hline
(1,1) & 1.4 $\times$ $10^{-2}$ &1,1) & 4.9 $\times$ $10^{-2}$ \\
(2,2) & 7.5 $\times$ $10^{-4}$ & (2,2) & 4.9 $\times$ $10^{-3}$ \\
(6,4) & 1.1 $\times$ $10^{-5}$ & (6,4) & 1.3 $\times$ $10^{-4}$ \\
(8,4)  & 6.5 $\times$ $10^{-7}$ & (8,4)  & 5.3 $\times$ $10^{-6}$ \\
\hline
\end{tabular}
\caption{Errors in the near field obtained from the QBX Nystr\"om quadrature discretization of the single layer formulation for real wavenumber $k=8$ and plane wave normal incidence impinging on the strip scatterer. We used the formulation involving the weighted unknown~\eqref{eq:weighted_phi}, 16 and respectively 64 dyadically refined panels around the end points of the strip,  Chebyshev grids of size $n_m=16$ and respectively $n_m=8$ on each panel, Clenshaw-Curtis quadratures~\eqref{eq:Clensh-Curtis-1} in the expansion~\eqref{eq:qbx}, and various QBX parameters $(p,\beta)$.\label{compQBXM}}
\end{center}
\end{table}

\subsubsection{CQ Numerical Experiments}\label{semidiscrete}

The numerical solution of a semi-discretization in time for the example in Table~\ref{compEig} showcases the advantages that can be garnered from using Nystr\"om discretizations based on Alpert and QBX quadratures when dealing with layer potentials featuring the complex wavenumbers neccessary in the CQ solution of time-domain problems. Specifically, we considered the ensemble of wavenumbers $s_\ell = \frac{P( \lambda \zeta^{-\ell}_{N+1}) }{\Delta t}$ where $\zeta_{N_t+1}=e^{2\pi i/(N_t+1)}$ and $P(\zeta)=\frac{1}{2}(\zeta^2-4\zeta+3)$ for $0\leq \ell\leq N$ with $\lambda=\max\{(\Delta t)^{3/N},{\rm{eps}}^{1/2N}\}$, $N_t=1024$ and $\Delta t=\frac{2}{N_t}$, associated with the CQ BDF2 method (here $\rm{eps}$ is the machine precision)~\cite{banjai2009rapid}. In the case when $\Gamma$ is the unit circle, the single layer operators are diagonalizable, that is
\begin{equation}\label{eq:eig}
  V(s_\ell)e^{imx}=\frac{i\pi}{2}J_{|m|}(s_\ell)H_m^{(1)}(s_\ell)e^{imx},\quad m\in\mathbb{Z}.
\end{equation}
We used the relation~\eqref{eq:eig} in the case $m=3$ to test the accuracy of various Nystr\"om discretizations of the single layer BIO $V(s_\ell),\ 0\leq \ell\leq N_t$. Specifically, we considered equispaced discretizations $x^M$ of size $M$ of the unit circle and various Nystr\"om discretizations $V^M(s_\ell)$ end we report the errors on the boundary
\begin{equation*}
  \varepsilon_M=\max_{0\leq \ell\leq N_t}\max_{x\in x^M}\left|V^M(s_\ell)e^{imx}-\frac{i\pi}{2}J_{|m|}(s_\ell)H_m^{(1)}(s_\ell)e^{imx}\right|
\end{equation*}
as the size $M$ of the discretization is increased. We investigated the behavior of two classes of Nystr\"om discretization: the Alpert and QBX quadratures that are applied directly to the kernel of the single layer BIO.  We considered discretizations that can deliver accuracies at the level of ${\rm eps}^{\frac{1}{2}}$ as this is the level achievable by CQ quadratures~\cite{banjai2009rapid}. It can be seen from the results in Table~\ref{compEig} that the Nystr\"om discretizations based on Alpert and QBX quadratures respectively achieve uniformly small errors at all values of $M$.
\begin{table}
  \begin{center}
\begin{tabular}{|c|c|c|}
\hline
$M$ & Alpert $a=6,\ m=10$ & QBX $M=1,\ p=12,\ \beta=4$  \\

\hline
256 & 7.4 $\times$ $10^{-6}$ & 7.1 $\times$ $10^{-7}$\\
512 & 2.7 $\times$ $10^{-7}$ & 1.1 $\times$ $10^{-7}$\\
\hline
\end{tabular}
\caption{Boundary errors in the eigenvalue test for circular geometries and a CQ BDF ensemble of complex wavenumbers.\label{compEig}}
\end{center}
\end{table}

\subsection{Axi-symmetric scatterers in three dimensions}

The Alpert quadratures technology can be extended to the solution of frequency domain BIE in three dimensions in the case when the scatterer $\Gamma$ is axisymmetric~\cite{young2012high,epstein2019high} taking advantage of the Fourier modal representation of the rotationally invariant Green functions~\cite{cohl1999compact,conway2010exact}. Indeed, we assume without loss of generality that the boundary $\Gamma$ of the body of revolution $\Omega$ is obtained from rotating the curve $\gamma:[0,2\pi]\to \mathbb{R}^3$ whose parametrization is given by $\gamma(t)=(r(t),0,z(t))$, that is the surface $\Gamma$ is parametrized in the form $x(t,\theta)=(r(t)\cos{\theta},r(t)\sin{\theta},z(t))$ with $0\leq\theta\leq 2\pi$. In this case, the single layer BIE formulation of the Helmholtz equation is expressed in the form 
\[
\int_\Gamma G_k(x-x')\varphi(x')ds(
x')=f(x),\quad x\in\Gamma
\]
in terms of the unknown functional density $\varphi$ on $\Gamma$, where $ds(x')=r(\tau)|\gamma'(\tau)|d\tau\ d\theta'$ if $x'(\tau,\theta')=(r(\tau)\cos{\theta'},r(\tau)\sin{\theta'},z(\tau))$. Since the Green's function $G_k(x-x')=G_k(r(t),z(t),r(\tau),z(\tau);\theta-\theta')$ is rotationally invariant, we look for the density $\varphi$ in the Fourier form
\[
\varphi^w(\tau,\theta'):=\varphi(\tau,\theta') r(\tau)|\gamma'(\tau)|=\sum_{m=-\infty}^\infty \varphi_m(\tau)e^{im\theta'}
\]
whose Fourier coefficients $\varphi_m$ satisfy the following integral equations along the generating curve $\gamma$:
\begin{equation}\label{eq:Fourier}
2\pi\int_\gamma G_m(t,\tau;k)\varphi_m(\tau)d\tau=f_m(t)
\end{equation}
where $f_m$ are the Fourier coefficients of the right hand side function $f$ and $G_m(t,\tau)$ are the modal coefficients of the rotationally invariant Green's function defined as
\begin{equation}
G_k(r(t),z(t),r(\tau),z(\tau);\theta-\theta')=\sum_{m=-\infty}^\infty G_m(t,\tau;k)e^{im(\theta-\theta')}.
\end{equation}
The modal Fourier coefficients $G_m(t,\tau)$ exhibit logarithmic singularities at $t=\tau$ and therefore the numerical solutions of the integral equations~\eqref{eq:Fourier} is amenable to the Alpert quadrature approach described in the previous sections, in both cases when $\gamma$ is a closed curve or an open arc. In practice the Fourier series are truncated $-M_F\leq m\leq M_F$ and consequently we solve $2M_F+1$ modal integral equations of the type~\eqref{eq:Fourier}. We follow the same strategy described in references~\cite{young2012high,epstein2019high} for the numerical evaluation of the modal coefficients $G_m(t,\tau;k)$. Specifically, we compute the modal coefficients $G_m(t,\tau;k)$ via convolutions of Fourier series corresponding to the the rotationally invariant quantities $e^{ik|x(t,\theta) -x'(\tau,\theta')|}$ and the modal coefficients $G_m(t,\tau;0)$ of the Laplace Green's function. The latter quantities $G_m(t,\tau;0)$ can be computed explcitly~\cite{cohl1999compact} in terms of half-integer degree Legendre functions of the second kind, which, in turn, can be evaluated accurately using recursion relations and Miller's algorithm per the prescriptions in~\cite{young2012high,epstein2019high}. The Fourier coefficients of the quantities $e^{ik|x(t,\theta) -x'(\tau,\theta')|}$, in turn, can be evaluated via FFTs of sufficiently large size which are needed in order to resolve the rapid decay of these quantities due to the possibly large magnitude of the imaginary part of the wavenumbers $k$ that are featured in CQ algorithms. More precisely, we have observed that using the Fourier truncation parameter $M_F=512$ and fine enough Alpert discretizations of the integral equations~\eqref{eq:Fourier}, errors of the order $10^{-7}$ were observed in the CQ solutions of the three dimensional wave equations considered in this paper.

\subsection{Neumann Boundary Conditions}

The CQ methodology can be also applied to the case when the wave equation~\eqref{eq:wave} has Neumann boundary conditions. In this case we would have to solve an ensemble of Helmholtz problems with Neumann boundary conditions 
\begin{equation}
    \begin{cases}
        \Delta u(x) + k^2 u(x) = 0 \quad x\in \Omega \\
        \frac{\partial u(x)}{\partial n} = f(x) \quad x\in \Gamma\\
        u(x) = \sum_{m=-\infty}^\infty u_m H_m^{(1)}(k|x|)e^{im\theta},\quad &|x|>r_0,\quad x=|x|(\cos{\theta},\sin{\theta}).
    \end{cases}
    \label{Nmann}
\end{equation}
for a family of complex wave numbers $k$ such that $\Im{k}\geq 0$.

\subsubsection{Single Layer Formulation}\label{SLform}
We begin with the case when $\Gamma$ is a closed boundary; seeking the solution of the Helmholtz problem~\eqref{Nmann} in the form of the single layer potential
$$u(x) = \int_{\Gamma}G_{k}(x-y)\varphi(y)ds(y),\ x\in\Omega^{+},  $$ 
we obtain the following BIE of the second kind
\begin{equation}\label{eq:DLT}
  -\frac{\varphi}{2}+K^\top(k)[\varphi]=f\quad {\rm on}\ \Gamma,\quad K^\top(k)[\varphi](x):=\int_\Gamma \frac{\partial G_k(x-y)}{\partial n(x)}\varphi(y)ds(y), \ x\in\Gamma,
\end{equation}
and note that~\eqref{eq:DLT} is uniquely solvable whenever $\Im{k}>0$. In order to make~\eqref{eq:DLT} amenable to Nystr\"om discretizations based on Alpert quadratures we use the parametrization $\gamma$ of the boundary curve $\Gamma$ that incorporates the sigmoid transforms and employ the following weighted formulation
\begin{equation}  
    -\frac{1}{2}\varphi^w(t)-\frac{ik}{4}\int_0^{2\pi}H_1^{(1)}(k|\gamma(t)-\gamma(\tau)|)\frac{(\gamma(t)-\gamma(\tau))\cdot (\gamma'(t))^\perp}{|\gamma(t)-\gamma(\tau)|}\varphi^w(\tau)d\tau=f(\gamma(t))|\gamma'(t)|
    \label{Nsingle}
\end{equation} 
where $(\gamma'(t))^\perp=(\gamma'_2(t),-\gamma'_1(t))$ and again here $\varphi^w(t):=|\gamma'(t)|\varphi(\gamma(t))$. On the other hand, the application of QBX quadratures to the evaluation of the BIO $K^\top(k)[\varphi]$ relies on the classical jump formulas for the gradients of single layer potentials
\begin{equation*}
  \nabla_x SL(k)[\varphi](x)=\lim_{x^{\pm}\to x}\nabla_{x^\pm} SL(k)[\varphi](x^\pm)\pm\frac{1}{2}n(x)\varphi(x),\ x\in\Gamma,\ x^{\pm}=x\pm\varepsilon n(x),
  \end{equation*}
which amounts to differentiating term by term the expressions
\begin{equation*}
  SL(k)[\varphi](x^\pm)\approx\sum_{\ell=-p}^{p}\alpha_\ell[\varphi] J_\ell(k|x^\pm-x|)e^{-i\ell\theta^\pm}
\end{equation*}
in the polar coordinate systems around the expansion centers~\cite{klockner2013quadrature}. Furthermore, per the prescriptions in~\cite{klockner2013quadrature}, we average the derivatives of the expansions above with respect to both $x^{+}$ and $x^{-}$ in order to evaluate the BIO expressions in equation~\eqref{eq:DLT}. Specifically, the QBX quadrature takes on the form 
\begin{eqnarray}\label{eq:dltqbx}
  K^\top(k)[\varphi](x)  &\approx&-\frac{1}{2}\sum_{\ell=-p}^{p}k\ \alpha_\ell [\varphi] \left(-J_{\ell+1}(k|x^+-x|+\frac{\ell}{k|x^+-x|}J_{\ell}(k|x^+-x|)\right)e^{-i\ell\theta^+}\nonumber\\
  &+&\frac{1}{2}\sum_{\ell=-p}^{p}k\ \alpha_\ell [\varphi]\left(-J_{\ell+1}(k|x^--x|+\frac{\ell}{k|x^--x|}J_{\ell}(k|x^--x|)\right)e^{-i\ell\theta^-}\nonumber\\
\end{eqnarray}
where the coefficients $\alpha_\ell[\varphi]$ are in turn computed via the formulas~\eqref{eq:Clensh-Curtis}. It is straightforward to modify the formulas above to accommodate for a weighted formulation whose unknown is defined as in formula~\eqref{eq:weighted_phi}. We present in Tables~\ref{compAlpDLT1} and~\ref{compQBXDLT1} numerical results concerning Alpert and QBX discretizations of the second kind formulation~\eqref{eq:DLT} using weighted unknowns for wavenumber $k=8$, using reference solutions produced by the high-order kernel-splitting Nystr\"om discretization in~\cite{turc_corner_N}. It is also possible to apply the dyadic refinement around corners strategy in connection to QBX discretizations in the case of the second kind formulation~\eqref{eq:DLT}, and we obtained qualitatively similar results to those  in Table~\ref{compQBXAlp}. On account of the BIE formulation used being of the second kind, both the condition numbers of the Nystr\"om QBX matrices and their spectral properties in connection with GMRES iteration counts are virtually independent of the number of panels $M$ (and thus their geometrical size). 
\begin{table}
  \begin{center}
\begin{tabular}{|c|c|c|c|c|}
\hline
$N$ & \multicolumn{2}{c|} {Teardrop} &  \multicolumn{2}{c|} {Boomerang} \\
\hline
 & $\varepsilon_\infty$ & e.o.c.  & $\varepsilon_\infty$ & e.o.c.\\
\hline
64 & 2.4 $\times$ $10^{-4}$ & & 6.2 $\times$ $10^{-3}$& \\
128 & 5.7 $\times$ $10^{-6}$ & 5.40 & 1.4 $\times$ $10^{-4}$ & 5.40 \\
256 & 1.9 $\times$ $10^{-7}$ & 4.87 & 4.7 $\times$ $10^{-6}$& 4.91  \\
512 & 7.7 $\times$ $10^{-9}$ & 4.67 & 1.8 $\times$ $10^{-7}$ & 4.71\\
\hline
\end{tabular}
\caption{Errors in the near field and estimated orders of convergence corresponding to the Alpert discretization of the weighted formulation~\eqref{Nsingle} for real wavenumber $k=8$ and plane wave normal incidence. We used the following Alpert quadrature parameters: $\sigma=4$ in the sigmoid transform, and respectively $a=2, m=3$.\label{compAlpDLT1}}
\end{center}
\end{table}

\begin{table}
  \begin{center}
\begin{tabular}{|c|c|c|c|c|}
\hline
$n_1$ & \multicolumn{2}{c|} {Teardrop $M=1$} &  \multicolumn{2}{c|} {Boomerang $M=1$} \\
\hline
 & $\varepsilon_\infty$ & e.o.c.  & $\varepsilon_\infty$ & e.o.c.\\
\hline
32 & 2.5 $\times$ $10^{-2}$ & & 3.3 $\times$ $10^{-1}$& \\
64 & 3.0 $\times$ $10^{-3}$ & 3.07 & 2.6 $\times$ $10^{-2}$ & 3.65 \\
128 & 3.6 $\times$ $10^{-4}$ & 3.07 & 1.6 $\times$ $10^{-3}$& 4.03  \\
256 & 3.1 $\times$ $10^{-5}$ & 3.54 & 8.9 $\times$ $10^{-5}$ & 4.16\\
\hline
\end{tabular}
\caption{Errors in the near field and estimated orders of convergence obtained with the QBX Nystr\"om quadrature discretization of the weighted single layer formulation for real wavenumber $k=8$ and plane wave normal incidence. We used a weighted unknown, one global Chebyshev grid, the QBX expansion~\eqref{eq:dltqbx} and the Clenshaw-Curtis quadratures~\eqref{eq:Clensh-Curtis-1} in the expansion~\eqref{eq:qbx}, as well as the QBX parameters $p=4$ and $\beta=4$.\label{compQBXDLT1}}
\end{center}
\end{table}

\subsubsection{Double Layer Formulation}\label{DLform}
We note that the BIE formulation~\eqref{eq:DLT} is restricted to closed boundaries $\Gamma$ only. An alternative formulation that is applicable to both closed and open boundaries can be derived if we seek the solution of~\eqref{Nmann} in the form of a double layer potential 
$$ u(x) = \int_{\Gamma}\frac{\partial G_{k} (x-y) }{\partial n(y) } \varphi(y) ds(y),\ x\in\Omega^{+}.$$
Then the enforcement of Neumann boundary conditions results in the following BIE
\begin{equation}
    N(k)[\varphi](x) = f(x),\ x\in\Gamma 
    \label{hypN}
\end{equation}
where the operator $N(k):H^{1/2}(\Gamma)\to H^{-1/2}(\Gamma)$ is a BIO which has to be understood in the sense of Hadamard finite parts when $\Gamma$ is a closed curve
\begin{equation}
    N(k)[\varphi](x) = f.p.\int_{\Gamma}\frac{\partial^2 G_{k} (x-y) }{\partial n(x)\partial n(y) } \varphi(y) ds(y),\ x\in\Gamma.
    \label{hyp_op}
\end{equation}
We can recast the operator $N(k)$ in the following form
\begin{align}
    N(k)[\varphi](x) &= \int_{\Gamma}\frac{\partial^2 (G_{k}-G_0) (x-y) }{\partial n(x)\partial n(y) } \varphi(y) ds(y) + N(0)[\varphi](x)\nonumber \\
     N(0)[\varphi](x)    &:= f.p. \int_{\Gamma}\frac{\partial^2 G_{0} (x-y) }{\partial n(x)\partial n(y) } \varphi(y) ds(y)
    \label{Hyp_massage}    
\end{align}
where $G_0$ is the Green's function of the Laplace equation: 
$$G_0(x) = -\frac{1}{2\pi} \log{|x|} . $$
The integral kernel of the operator $N(k)-N(0)$ has a logarithmic singularity~\cite{turc_corner_N}, and as such, the Nystr\"om discretization of the operator $N(k)-N(0)$ is amenable to Alpert quadratures. On the other hand, the Hadamard finite parts integral operator term in the definition of the hypersingular BIO, $N(0)$, can be manipulated using integration by parts techniques starting from the identity
\begin{equation*}
    \frac{\partial^2 G_0(x-y) }{\partial n(x)\partial n(y) } = -\frac{\partial^2 G_0(x-y)}{\partial t(x) \partial t(y)}, 
\end{equation*}
where $t(x)$ is the unit tangent on $\Gamma$ at $x$. In the case when $\Gamma$ is at least $C^3$, the discretization of the Hadamard finite parts integral operator $N(0)$ can be performed using the fact that the kernel
      \[
      N_0(t,s):=\frac{\partial^2}{\partial t\ \partial s}\left(\log(|\gamma(t)-\gamma(s)|^2)-\log \left(4\sin^2\left(\frac{t-s}{2}\right)\right)\right)
      \]
is a bounded function in the variables $(t,s)$, assuming a $\gamma:[0,2\pi]\to \mathbb{R}^2$ parametrization of the closed curve $\Gamma$~\cite{KressH}. Indeed, denoting
      \[
      \widetilde{N}(t,s)=\frac{(\gamma(t)-\gamma(s))\cdot \gamma'(t)\ (\gamma(t)-\gamma(s))\cdot \gamma'(s)}{|\gamma(t)-\gamma(s)|^2}
      \]
we have
      \[
      N_0(t,s)= -\frac{2}{|\gamma(t)-\gamma(s)|^2}(\gamma'(t)\cdot \gamma'(s) -2 \widetilde{N}(t,s))-\frac{1}{2\sin^2\left(\frac{t-s}{2}\right)}
      \]
      with
      \[
      N_0(t,t) =-\frac{( \gamma'(t)\cdot \gamma''(t))^2}{(\gamma'(t)\cdot \gamma'(t))^2}+\frac{1}{2}\frac{\gamma''(t)\cdot \gamma''(t)}{\gamma'(t)\cdot \gamma'(t)}+\frac{1}{3}\frac{\gamma'(t)\cdot \gamma'''(t)}{\gamma'(t)\cdot \gamma'(t)}-\frac{1}{6}.
      \]
Finally, the Hilbert transform operators
      \[
      H[g](t):=\frac{1}{2\pi}\int_0^{2\pi}\cot{(\frac{t-s}{2})}{g}'({ x}(s))ds
      \]
can be discretized using $2\pi$ trigonometric interpolation~\cite{KressH}. In the case when $\Gamma$ is an open smooth arc, the numerical scheme above has to be slightly modified in order to accommodate for the different mapping property of the BIO $N(k)$, that is $N(k):\widetilde{H}^{1/2}(\Gamma)\to H^{-1/2}(\Gamma)$. Thus, the density function $\varphi$ has to vanish at the end points of the arc $\Gamma$, which has to be explicitly accounted for in the numerical scheme. Following the procedure in~\cite{chapko2000numerical}, we assume a parametrization of $\Gamma$ of the form $\gamma:[-1,1]\to\mathbb{R}^2$, followed by the cosine change of variables and an odd periodic extension of the density and the parametrized integral equations on the interval $[0,2\pi]$. Finally, the ensuing discretized integral equations are projected back onto the grid points on the interval $(0,\pi)$.  We present in Table~\ref{compAlpertN} numerical results concerning Alpert quadratures Nystr\"om discretizations of the hypersingular BIE equation~\eqref{hypN} for smooth closed (specifically the kite geometry in~\cite{ColtonKress} and the starfish geometry~\cite{hao2014high}) and open boundaries $\Gamma$ in the case of real wavenumbers; the errors were computed using reference solutions produced via the high-order kernel-splitting Nystr\"om discretization presented in~\cite{turc_corner_N}. Similar orders of convergence are observed in the case when the wavenumbers are complex.

In the case when $\Gamma$ is Lipschitz, the discretization above used in conjunction with sigmoid transforms~\eqref{eq:cov_w} does not result in high-order convergence. The lack of high-order convergence is attributable to the fact that the expressions $N_0(t,s)$ become unbounded for $s=t$. A similar situation occurs when Maue's type integration by parts formulas are employed to recast the operators $N(k)$ in terms of Cauchy Principal Value integral operators. We mention that in the case when $\Gamma$ is closed, the evaluation of hypersingular operators can be entirely bypassed through the use of Calder\'on identities when dealing with the BIO $N(k)$~\cite{anand2012well}. However, the extension of the Calder\'on identities to the case when $\Gamma$ is an open arc requires that the singularities of density functions near the end points be explicitly accounted for~\cite{bruno2012second}. We present in what follows a reformulation of the hypersingular BIE~\eqref{hypN} whose discretization does not require numerical differentiation, and which is amenable to QBX discretizations that are oblivious to whether $\Gamma$ is open or closed.

\begin{table}
  \begin{center}
\begin{tabular}{|c|c|c|c|c|c|c|}
\hline
$N$ & \multicolumn{2}{c|} {Kite} &  \multicolumn{2}{c|} {Starfish} &  \multicolumn{2}{c|} {Strip}\\
\hline
 & $\varepsilon_\infty$ & e.o.c.  & $\varepsilon_\infty$ & e.o.c.& $\varepsilon_\infty$ & e.o.c. \\
\hline
64 & 9.8 $\times$ $10^{-4}$ & & 3.6 $\times$ $10^{-4}$& & 1.5 $\times$ $10^{-5}$ &  \\
128 & 2.7 $\times$ $10^{-6}$ & 5.07 & 1.3 $\times$ $10^{-5}$ & 4.71 & 5.0 $\times$ $10^{-7}$ & 4.94\\
256 & 1.1 $\times$ $10^{-6}$ & 4.70 & 5.4 $\times$ $10^{-7}$& 4.66  & 1.9 $\times$ $10^{-8}$ & 4.71\\
512 & 7.7 $\times$ $10^{-9}$ & 4.63 & 2.2 $\times$ $10^{-8}$ & 4.62 & 7.6 $\times$ $10^{-10}$ & 4.66 \\
\hline
\end{tabular}
\caption{Errors in the near field and estimated orders of convergence obtained when the Alpert discretization is applied to the weighted formulation~\eqref{hypN} in the case of real wavenumber $k=8$ and normal plane wave incidence. We used the following Alpert quadrature parameteres: $\sigma=4$ in the sigmoid transform,and respectively $a=2, m=3$.\label{compAlpertN}}
\end{center}
\end{table}

\subsubsection{Current and Charge (CC) Formulations and QBX Discretizations}\label{CCform}
We start with the Maue's integration by parts recasting of the BIO $N(k)$ in the form
\begin{equation}
    N(k)[\varphi](x) = k^2\int_{\Gamma}G_{k} (x-y) \hat{n}(x)\cdot \hat{n}(y) \varphi(y) ds(y) + \partial_s^x\int_{\Gamma}G_{k} (x-y)\partial_s\varphi(y) ds(y),\ x\in\Gamma
    \label{hyp_opref}
\end{equation}
where $\partial_s$ denote tangential differentiation, where $\varphi\in H^{1/2}(\Gamma)$ when $\Gamma$ is closed and $\varphi\in \widetilde{H}^{1/2}(\Gamma)$ when $\Gamma$ is open (we note that integration by parts is still justified in the latter case since $\varphi$ vanishes at the endpoints of $\Gamma$). We shall use in what follows the $\hat{\cdot}$ notation to denote vector quantities. In what follows, the Current and Charge formulation of~\eqref{hypN} is employed starting from equation~\eqref{hyp_opref}. Our use of the CC formulation is a variant of the formulations pioneered in~\cite{taskinen2006current} for 3D Maxwell equations and subsequently studied in the excellent contributions~\cite{bendali2012extension,epstein2010debye}. We mention that extended BIE formulations that use charges as unknowns were recently analyzed in the context of transmission problems with complex wavenumbers~\cite{helsing2020extended}. First, we introduce the \emph{charge} unknown $\rho$ which is related to the \emph{current} unknown $\varphi$ in the following manner 
\begin{equation}\label{eq:CC}
  \partial_s\varphi=ik\rho.
\end{equation}
Replacing the tangential derivative $\partial_s\varphi$ by the new unknown $\rho$ in the Maue formula~\eqref{hyp_opref} and enforcing the constraint~\eqref{eq:CC} in an integral form, we arrive at the system of BIE
\begin{eqnarray*}
  -ik\int_{\Gamma}G_{k} (x-y) \hat{n}(x)\cdot \hat{n}(y) \varphi(y) ds(y) + \partial_s\int_{\Gamma}G_{k} (x-y)\rho(y)ds(y)&=&\frac{1}{ik}f(x)\\
    \int_{\Gamma}G_{k} (x-y) \partial_s\varphi(y) ds(y)-ik\int_{\Gamma}G_{k} (x-y)\rho(y) ds(y)&=&0.
\end{eqnarray*}
Performing an integration by parts on the integral term above that contains the quantity $\partial_s\varphi$ (this is justified even in the case when $\Gamma$ is an open arc since $\varphi$ vanishes at the open ends) we obtain the following Current and Charge BIE formulation of the Helmholtz problem ~\eqref{Nmann} whose unknowns are the densities $(\varphi,\rho)\in H^{1/2}(\Gamma)\times H^{-1/2}(\Gamma)$ in the case when $\Gamma$ is closed and respectively $(\varphi,\rho)\in \widetilde{H}^{1/2}(\Gamma)\times \widetilde{H}^{-1/2}(\Gamma)$ when $\Gamma$ is an open arc
\begin{eqnarray}\label{eq:CC}
  -ik\int_{\Gamma}G_{k} (x-y) \hat{n}(x)\cdot \hat{n}(y) \varphi(y) ds(y) + \partial_s\int_{\Gamma}G_{k} (x-y)\rho(y)ds(y)&=&\frac{1}{ik}f(x)\nonumber\\
    \nabla\cdot\int_{\Gamma}G_{k} (x-y) \hat{t}(y)\varphi(y) ds(y)-ik\int_{\Gamma}G_{k} (x-y)\rho(y) ds(y)&=&0.
\end{eqnarray}
We note that the CC formulation~\eqref{eq:CC} above features BIOs whose singularities are no longer integrable (i.e they exhibit Cauchy p.v. operators), and as such, they cannot be discretized directly using Alpert quadratures. On the other hand, the QBX methods can handle relatively seamlessly these new operators. In the QBX approach, the derivatives of the layer potentials that feature in the CC formulation~\eqref{eq:CC} are performed by differentiating term by term the expansions around the centers, via the formulas
\begin{equation*}
  \partial_sSL(k)[\rho](x)  \approx-\sum_{\ell=-p}^{p}\ell\alpha_\ell[\rho] \frac{J_{\ell}(k|x^+-x|)}{|x^+-x|}e^{-i\ell\theta^+}
\end{equation*}
and respectively
\begin{eqnarray*}
  \nabla^x \cdot SL(k)[\hat{t}\varphi](x)  &\approx&-\sum_{\ell=-p}^{p}k\hat{\alpha}_\ell[\varphi]\cdot\hat{n}(x)\left(-J_{\ell+1}(k|x^+-x|+\frac{\ell}{k|x^+-x|}J_{\ell}(k|x^+-x|)\right)e^{-i\ell\theta^+}\\
  &-&i\sum_{\ell=-p}^{p}\ell\hat{\alpha}_\ell[\varphi]\cdot \hat{t}(x) \frac{J_{\ell}(k|x^+-x|)}{|x^+-x|}e^{-i\ell\theta^+}
\end{eqnarray*}
where
\begin{equation*}
  \hat{\alpha}_\ell[\varphi]:=\frac{i}{4}\int_\Gamma H_\ell^{(1)}(k|x'-x^+|)e^{i\ell\theta'}\hat{t}(x')\varphi(x')ds(x').
\end{equation*}
Clearly, we can express the QBX expansion coefficients $\alpha_\ell[\rho]=\sum_{m}\alpha_{\ell,m}[\rho]$  and $\hat{\alpha}_\ell[\varphi]=\sum_m \hat{\alpha}_{\ell,m}[\varphi]$ by accounting for the panel contributions. Just like in the Dirichlet case, the incorporation of weighted unknowns in the CC formulation is beneficial with regards to both the condition numbers and the spectral properties vis a vis iterative solvers of the QBX Nystr\"om matrices. Indeed, instead of using the more singular functional density $\rho$, we favor using its weighted version of defined on each panel $\Gamma_m$ as
\begin{equation}\label{eq:weighted_phi_rho1}
  \rho_m^w(t):=\rho(\gamma_m(t))|\gamma_m'(t)|\sqrt{1-t^2}, \quad -1\leq t\leq 1,\ 1\leq m\leq M,
\end{equation}
together with the Crenshaw-Curtis quadratures~\eqref{eq:Clensh-Curtis-1} for the evaluation of the expansion coefficients $\alpha_{\ell,m}[\rho^w]$. Given that the functional density $\varphi$ is more regular than $\rho$, we do not choose to use a weighted version of it in the case when $\Gamma$ is closed, and thus we use the Crenshaw-Curtis quadratures~\eqref{eq:Clensh-Curtis} to evaluate the coefficients $\hat{\alpha}_{\ell,m}[\varphi]$. 

In the case when $\Gamma$ is an open arc, we can incorporate the explicit asymptotic behavior of the more regular density $\varphi$ in the neighborhood of the end points into our numerical scheme. We start with the simpler case when $\Gamma$ is a smooth open arc with the usual $\gamma:[-1,1]\to\mathbb{R}^2$, and we employ regularized unknowns that take into account the singular behavior of the functional densities $(\varphi,\rho)$ in the BIE~\eqref{eq:CC} as proposed in~\cite{bruno2012second}. Specifically, we use the asymptotic behavior $\varphi\sim \mathcal{O}(\sqrt{d})$ and $\rho\sim \mathcal{O}(1/\sqrt{d})$ where $d$ denotes the distance to the open ends, to define weighted densities that are globally regular
\begin{equation}\label{eq:weighted_phi_rho}
  \varphi(\gamma(t))|\gamma'(t)|=\varphi^w(t)\sqrt{1-t^2},\ \rho^w(t):=\rho(\gamma(t))|\gamma'(t)|\sqrt{1-t^2}, \quad -1\leq t\leq 1.
\end{equation}
However, the weighted approach~\eqref{eq:weighted_phi_rho} has to be slightly modified if we represent the open arc $\Gamma$ via panels, that is $\Gamma=\bigcup_{m=1}^M \Gamma_m$. In this case we make a distinction between the treatment of densities defined on panels that contain the end points of $\Gamma$---we assume without loss of generality that those are $\Gamma_1$ and respectively $\Gamma_M$, and densities defined on panels that do not contain the open ends of $\Gamma$. The latter case is straightforward, as we assume $[-1,1]$ parametrizations of those panels equipped with Chebyshev meshes, and we use the weighted version $\rho_m^w, 2\leq m\leq M-1$ defined in equation~\eqref{eq:weighted_phi_rho1}, and the unweighted density $\varphi$. In the former case, since the densities defined on the panels $\Gamma_1$ and $\Gamma_M$ exhibit the asymptotic behavior $\varphi\sim \mathcal{O}(\sqrt{d})$ and $\rho\sim \mathcal{O}(1/\sqrt{d})$ where $d$ denotes the distance to the open ends of $\Gamma$ (and implicitly of $\Gamma_1$ and $\Gamma_M$), we consider in turn $[0,1]$ parametrizations of the open end panels $\Gamma_1$ and $\Gamma_M$ such that the open ends correspond in the parameter space to $t=0$. We then use on $\Gamma_1$ and $\Gamma_M$ Chebyshev quadrature points $v_n$ defined as
\begin{equation*}
  v_n:=\frac{1}{2}(1 + \cos(\vartheta_n)),\quad \vartheta_n:=\frac{(2n-1)\pi}{2N},\quad n=1,\ldots,N,
\end{equation*}
weighted densities $\varphi_m^w$ and $\rho_m^w$ for $m\in\{1,M\}$ that incorporate explicitly the asymptotic behavior of the densities $(\varphi,\rho)$ around the open end points, and the following modified Fej\'er quadratures to evaluate the QBX expansion coefficients
\begin{equation}\label{eq:fej1}
  \int_0^1 \frac{f(v)}{\sqrt{v}}dv\approx\sum_{n=1}^{N}{}^{'}\omega_n^{(1)} f(v_n),\ \omega_n^{(1)}:=\frac{2}{N}\sum_{q=1}^{N}\frac{4}{1-4(q-1)^2}\cos((q-1)\vartheta_n),\quad n=1,\ldots,N.
\end{equation}
and respectively
\begin{equation}\label{eq:fej2}
  \int_0^1 f(v)\sqrt{v}\ dv\approx\sum_{n=1}^{N}{}^{'}\omega_n^{(2)} f(v_n),\ \omega_n^{(2)}:=\frac{2}{N}\sum_{q=1}^{N}\left(\frac{1}{1-4(q-1)^2}+\frac{3}{9-4(q-1)^2}\right)\cos((q-1)\vartheta_n).
\end{equation}
where the prime notation denotes that the first term is halved. We note that the panel approach presented above can also handle the case when the open arc $\Gamma$ is piecewise smooth. 

We present in Table~\ref{compQBXCS1} numerical results concerning QBX discretizations of the CC formulations~\eqref{eq:CC} in the case when the wavenumber $k$ is real; reference solutions were obtained via (a) the high-order kernel-splitting Nystr\"om methods in~\cite{turc_corner_N,bruno2012second} in the case of closed Lipschitz curves and smooth arcs and (b) refined QBX discretizations of the CC formulations with the modified Fej\'er quadratures~\eqref{eq:fej1} and~\eqref{eq:fej2}. Similar orders of convergence were observed in the case when the wavenumbers are complex. We used one panel $M=1$ for the first three scatterers and two panels $M=2$ for the V-shaped arc, a fixed QBX expansion parameter $p=4$ and up-sampling parameter $\beta=4$ and we present errors in the near field together with estimated orders of convergence as the size of the Chebyshev meshes is refined. We have observed in practice that using larger numbers of panels (of smaller size) leads to QBX matrices with larger condition numbers that also require larger GMRES iterations counts for the solution of ensuing linear systems. These findings can be attributed to the fact that the CC formulation is a formulation of the first kind; the condition numbers and GMRES iteration counts do not vary with the size of the panels on the case of the QBX Nystr\"om matrices corresponding to the second kind formulation~\eqref{eq:DLT}. Nevertheless, the CC formulations~\eqref{eq:CC} are advantageous because they are universally applicable to all types of boundaries $\Gamma$, including piecewise smooth closed boundaries as well as open arcs. In addition, we mention that analytical preconditioners can be employed in order to speed up the iterative convergence of the discrete linear systems corresponding to Nystr\"om discretizations of the first kind CC formulations. Indeed, assuming that the Nystr\"om CC system is expressed in the block matrix form
\begin{equation*}
  \begin{bmatrix}\mathcal{C}\mathcal{C}_{11}^n & \mathcal{C}\mathcal{C}_{12}^n\\ \mathcal{C}\mathcal{C}_{21}^n & \mathcal{C}\mathcal{C}_{22}^n\end{bmatrix}\begin{bmatrix}\varphi_{n}\\ \rho_n\end{bmatrix}=\frac{1}{ik}\begin{bmatrix}f_n\\ 0\end{bmatrix}
\end{equation*}
we use the following preconditioned CC formulation
\begin{equation}\label{eq:prec}
  \begin{bmatrix}\mathcal{C}\mathcal{C}_{22}^n & -\mathcal{C}\mathcal{C}_{12}^n\\ -\mathcal{C}\mathcal{C}_{21}^n & \mathcal{C}\mathcal{C}_{11}^n\end{bmatrix}\begin{bmatrix}\mathcal{C}\mathcal{C}_{11}^n & \mathcal{C}\mathcal{C}_{12}^n\\ \mathcal{C}\mathcal{C}_{21}^n & \mathcal{C}\mathcal{C}_{22}^n\end{bmatrix}\begin{bmatrix}\varphi_{n}\\ \rho_n\end{bmatrix}=\frac{1}{ik}\begin{bmatrix}\mathcal{C}\mathcal{C}_{22}^n & -\mathcal{C}\mathcal{C}_{12}^n\\ -\mathcal{C}\mathcal{C}_{21}^n & \mathcal{C}\mathcal{C}_{11}^n\end{bmatrix}\begin{bmatrix}f_n\\ 0\end{bmatrix}.
\end{equation}
The preconditioned CC formulation~\eqref{eq:prec} can be justified heuristically on account of the pseudodifferential calculus, whose details we leave out for now. We have observed in practice that this preconditioned formulation appears to be of the second kind. 

\begin{table}
  \begin{center}
\begin{tabular}{|c|c|c|c|c|c|c|c|c|}
\hline
$n_1$ & \multicolumn{2}{c|} {Teardrop $M=1$} &  \multicolumn{2}{c|} {Boomerang $M=1$} & \multicolumn{2}{c|}{Strip $M=1$}& \multicolumn{2}{c|}{V-shaped Strip $M=2$}\\
\hline
 & $\varepsilon_\infty$ & e.o.c.  & $\varepsilon_\infty$ & e.o.c. & $\varepsilon_\infty$ &  e.o.c. & $\varepsilon_\infty$ &  e.o.c.\\
\hline
32 & 2.1 $\times$ $10^{-1}$ & & 8.4 $\times$ $10^{-2}$& & 1.2 $\times$ $10^{-2}$ & & 8.5 $\times$ $10^{-3}$ & \\
64 & 2.6 $\times$ $10^{-2}$ & 3.03 & 7.6 $\times$ $10^{-3}$ & 3.46 & 3.4 $\times$ $10^{-3}$ & 1.89 & 2.2 $\times$ $10^{-3}$ & 1.95 \\
128 & 2.5 $\times$ $10^{-3}$ & 3.32 & 5.1 $\times$ $10^{-4}$& 3.88 & 9.1 $\times$ $10^{-4}$& 1.89 & 5.4 $\times$ $10^{-4}$ & 2.01\\
256 & 2.4 $\times$ $10^{-4}$ & 3.42 & 6.1 $\times$ $10^{-5}$ & 3.08 & 2.4 $\times$ $10^{-4}$ & 1.89 & 1.3 $\times$ $10^{-4}$ & 1.99\\
\hline
\end{tabular}
\caption{Errors in the near field and estimated orders of convergence obtained when the QBX Nystr\"om discretization is applied to the CC formulation~\eqref{eq:CC} with real wavenumber $k=8$ and normal plane wave incidence. We used the QBX parameters $p=4$ and $\beta=4$.\label{compQBXCS1}}
\end{center}
\end{table}

\section{Time-Domain Experiments}\label{Numerical_results}
We present in this section various numerical experiments concerning solutions of the wave equation in both 2D and 3D using BIE based CQ methods. 
\subsection{2D Convergence Tests}\label{convergencetests}
\begin{figure}
\centering
\includegraphics[height=60mm]{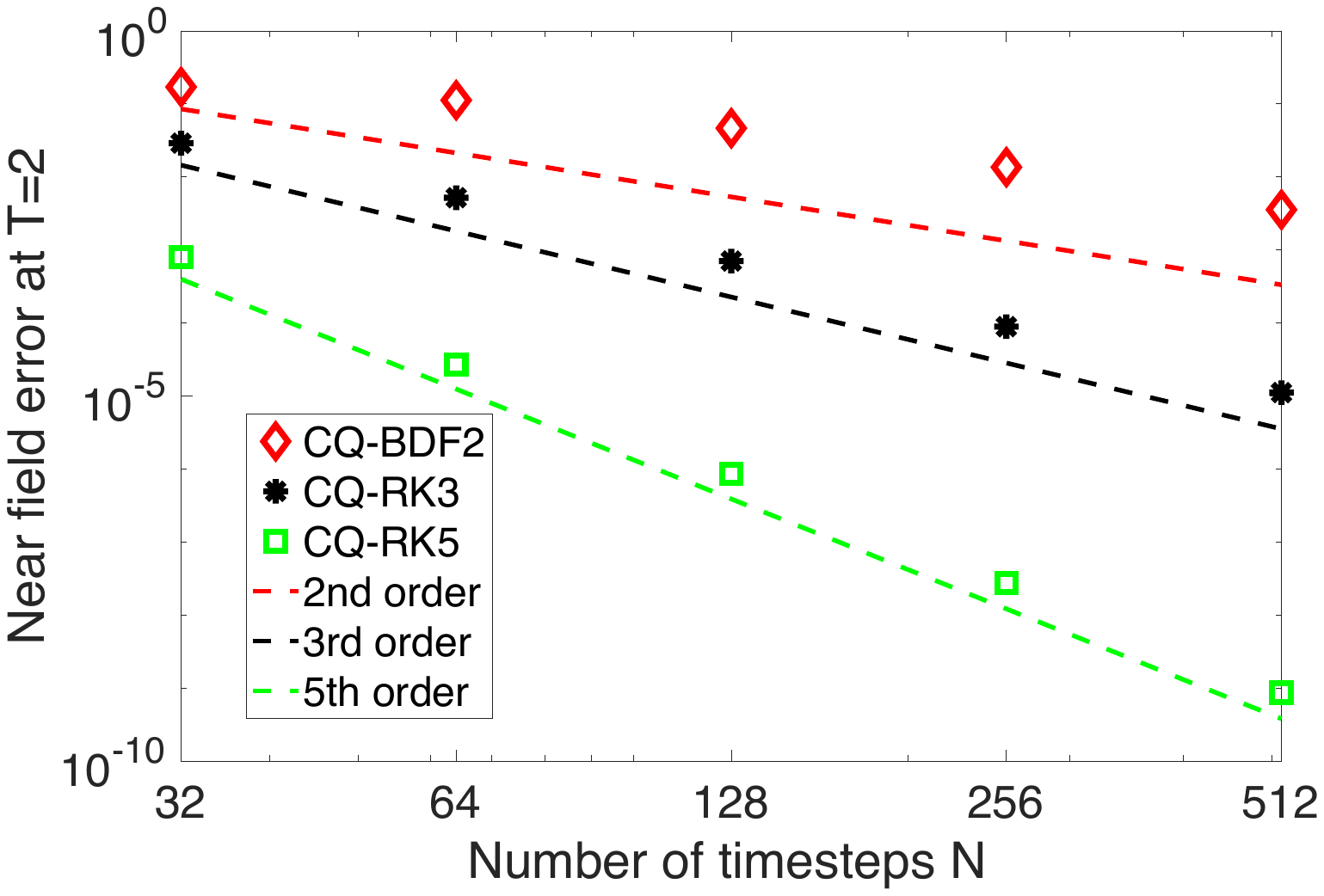}\includegraphics[height=60mm]{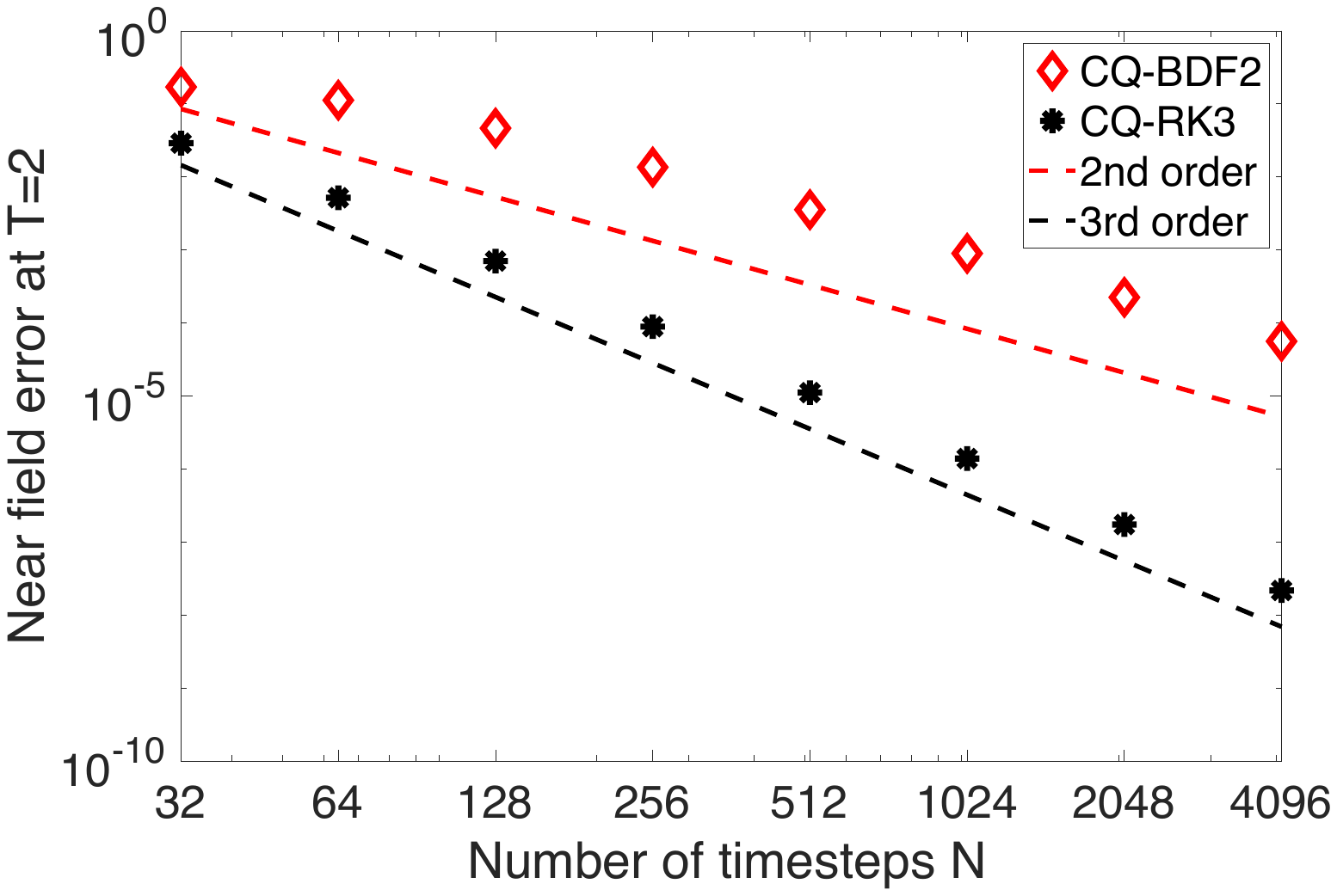}
\caption{Orders of convergence of the near field CQ solutions of the wave equation in the exterior of the teardrop domain at time $T=2$ with Dirichlet boundary conditions. We used Alpert quadratures with the same parameters $\sigma=4$, $a=2$ and $m=3$ for the solution of the ensemble of weighted single layer BIE formulations~\eqref{eq:SLw} in the Laplace domain using $N=256$ discretization points on the boundary. The reference solution was produced using the CQ-RK5 discretization with $2048$ time steps and $N=512$ boundary discretization points.}
\label{fig:CQ_D_teardrop}
\end{figure}

\begin{figure}
\centering
\includegraphics[height=60mm]{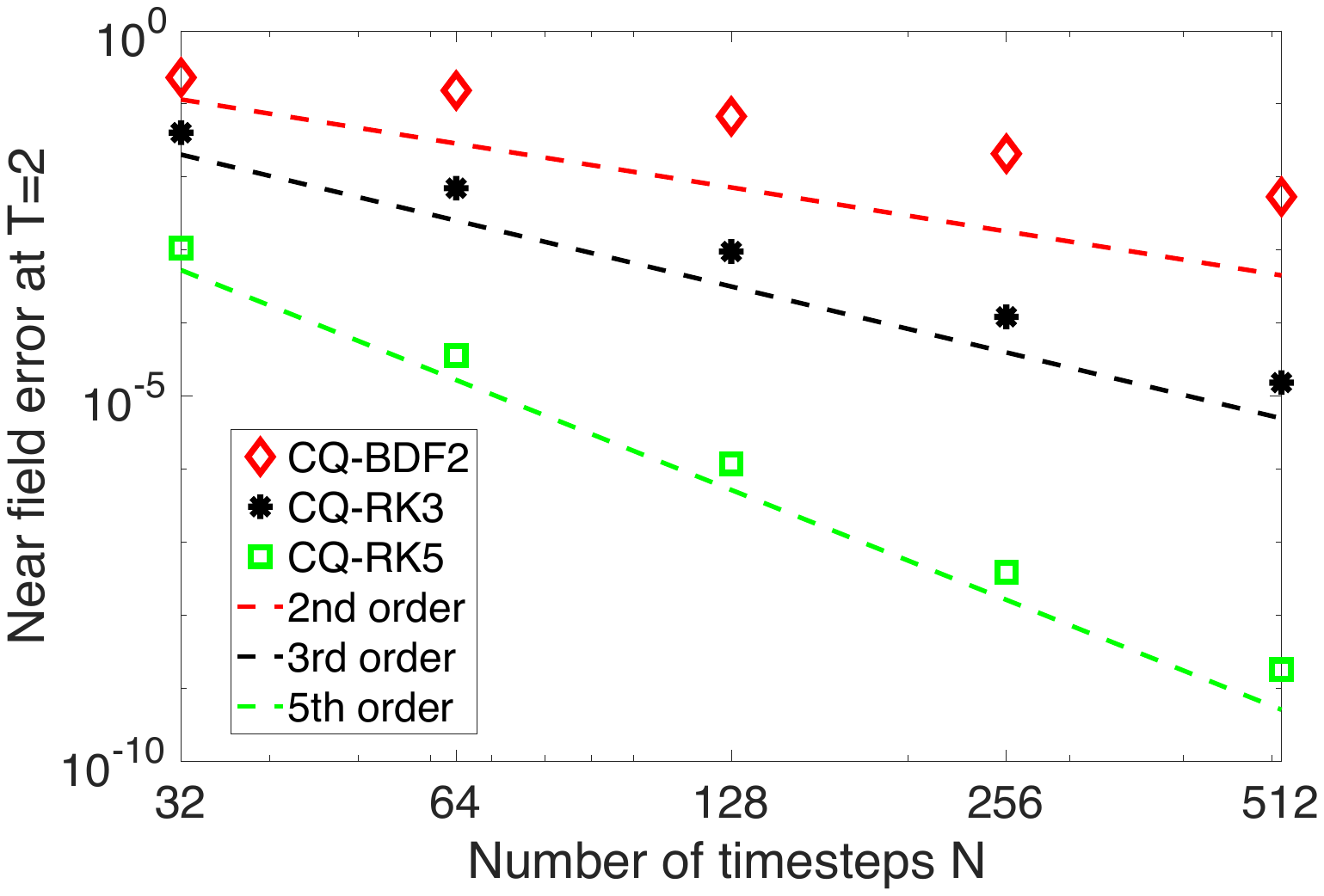}\includegraphics[height=60mm]{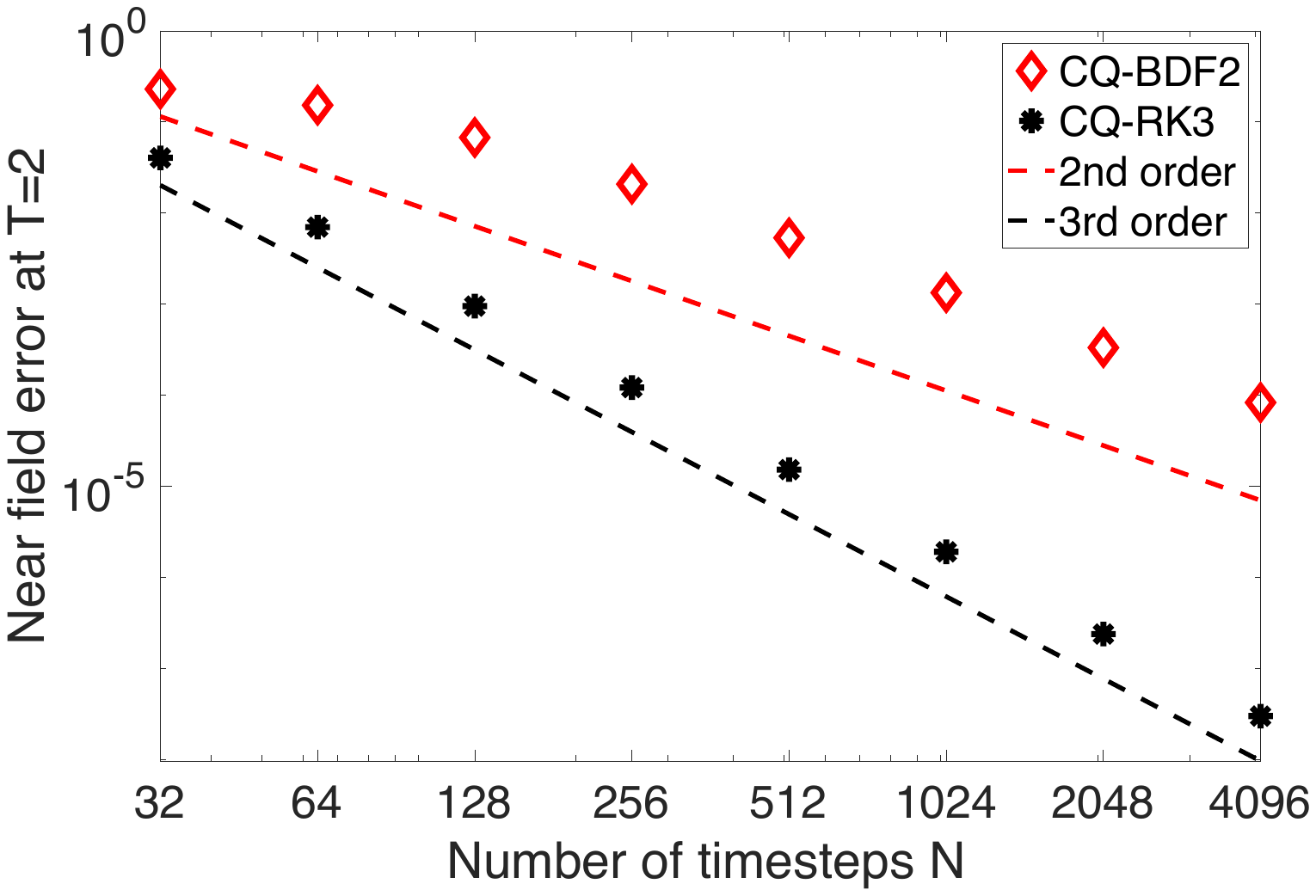}
\caption{Orders of convergence of the near field CQ solutions of the wave equation in the exterior of the boomerang domain at time $T=2$ with Dirichlet boundary conditions. We used Alpert quadratures with the same parameters $\sigma=4$, $a=2$ and $m=3$ for the solution of the ensemble of weighted single layer BIE formulations~\eqref{eq:SLw} in the Laplace domain using $N=256$ discretization points on the boundary. The reference solution was produced using the CQ-RK5 discretization with $2048$ time steps and $N=512$ boundary discretization points.}
\label{fig:CQ_D_boomerang}
\end{figure}

\begin{figure}
\centering
\includegraphics[height=60mm]{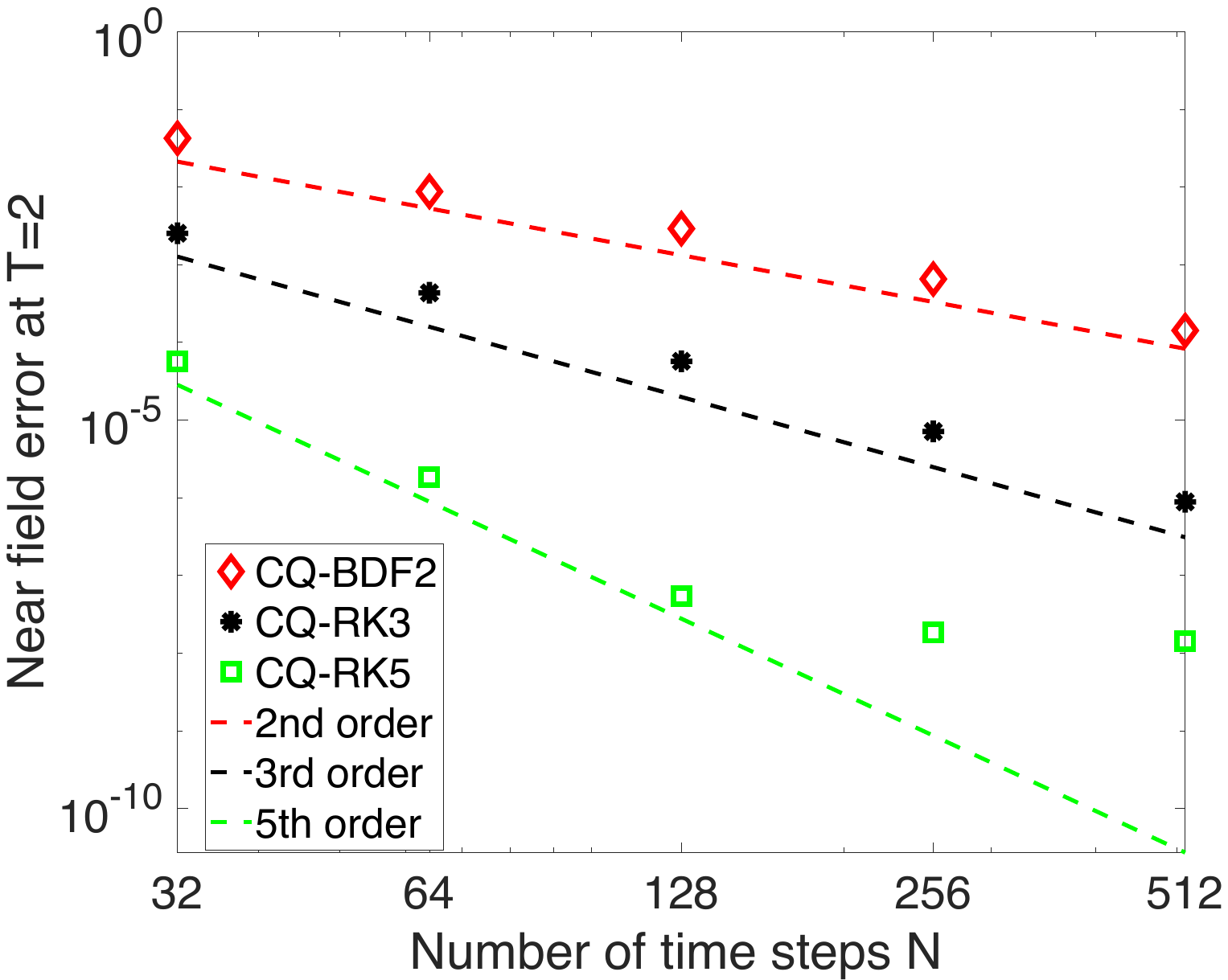}\includegraphics[height=60mm]{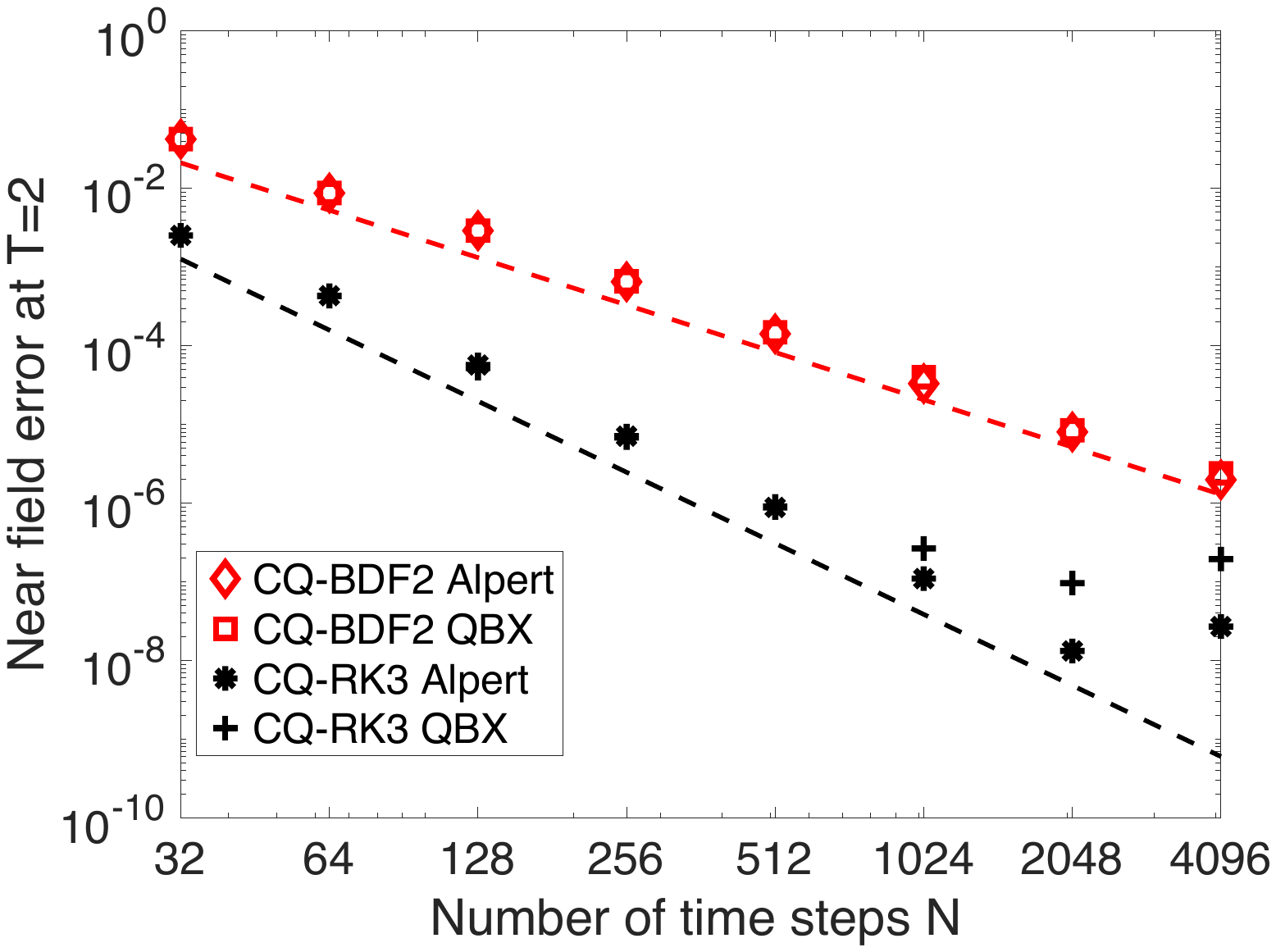}
\caption{Orders of convergence of the near field CQ solutions of the wave equation in the exterior of the strip domain at time $T=2$ with Dirichlet boundary conditions. We used (i) Alpert quadratures with the same parameters $\sigma=4$, $a=2$ and $m=3$ and (ii) Chebyshev based QBX discretizations~\eqref{eq:Clensh-Curtis-1} with parameters $M=1$ (that is one global Chebyshev mesh) $p=8$ and $\beta=4$ for the solution of the ensemble of BIE formulations~\eqref{eq:SLw} in the Laplace domain using $N=256$ discretization points on the boundary for both quadratures. The reference solution was produced using the CQ-RK5 discretization with $2048$ time steps and $N=512$ boundary discretization points.}
\label{fig:CQ_D_arc}
\end{figure}

We consider the time-domain incident plane wave
\begin{equation*}
    g(x,t) = \cos(5t-x \cdot (1,0) )\exp( -1.5(5t-x\cdot (1,0) -5)^2 ),
\end{equation*}
and present results concerning the convergence of CQ numerical schemes that evaluate the scattered fields at $T=2$ at $512$ points placed on circles located at distance $1$ from the scatterers. The incident field considered is smooth enough (in the causality sense, that is the field and its higher order derivatives vanish on the scatterers for negative times all the way to time $t=0$) so that the predicted theoretical rates of convergence rates of CQ methods~\cite{banjai2010multistep,banjai2011runge} are validated in our numerical experiments. The smoothness of the incident fields also manifests itself in the sparse frequency content of the boundary data in the Laplace domain~\cite{banjai2009rapid}, e.g. most boundary data defined in equations~\eqref{cq23} and respectively~\eqref{eq:Gl} are small in the infinity norm. As a result, we follow the common thresholding practice in CQ and we solve only those Laplace domain Helmholtz equations whose boundary data $\hat{g}_\ell$ in equations~\eqref{cq23} is larger than $10^{-10}$ when measured in the infinity norm. In all the time domain CQ experiments based on QBX quadratures for the discretization of the Laplace domain modified Helmholtz problems we considered representations of the scatterers in terms of geometrically large Chebyshev panels, whose number depended only on the number of corners or open ends present on the scatterers' boundaries. We show results in Figure~\ref{fig:CQ_D_teardrop}--\ref{fig:CQ_D_arc} convergence plots of the BDF2 and RK3/RK5 CQ methods in the case of Dirichlet boundary conditions for Lipschitz closed scatterers as well as open scatterers (strips). The Laplace domain problems that enter the CQ methods were solved using the weighted single layer formulation whose associated BIOs were discretized using both Alpert quadratures and QBX (based on the weighted unknown defined in equation~\eqref{eq:weighted_phi} and the Fej\'er quadratures~\eqref{eq:Clensh-Curtis-1}) of appropriately high order. As it can be seen from the results presented in Figure~\ref{fig:CQ_D_teardrop}--\ref{fig:CQ_D_arc} the expected second, third, and respectively fifth orders of convergence were observed for the BDF2, RK3 and respectively RK5 CQ methods. We mention that qualitatively similar results can be observed for smooth scatterers; we simply chose to report time domain results for non-smooth scatterers in this paper as the behavior of CQ solvers for such problems was not discussed in great detail in the delta-bem literature~\cite{sayas2013retarded}. 

\begin{figure}
\centering
\includegraphics[height=60mm]{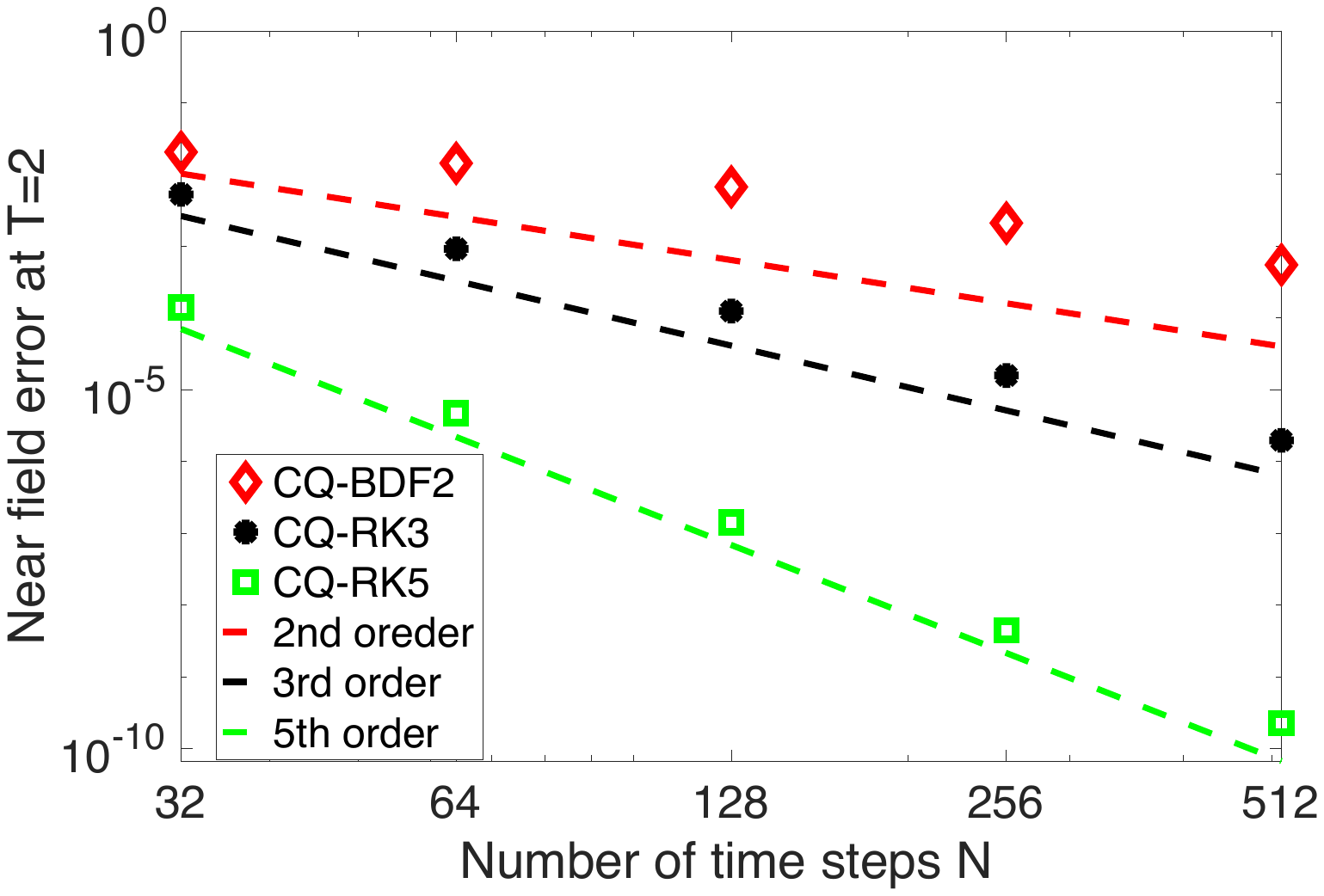}\includegraphics[height=60mm]{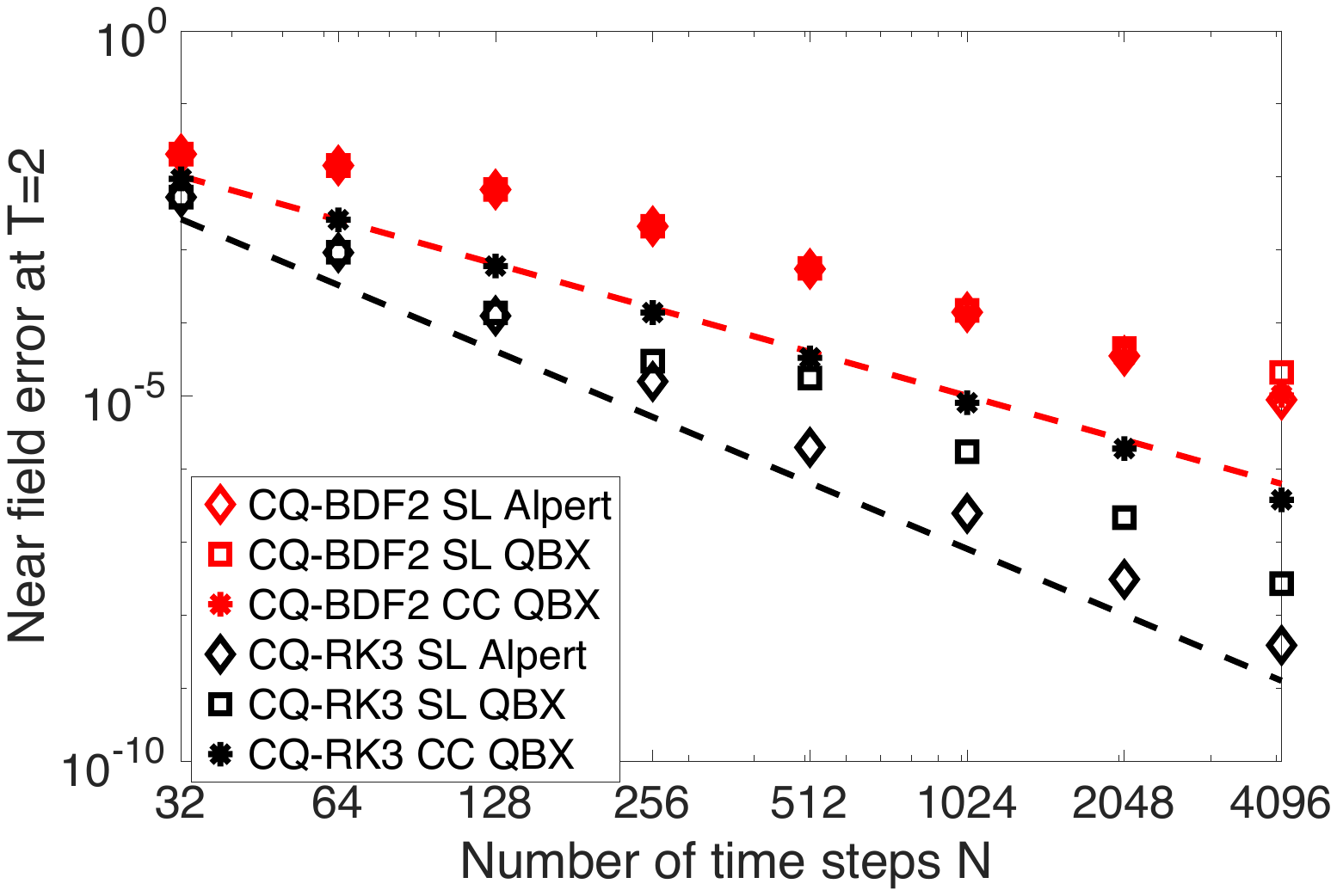}
\caption{Orders of convergence of the near field CQ solutions of the wave equation in the exterior of the teardrop domain at time $T=2$ with Neumann boundary conditions. We used (i) Alpert quadratures with the same parameters $\sigma=4$, $a=2$ and $m=3$ for the solution of the ensemble of second kind BIE formulations~\eqref{Nsingle} in the Laplace domain using $N=256$ discretization points on the boundary and (ii) Chebyshev based QBX discretizations with parameters $M=1,\ n_1=256$, $p=8$, and $\beta=4$ for the solution of the ensemble of CC BIE formulations~\eqref{eq:CC}. The reference solution was produced using the CQ-RK5 discretization with $2048$ time steps and $N=512$ boundary discretization points.}
\label{fig:CQ_N_teardrop}
\end{figure}

\begin{figure}
\centering
\includegraphics[height=60mm]{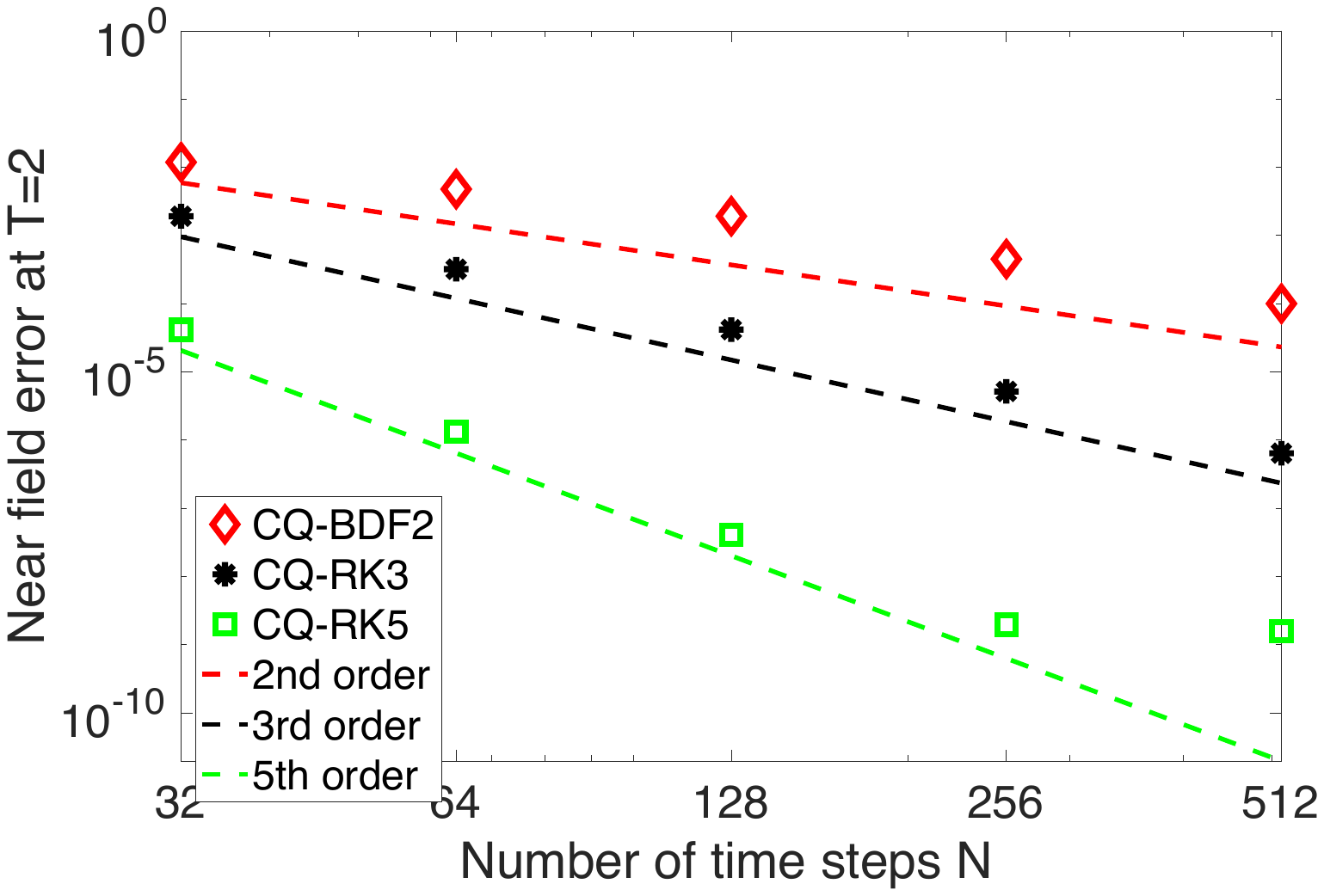}\includegraphics[height=60mm]{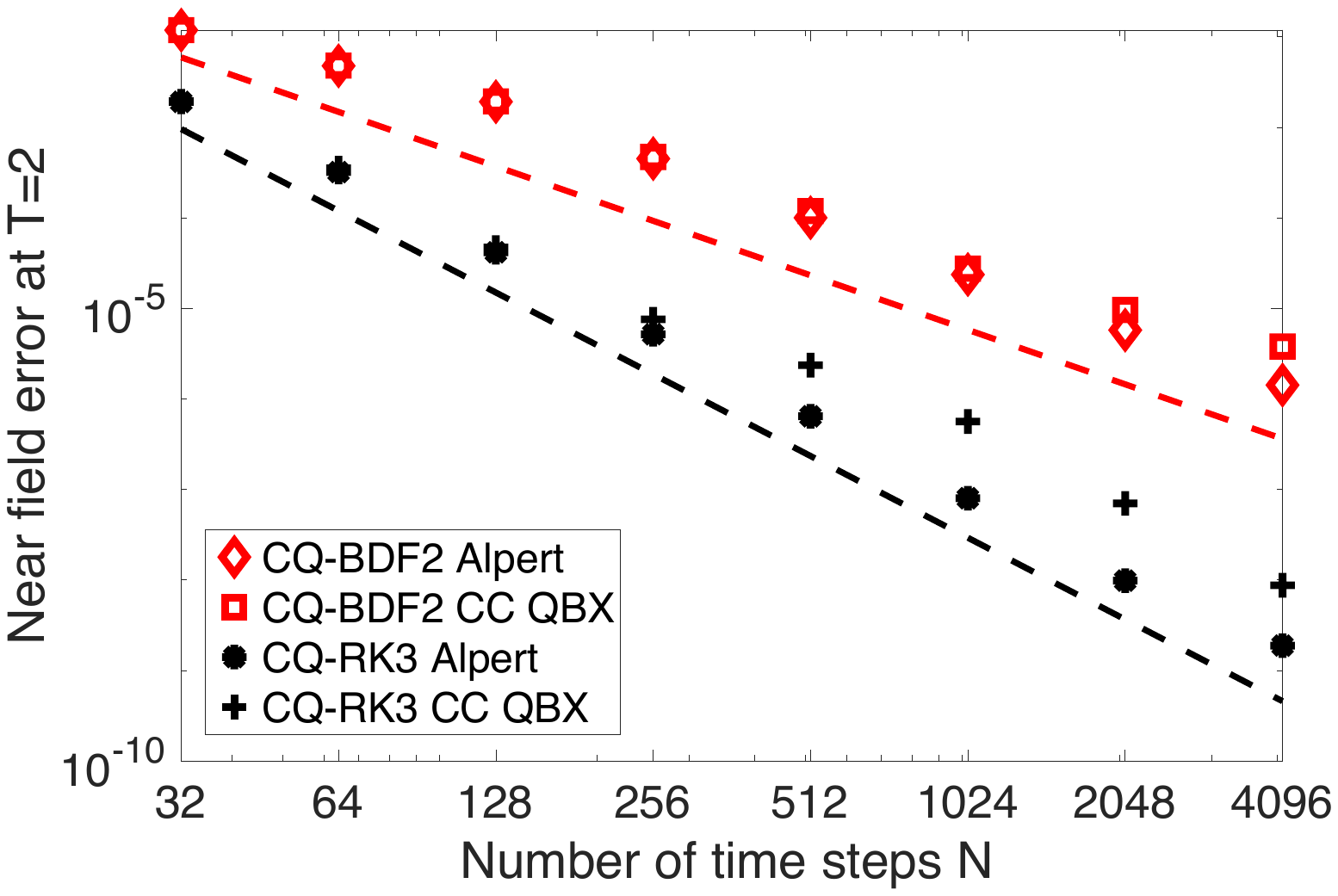}
\caption{Orders of convergence of the near field CQ solutions of the wave equation in the exterior of the strip at time $T=2$ with Neumann boundary conditions. We used (i) Alpert quadratures with the same parameters $\sigma=4$, $a=2$ and $m=3$ for the solution of the ensemble of the hypersingular BIE formulations~\eqref{hypN}in the Laplace domain using $N=256$ discretization points on the boundary and(ii) Chebyshev based QBX discretizations with parameters $M=1,\ n_1=256$, $p=8$ and $\beta=4$ for the solution of the ensemble of CC BIE formulations~\eqref{eq:CC}. The reference solution was produced using the CQ-RK5 discretization with $2048$ time steps and $N=512$ boundary discretization points.}
\label{fig:CQ_N_arc}
\end{figure}

\begin{figure}
\centering
\includegraphics[height=60mm]{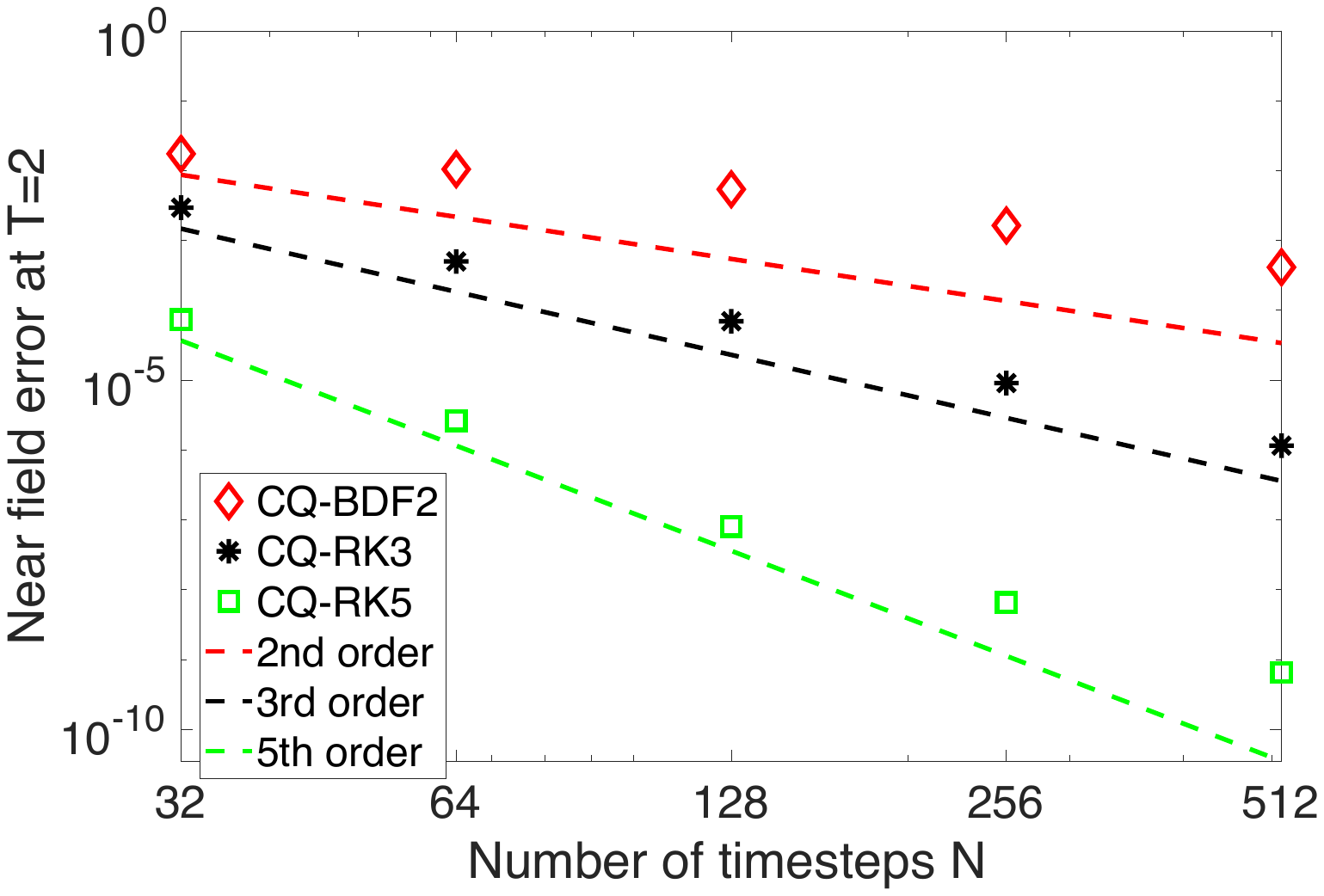}\includegraphics[height=60mm]{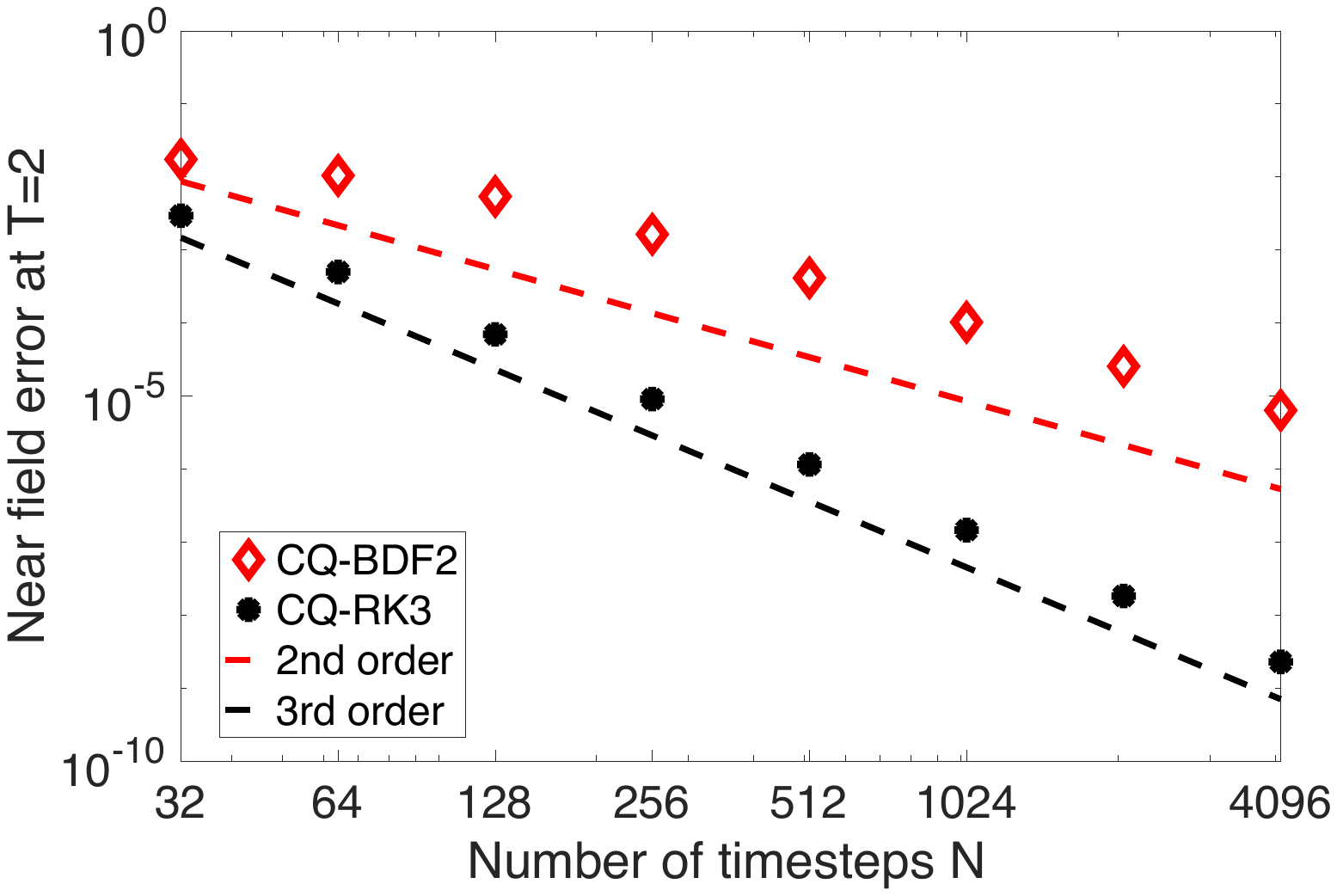}\\
\includegraphics[height=60mm]{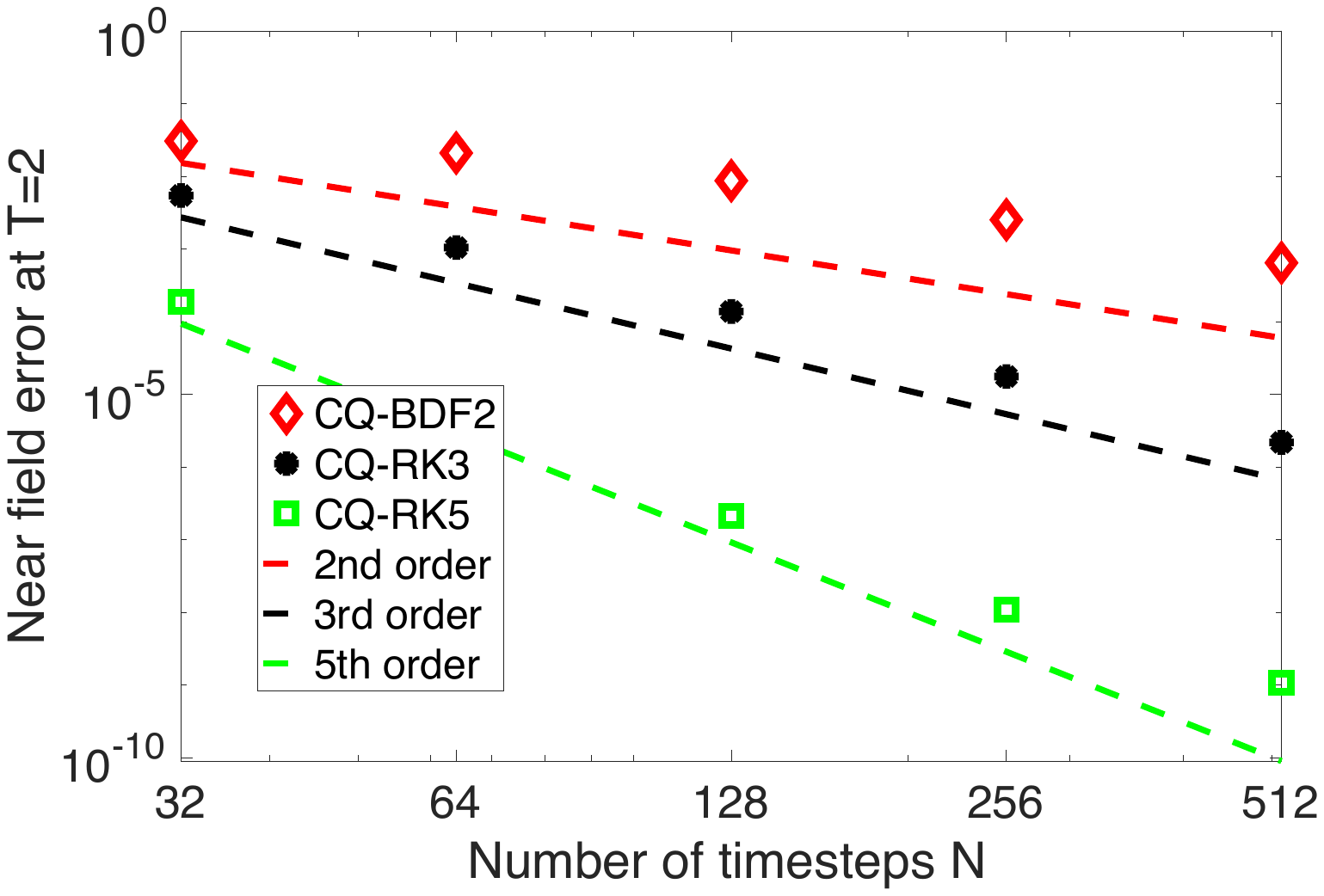}\includegraphics[height=60mm]{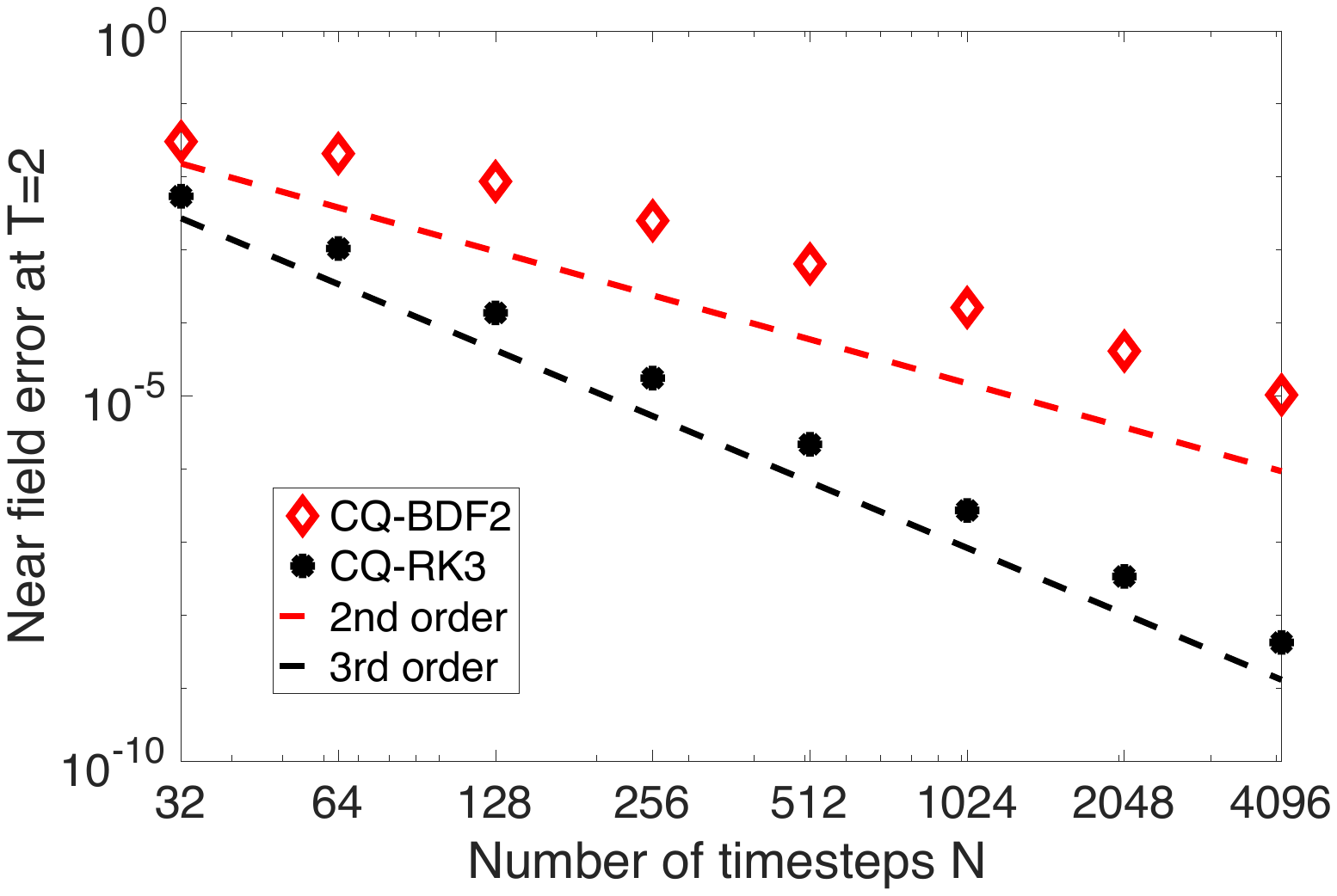}
\caption{Orders of convergence of the near field CQ solutions of the wave equation in the exterior of one V-shaped scatterer (top panel) and two V-shaped scatterers arranged in a chevron pattern (bottom panel) at time $T=2$ with Neumann boundary conditions. We used the  Chebyshev based QBX discretizations with parameters $M=2,\ n_1=n_2=128$, $p=8$ and $\beta=4$ that incorporated the modified Fej\'er quadratures~\eqref{eq:fej1} and~\eqref{eq:fej2} on each V-shaped scatterer for the solution of the ensemble of CC BIE formulations~\eqref{eq:CC}. The reference solution was produced using the CQ-RK5 discretization with $2048$ time steps and $N=512$ boundary discretization points per scatterer.}
\label{fig:CQ_N_vshape}
\end{figure}

\begin{figure}
\centering
\includegraphics[height=60mm]{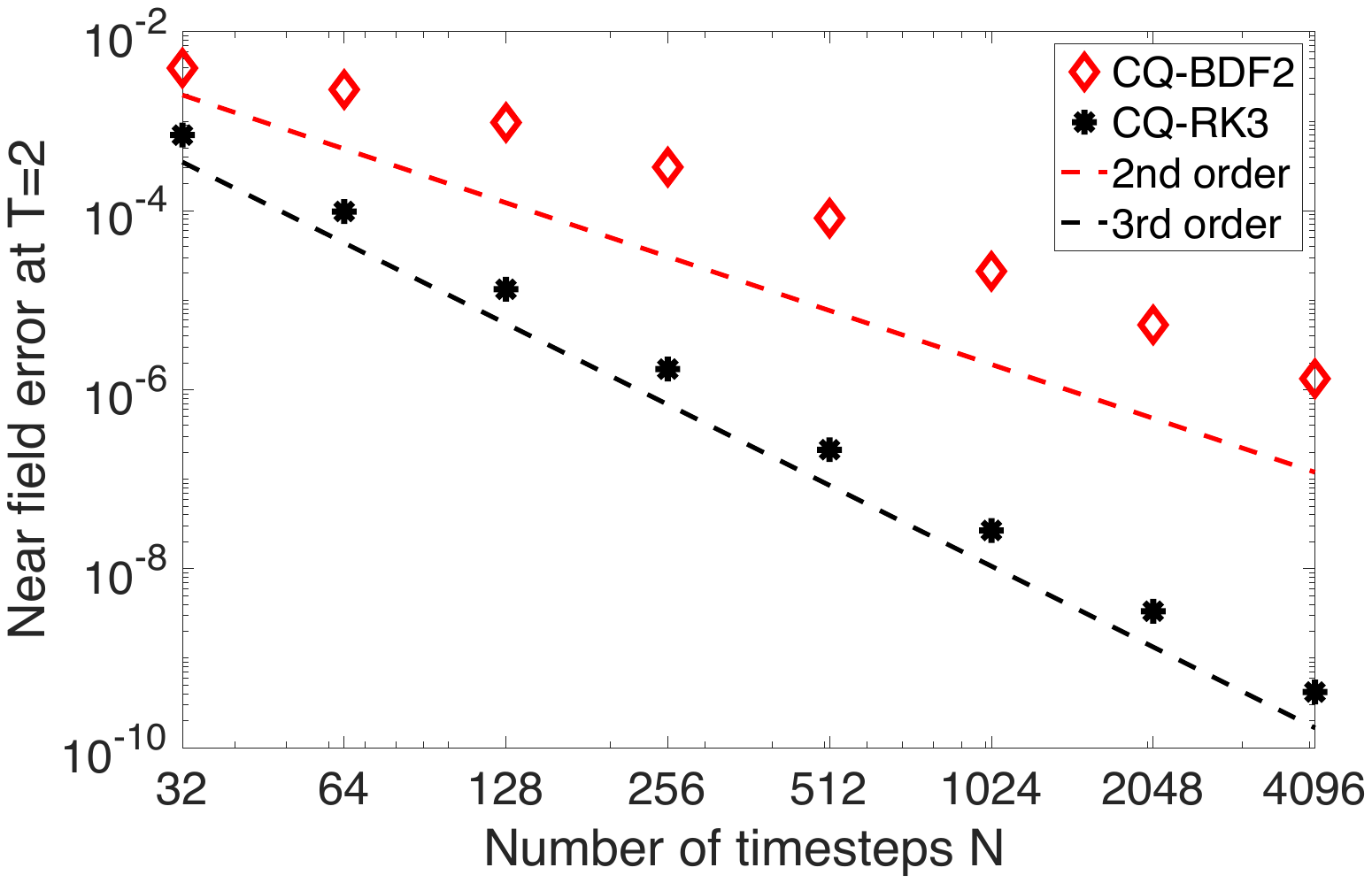}\includegraphics[height=60mm]{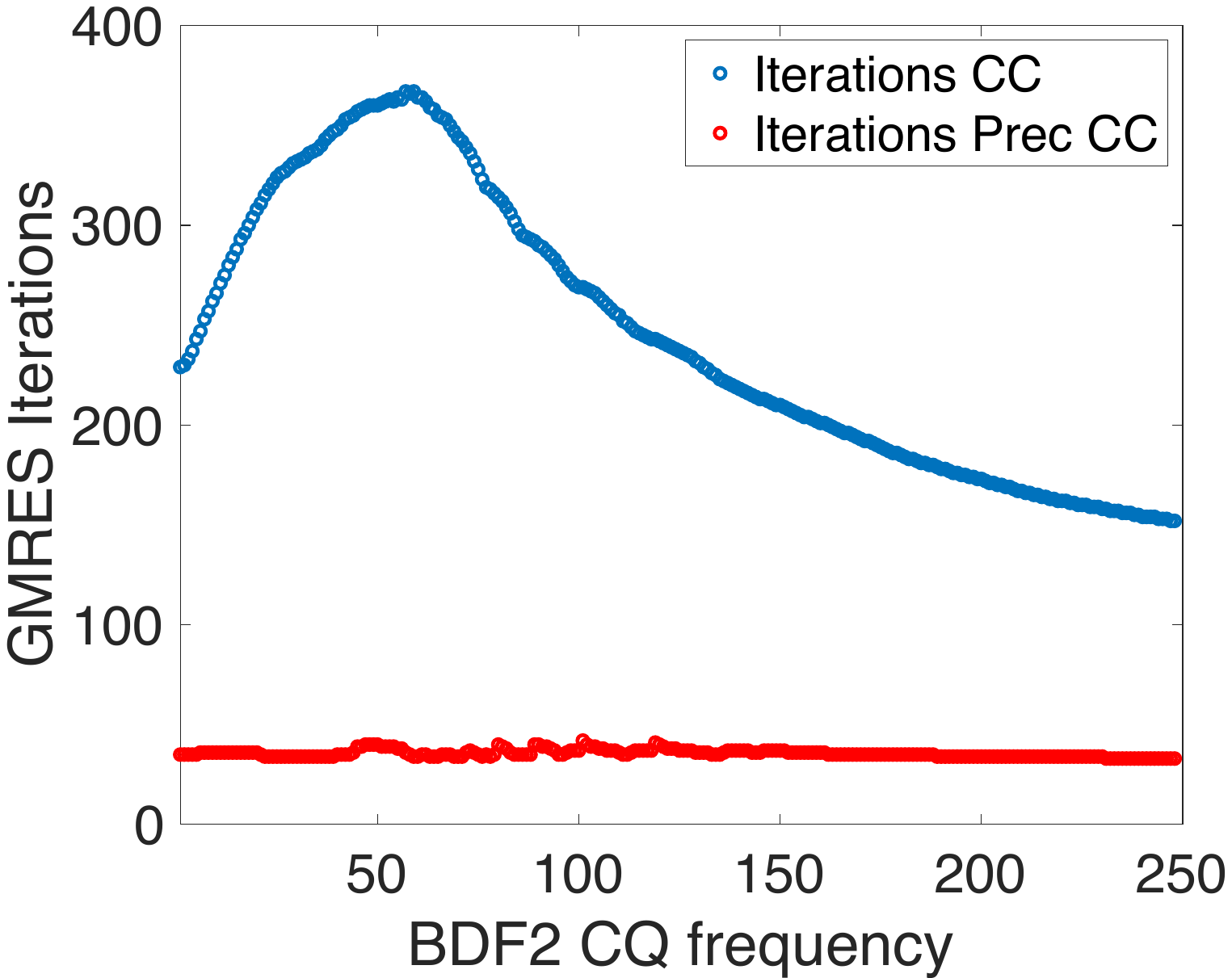}
\caption{Left panel: Orders of convergence of the near field CQ solutions of the wave equation in the exterior of a W-shaped scatterer at time $T=2$ with Neumann boundary conditions. We used the  Chebyshev based QBX discretizations with parameters $M=4,\ n_j=128, 1\leq j\leq 4$, $p=8$ and $\beta=4$ that incorporated the modified Fej\'er quadratures~\eqref{eq:fej1} and~\eqref{eq:fej2} for the solution of the ensemble of CC BIE formulations~\eqref{eq:CC}. The reference solution was produced using the CQ-RK5 discretization with $2048$ time steps and $N=1024$ Chebyshev points on the W-shaped boundary. Right panel: GMRES iterations required by the CC formulation and its preconditioned version to reach residuals of $10^{-7}$ for the Laplace domain modified Helmholtz problems corresponding to BDF2 frequencies for $T=2$ and $1024$ time steps; we considered only those frequencies whose Laplace domain BDF2 boundary data has a norm larger than $10^{-10}$.}
\label{fig:CQ_N_W}
\end{figure}

\begin{figure}
\centering
\includegraphics[height=60mm]{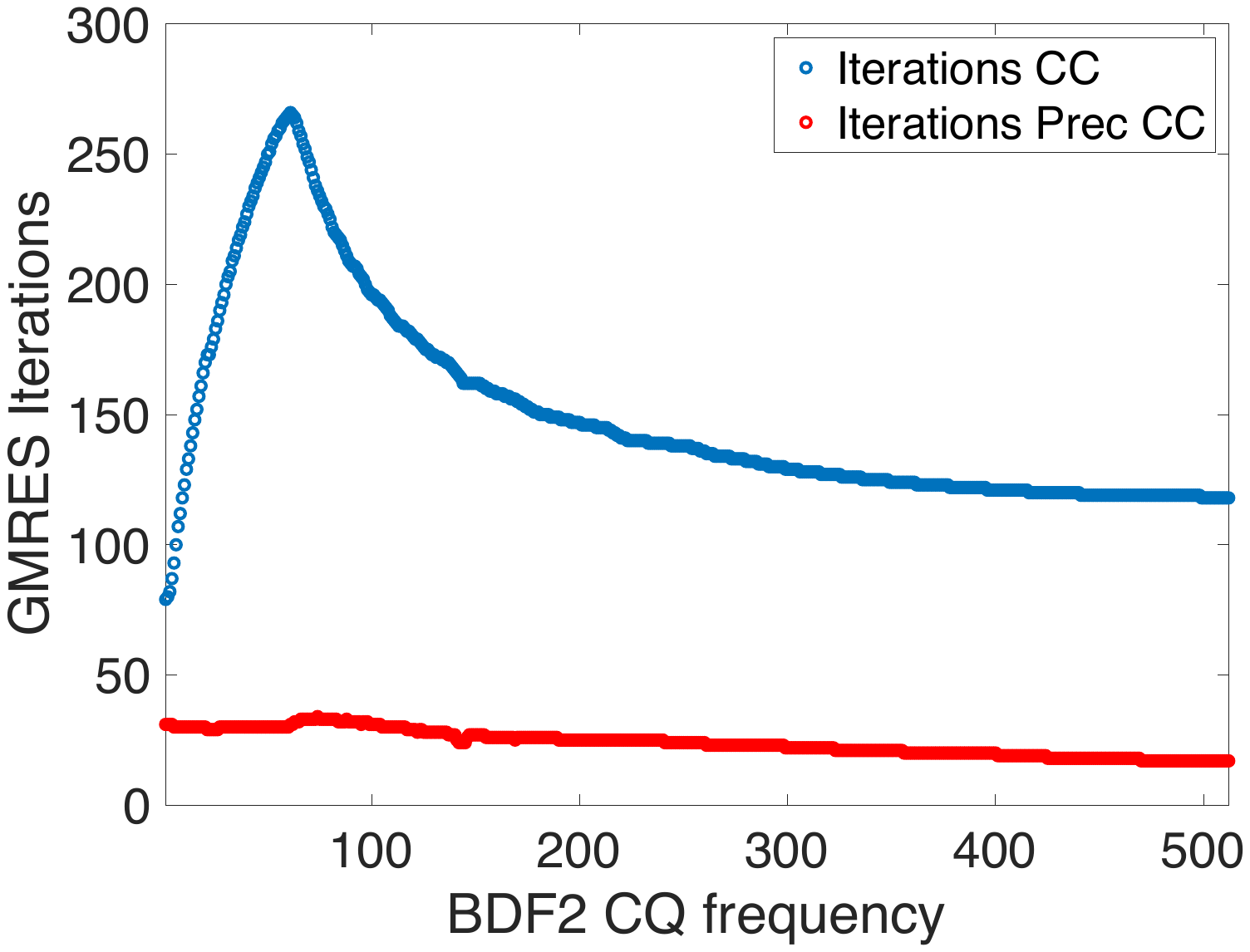}\includegraphics[height=60mm]{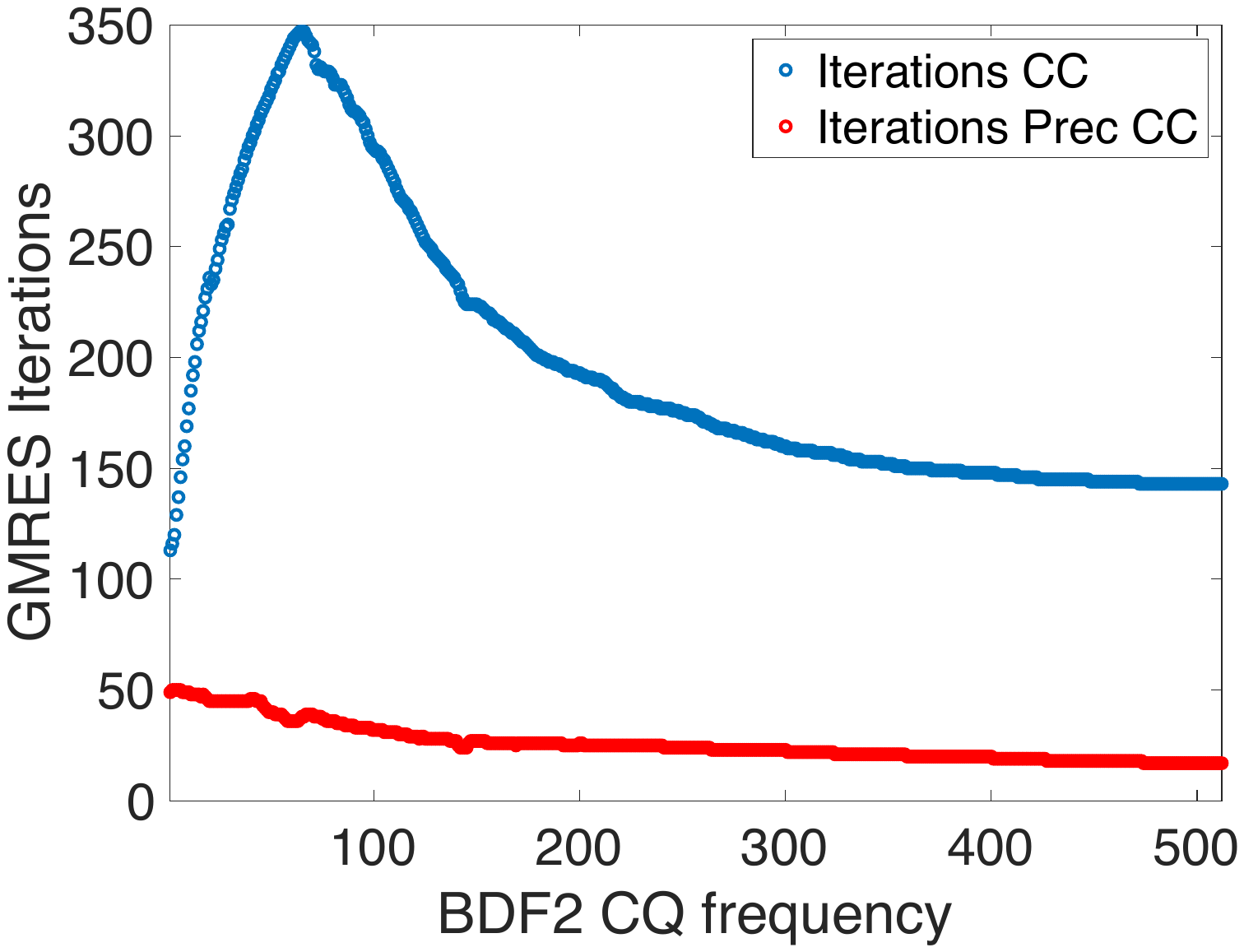}
\caption{Numbers of GMRES iterations to reach $10^{-7}$ residuals of the Nystr\"om CC and preconditioned CC~\eqref{eq:prec} systems for the V-shaped scatterer and the whole ensemble of Laplace domain frequencies associated with the CQ BDF2 formulation with $N=1024$ time steps and final time $T=2$. We considered only those frequencies whose Laplace domain BDF2 boundary data has a norm larger than $10^{-10}$}
\label{fig:CQ_N_iter}
\end{figure}

We continue in Figure~\ref{fig:CQ_N_teardrop} with an illustration of the high orders of convergence achieved by various CQ methods considered in this paper---up to and including fifth order, in the case of a time domain scattering problem off a closed non-smooth obstacle with Neumann boundary conditions. We considered BIE solvers based on both the weighted single layer and the double layer formulations of the CQ Laplace domain problems with Neumann boundary conditions, leading to BIE of the second kind~\eqref{Nsingle} in the former case and to the hypersingular BIE alternative CC formulation~\eqref{eq:CC} in the latter case. The ensemble of CQ second kind formulations~\eqref{Nsingle} was discretized using both Alpert quadratures and QBX methods of high enough orders, whereas the CC formulations~\eqref{Nsingle} were discretized using QBX. We continue in Figure~\ref{fig:CQ_N_arc} with an illustration of the high-order convergence achieved by CQ solvers in the case of a time domain scattering problem off a strip with Neumann boundary conditions. In this case only the hypersingular BIE formulation~\eqref{hypN} is available for the solution of Laplace domain problems, and we used both Alpert quadratures and QBX for its discretization---the former in conjunction with cosine changes of variables that resolve the end point singularities and the latter in the context of the CC reformulation~\eqref{eq:CC} of~\eqref{hypN}. Finally, we present in Figure~\ref{fig:CQ_N_vshape} and Figure~\ref{fig:CQ_N_W} time domain scattering experiments involving non-smooth open arcs with Neumann boundary conditions: one (top panel) and two V-shaped (bottom panel) obstacles in the former figure and a W-shaped obstacle in the latter figure. We used the modified Fej\'er quadratures~\eqref{eq:fej1} and~\eqref{eq:fej2} in the QBX discretizations for the solutions of the frequency domain CQ modified Helmholtz problems.  Since the CC formulation~\eqref{eq:CC} is of the first kind, we investigated the efficacy of the simple preconditioner~\eqref{eq:prec} for the Laplace domain problems corresponding to the BDF2 frequencies. As it can be seen from the plots in Figure~\ref{fig:CQ_N_W} and Figure~\ref{fig:CQ_N_iter} of GMRES iteration numbers required to reach $10^{-7}$ residuals for the unpreconditioned and preconditioned CC formulations, the preconditioner~\eqref{eq:prec} is highly effective. 

\begin{figure}
\centering
\includegraphics[height=105mm]{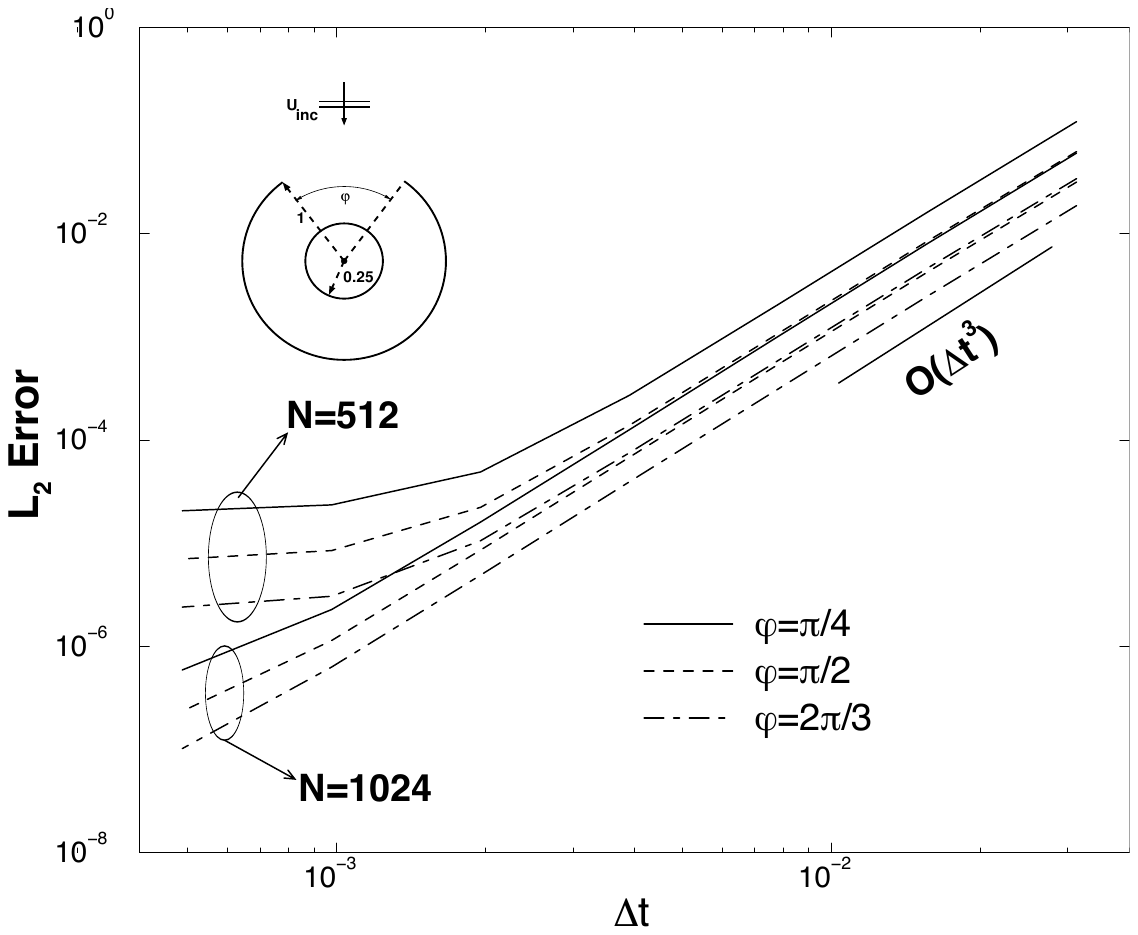}\includegraphics[height=115mm]{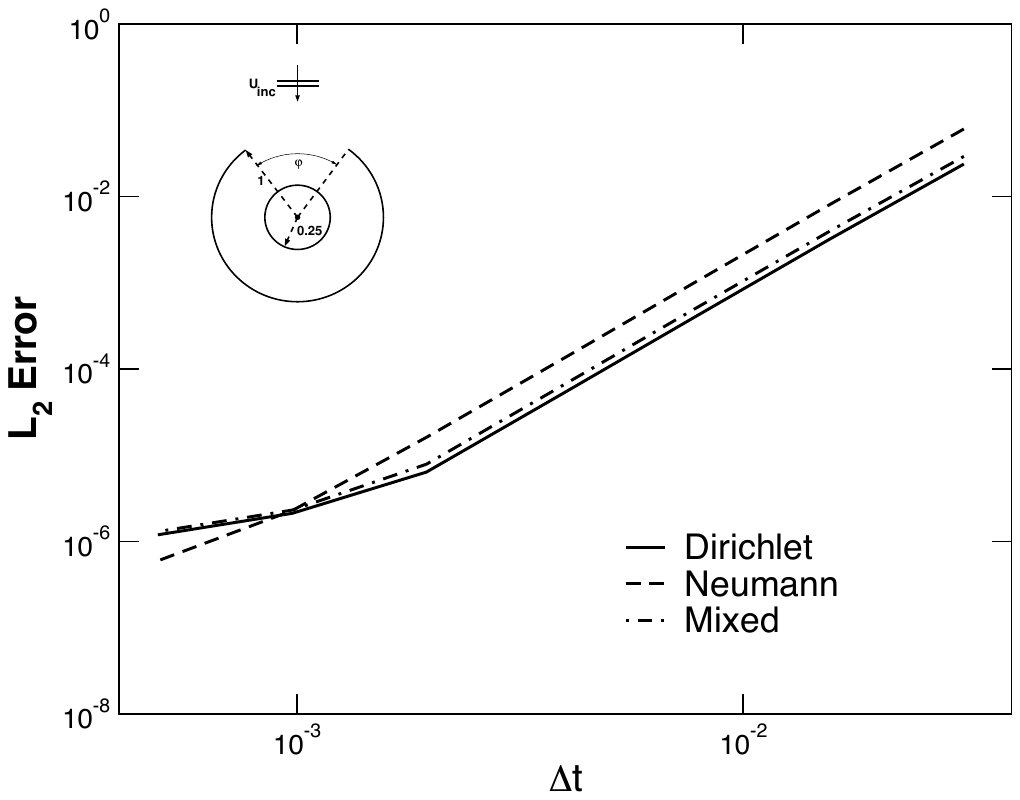}
\caption{Convergence orders of various CQ time integrators for scattering experiments involving slotted cylinders.}
\label{fig:cylinder-convergence}
\end{figure}

The next scattering experiments concern slotted cylinders, that is cavity like structures comprised of a circular scatterer enclosed by an open arc of small aperture, with both Dirichlet and Neumann boundary conditions. We present in Figure~\ref{fig:cylinder-convergence} the orders of convergence achieved by the various CQ time integrators (the final time was $T=16$ in these experiments) as the aperture angle varies when Dirichlet boundary conditions are imposed on both disconnected interfaces (left panel) as well as orders of convergence in the case when various types of boundary conditions are imposed on the two interfaces of the slotted cylinders (right panel). 

\begin{figure}
\centering
\includegraphics[height=50mm]{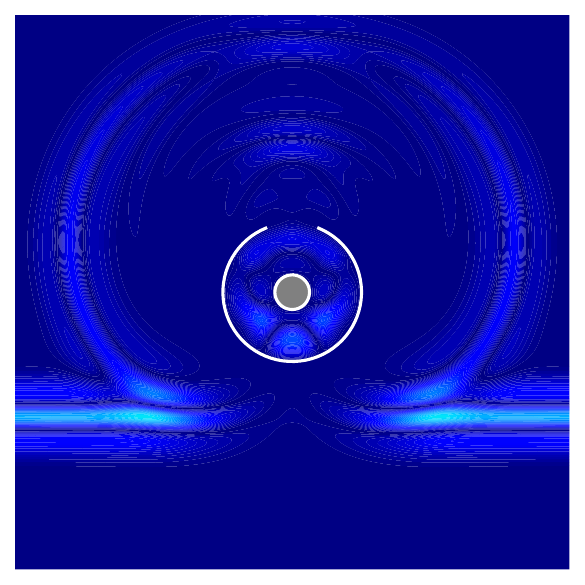}\includegraphics[height=50mm]{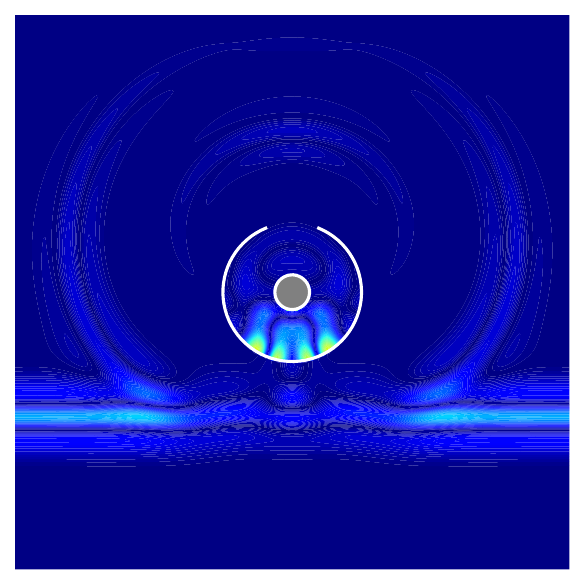}\includegraphics[height=50mm]{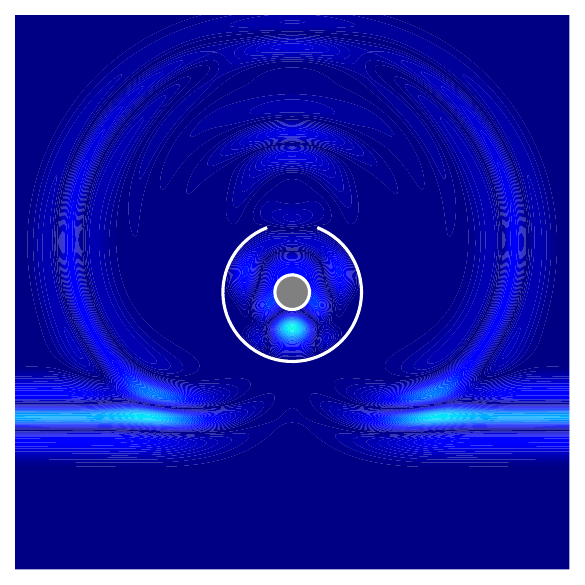}\\
\includegraphics[height=50mm]{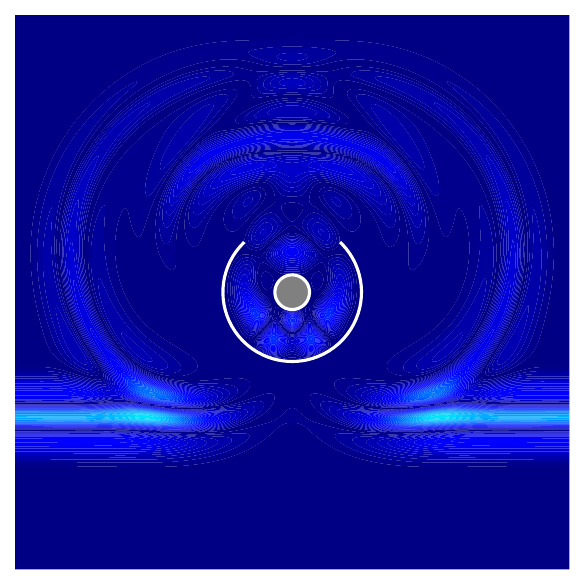}\includegraphics[height=50mm]{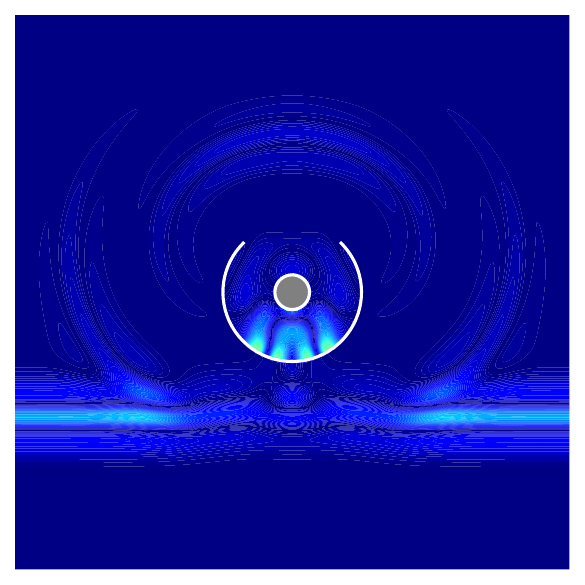}\includegraphics[height=50mm]{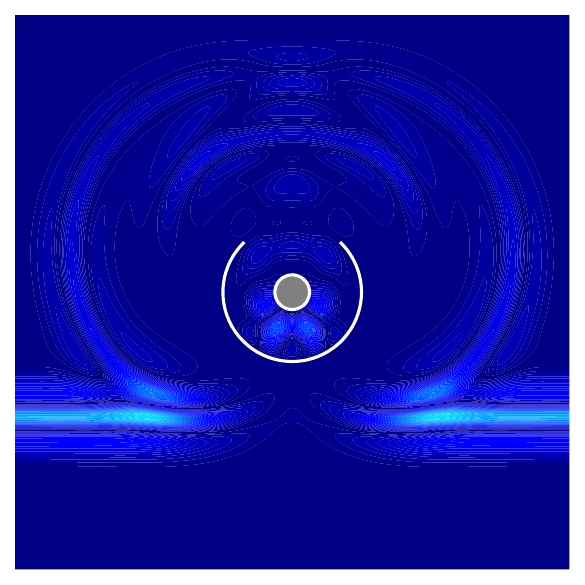}\\
\includegraphics[height=50mm]{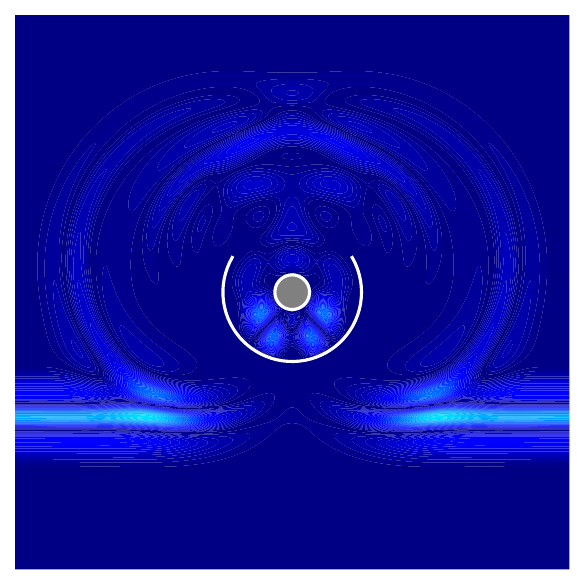}\includegraphics[height=50mm]{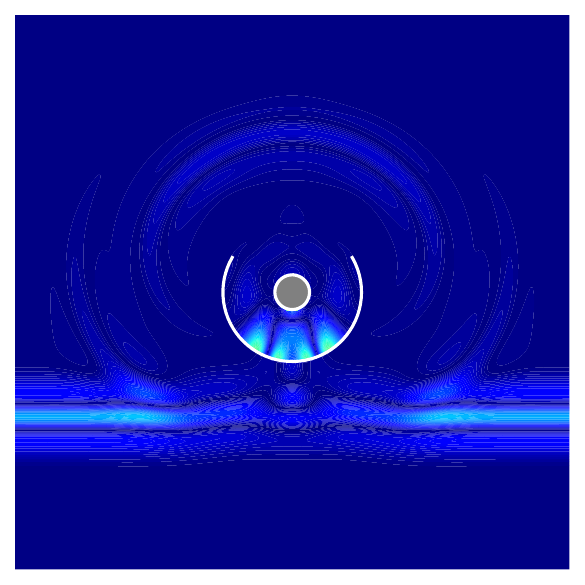}\includegraphics[height=50mm]{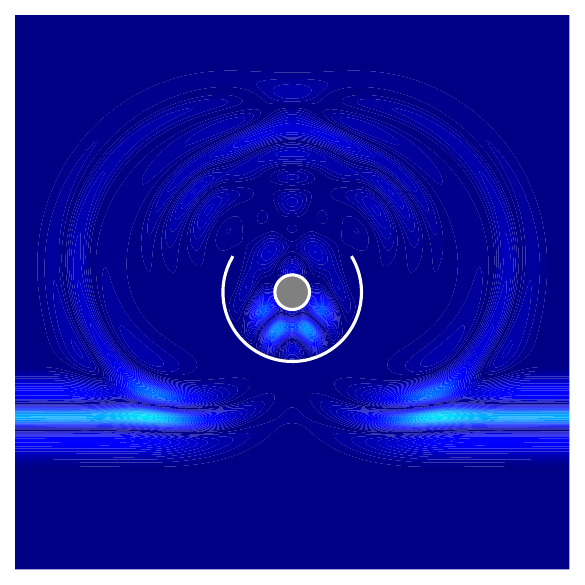}
\caption{Effect on total field of boundary condition and slot opening angle at a fixed time.}
\label{fig:apertureD}
\end{figure}

Finally, we present in Figures~\ref{fig:apertureD} time domain multiple scattering experiments involving a slotted cylinder over long periods of time (e.g., final time is $T=16$). The numerical simulations presented in these Figures were produced with CQ-RK3 methods using $2048$ time steps (that is $\Delta t=7.8\times 10^{-3}$) and $512$ discretization points per boundary for the discretization of the BIE formulation of the modified Helmholtz equations in the Laplace domain, parameters which result in accuracies at the level of $10^{-4}$ in the near field.

\begin{figure}
\centering
\includegraphics[height=105mm]{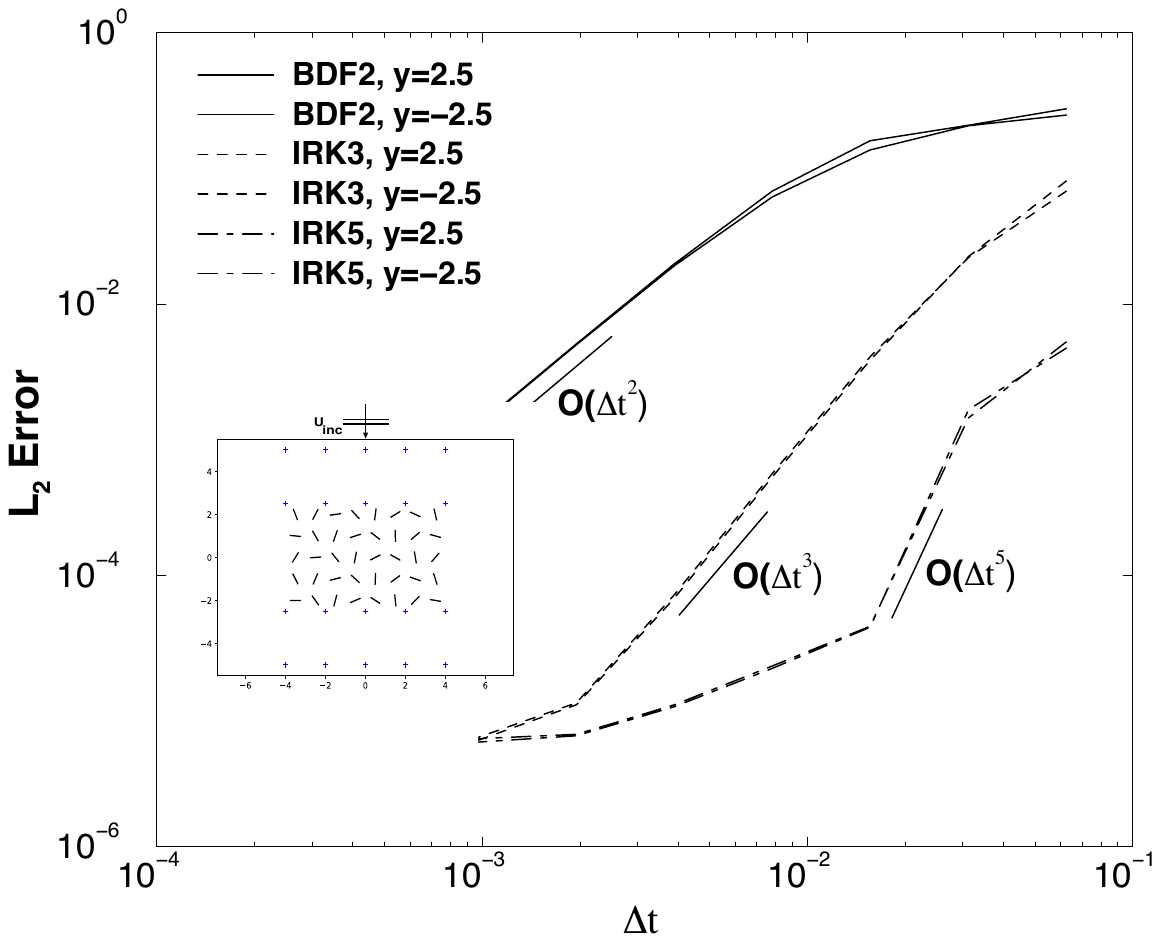}\includegraphics[height=105mm]{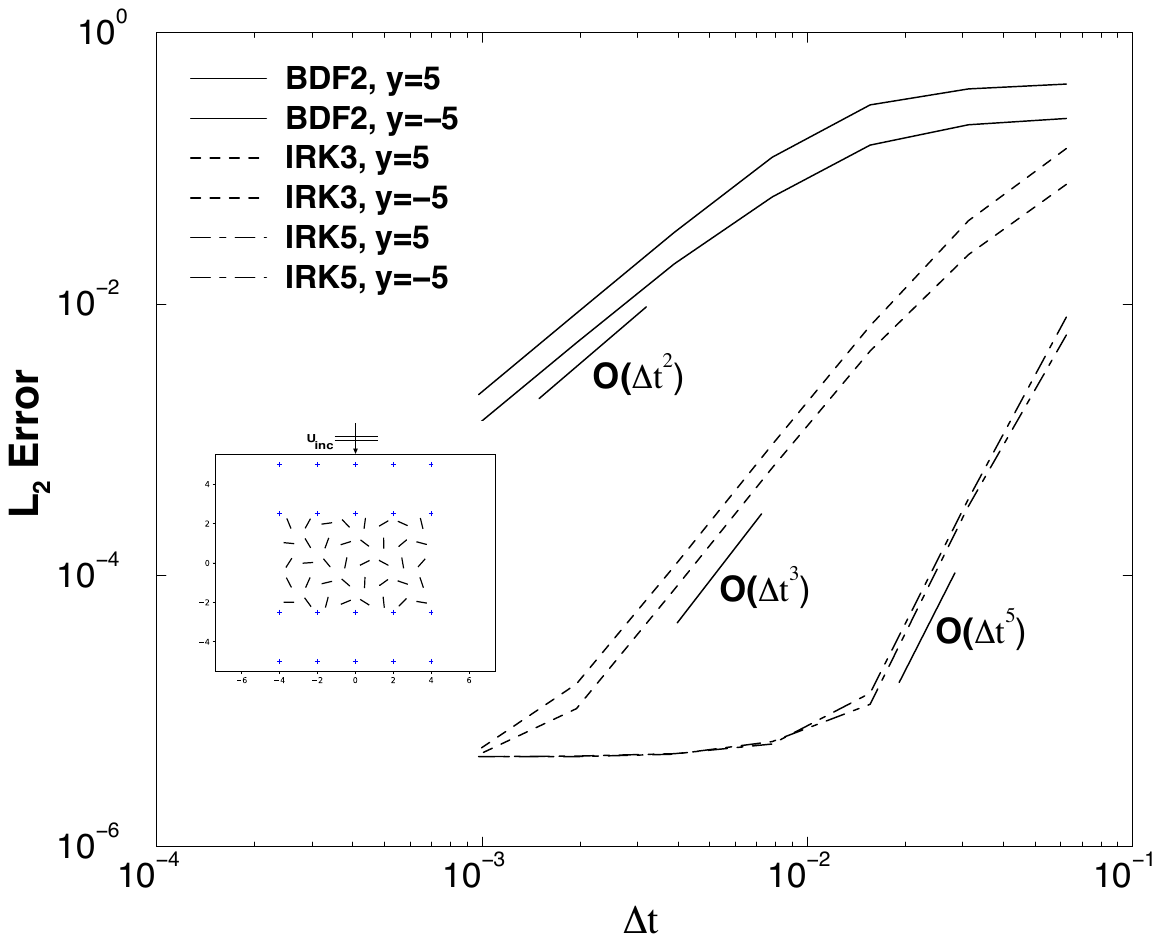}
\caption{Near field errors of CQ integrators for scattering experiments involving 40 randomly oriented strips.}
\label{fig:many-strips}
\end{figure}

We conclude the two dimensional experiments with an illustration in Figure~\ref{fig:many-strips} of the convergence rates of CQ solvers applied to time domain scattering by a collection of 40 randomly oriented strips. The errors reported were computed at two locations in the near field.

\subsection{3D Axisymmetric Scatterers}\label{3daxisym}

\begin{figure}
\centering
\includegraphics[height=90mm]{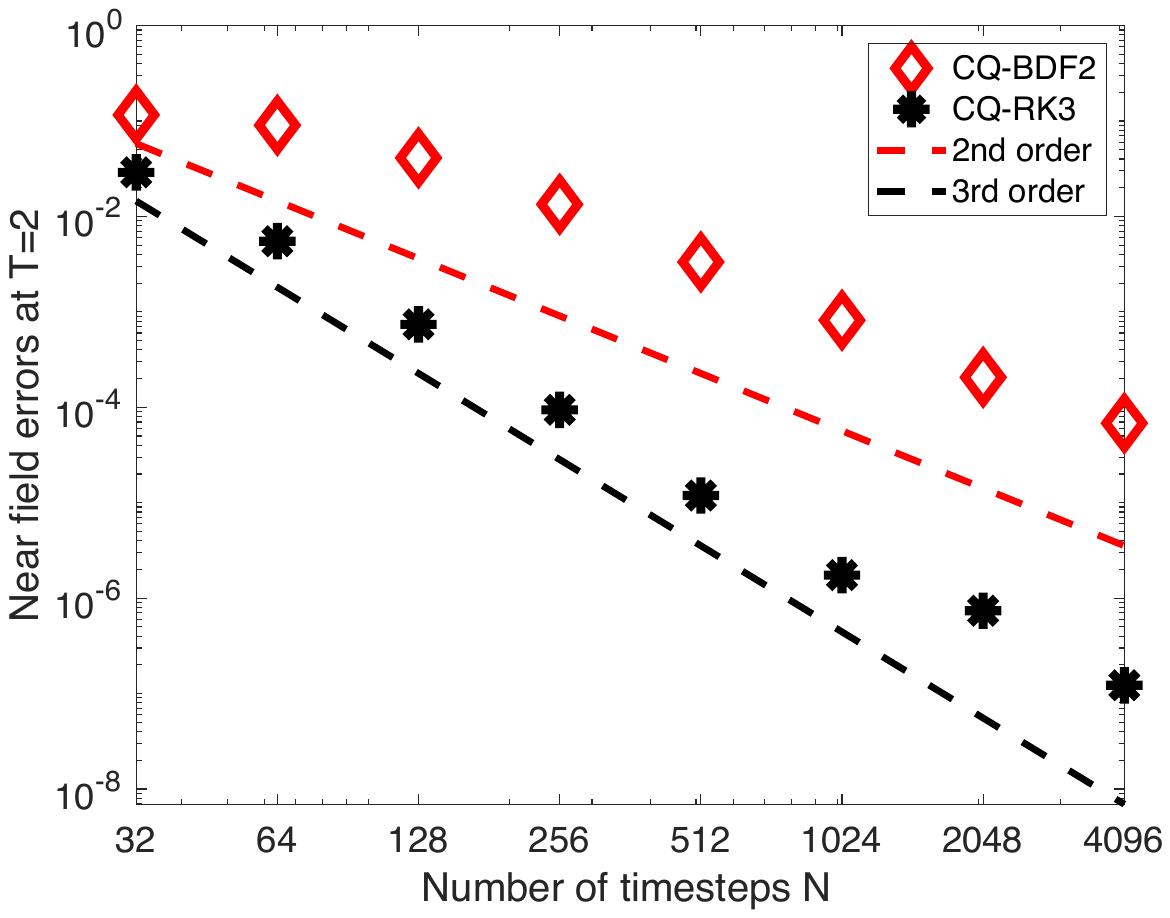}\includegraphics[height=90mm]{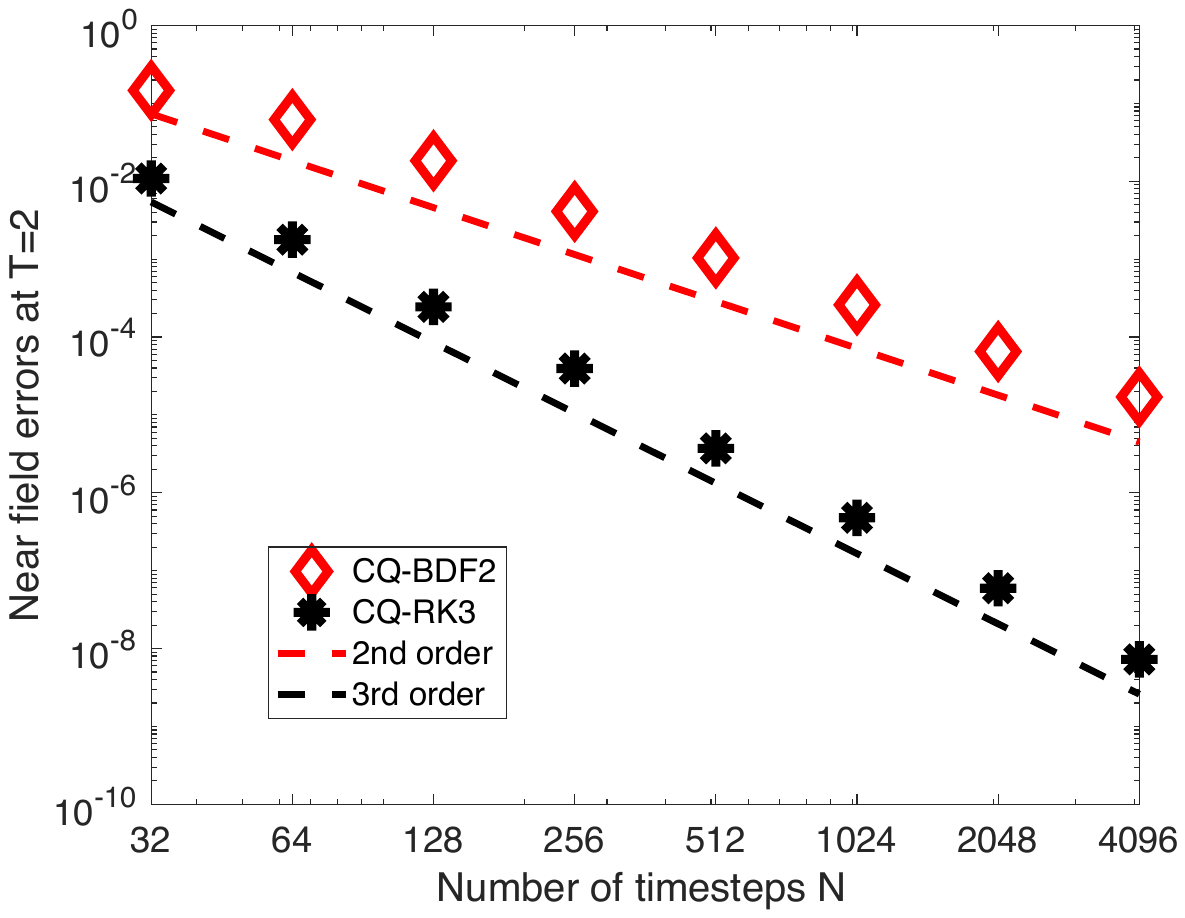}
\caption{Orders of convergence of the near field CQ solutions of the wave equation in the exterior of a torus (left panel) and an annular disc domain at time $T=2$ with Dirichlet boundary conditions. We used Alpert quadratures with the same parameters $\sigma=4$, $a=2$ and $m=3$ for the solution of the Fourier modal integral equations along the generated curve~\eqref{eq:Fourier} with a Fourier truncation parameter $M_F=512$ and all frequencies in the Laplace domain. In the case of the torus we used a double layer formulation analogue of formulations~\eqref{eq:Fourier} and global spacial equispaced mesh with $N=256$ discretization points on the boundary of the generating circle. In the case of the flat annulus domain we used the weighted single layer formulation~\eqref{eq:Fourier} and a graded mesh with parameter $p=3$ and $N=128$ points along the generating arc segment. The reference solution was produced using the CQ-RK5 discretization with $2048$ time steps and $N=256$ boundary discretization points for each geometrical configuration.}
\label{fig:CQ_D_3d}
\end{figure}

We present in this section numerical experiments concerning time domain scattering from axisymmetric obstacles in three dimensions. Specifically, we consider the following objects: a torus whose generating curve is a circle of radius 1 centered at $(2,0,0)$, a flat annulus (an open surface in three dimensions) whose generating curve is the line segment (arc) in the $xz$-plane joining the points $(1,0,0)$ to $(1/2,0,0)$, and a spherical cavity shell (an open surface as well) of aperture angle $pi/4$ whose generating curve is the arc of the circle of radius one connecting the points $(0,0,-1)$ to $(\sin(\pi/8),0,\cos(\pi/8))$. In all the scattering experiments we considered the following incident field 
\begin{equation}\label{eq:3d_inc}
    g(x,t) = \cos(5t-x \cdot \mathbf{d} )\exp( -1.5(5t-x\cdot \mathbf{d} -5)^2 ),
\end{equation}
where the direction of the incident Gaussian modulated plane wave incident field is given by $\mathbf{d}=(\sin(\pi/4)\cos(\pi/3),\sin(\pi/4)\sin(\pi/3),-\cos(\pi/4))$.  We present in Figure~\ref{fig:CQ_D_3d} the orders of convergence achieved by CQ BDF2 and CQ RK3 time domain solvers for the torus and annular geometries with Dirichlet boundary conditions. The numerical errors are near field errors at time $T=2$ where the $256\times 256$ observation points were placed on spheres of radii 2 and respectively 4 surrounding the obstacles. In each case we solved all of the Fourier modal integral equations along the generated curve~\eqref{eq:Fourier} using a Fourier truncation parameter $M_F=512$ for all complex frequencies in the Laplace domain required by the CQ algorithms. Each integral equation on the generated curve was solved using Alpert quadratures Nystr\"om methods based on (a) equispaced meshes in the case of the circular geometry that generates the torus through rotation around the $z$ axis and (b) sigmoid meshes with parameter $p=3$ in the case of the line segment that generates the annulus via rotation around the $z$ axis. The reference solutions were computed using CQ-RK5 with $2048$ time steps and fine meshes on the generating curve for the Nystr\"om Alpert discretization of the ensemble of modal Fourier integral equations for all frequencies in the Laplace domain. As it can be seen in Figure~\ref{fig:CQ_D_3d}, the CQ solvers in conjunction with Nystr\"om Alpert axisymmetric discretizations of the BIE formulations of frequency domain problems converge to high-order and are capable to deliver accuracies at the level $10^{-7}$ (and better) for both smooth surfaces as well as open surfaces with edges.

\begin{figure}
\centering
\includegraphics[height=50mm]{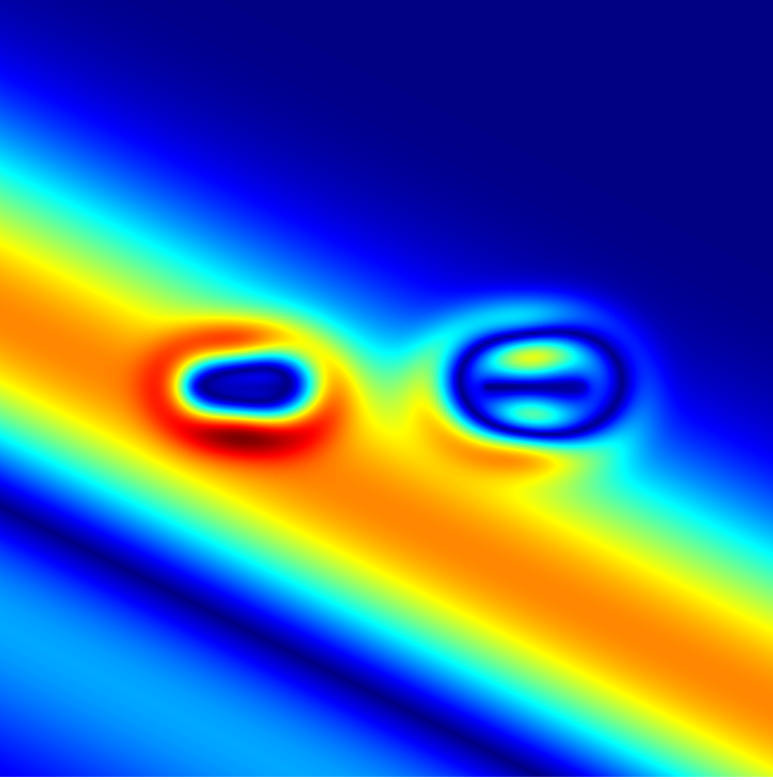}\includegraphics[height=50mm]{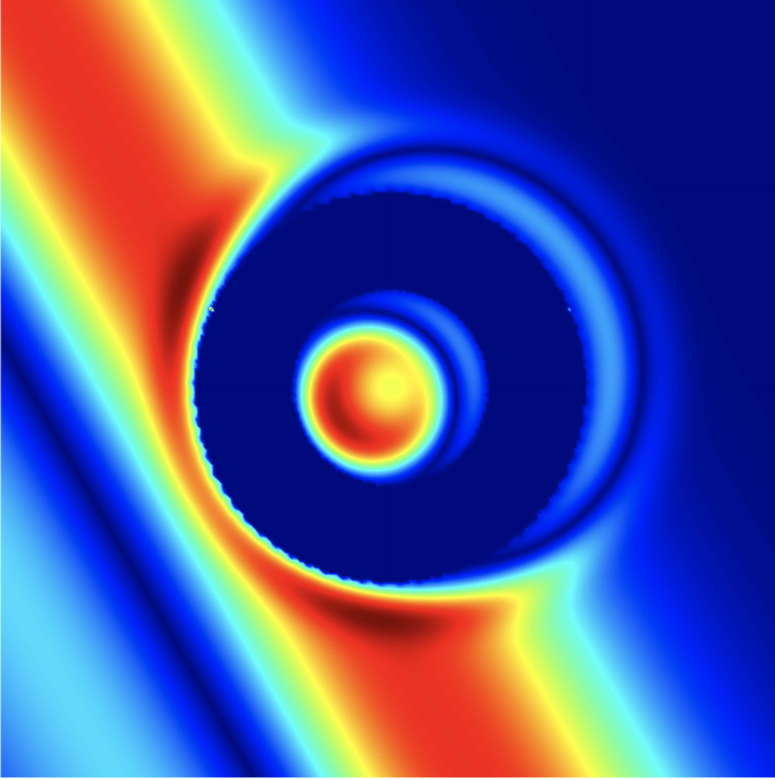}\includegraphics[height=50mm]{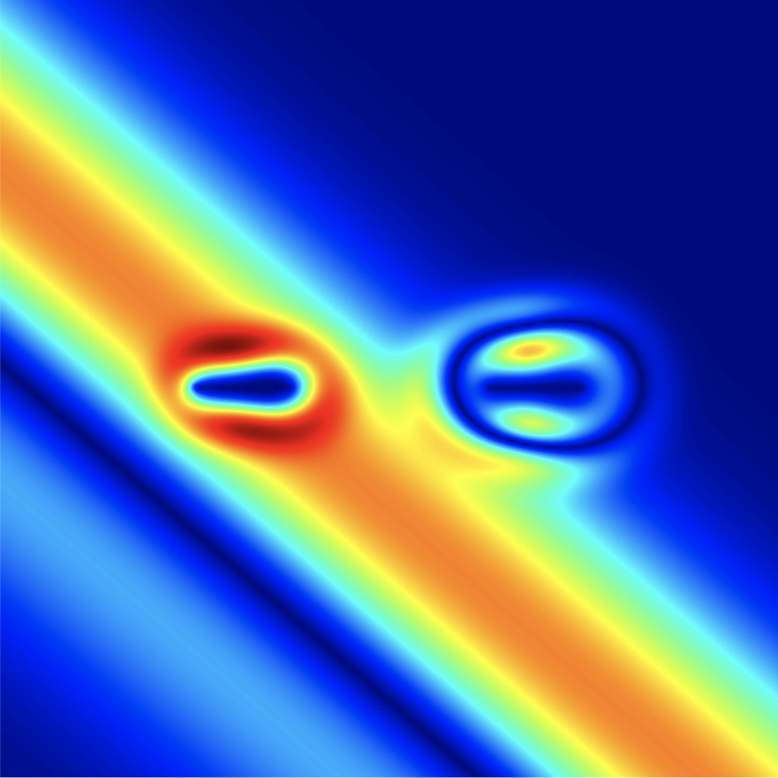}\\
\includegraphics[height=50mm]{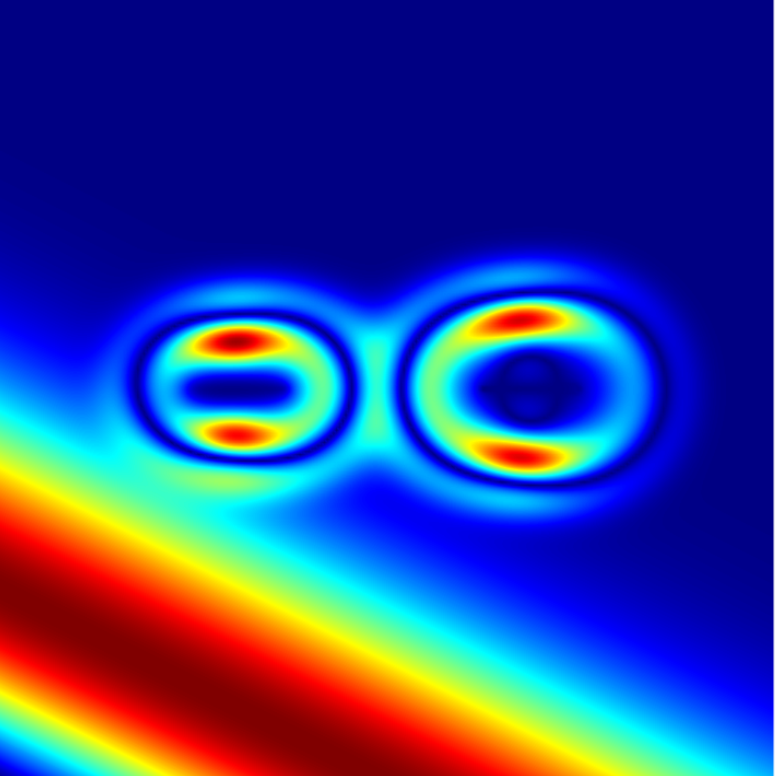}\includegraphics[height=50mm]{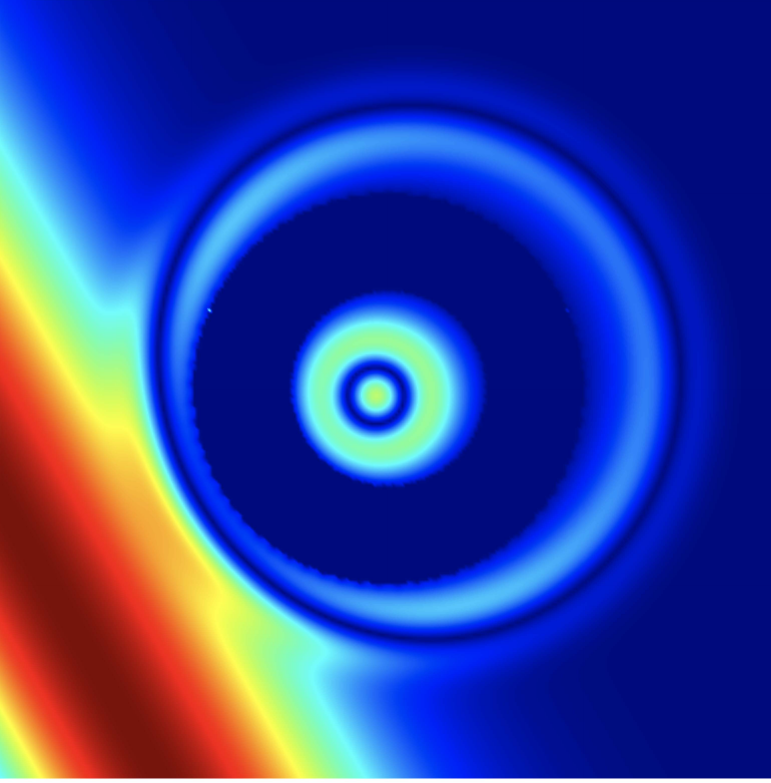}
\includegraphics[height=50mm]{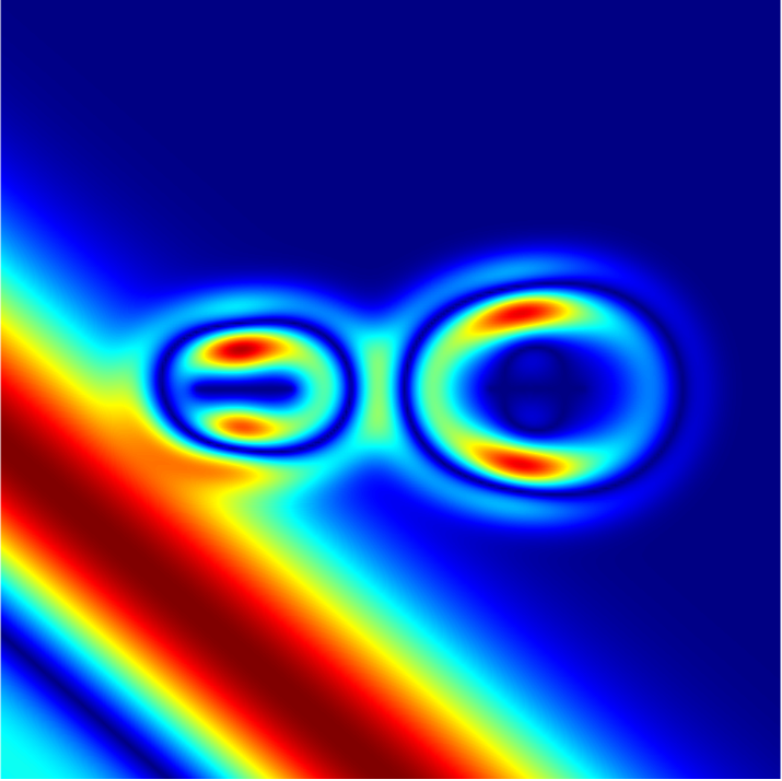}\\
\includegraphics[height=50mm]{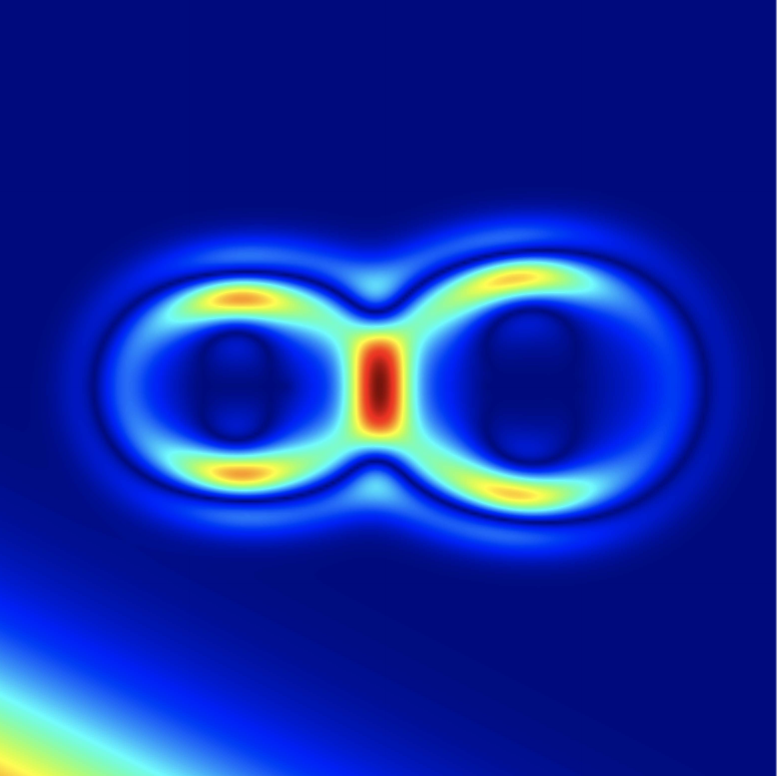}\includegraphics[height=50mm]{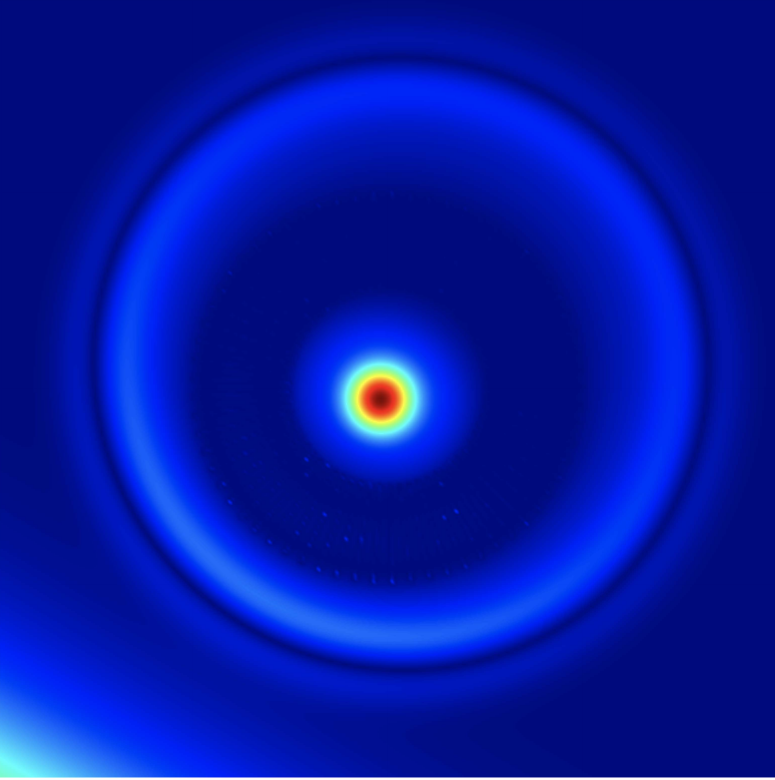}
\includegraphics[height=50mm]{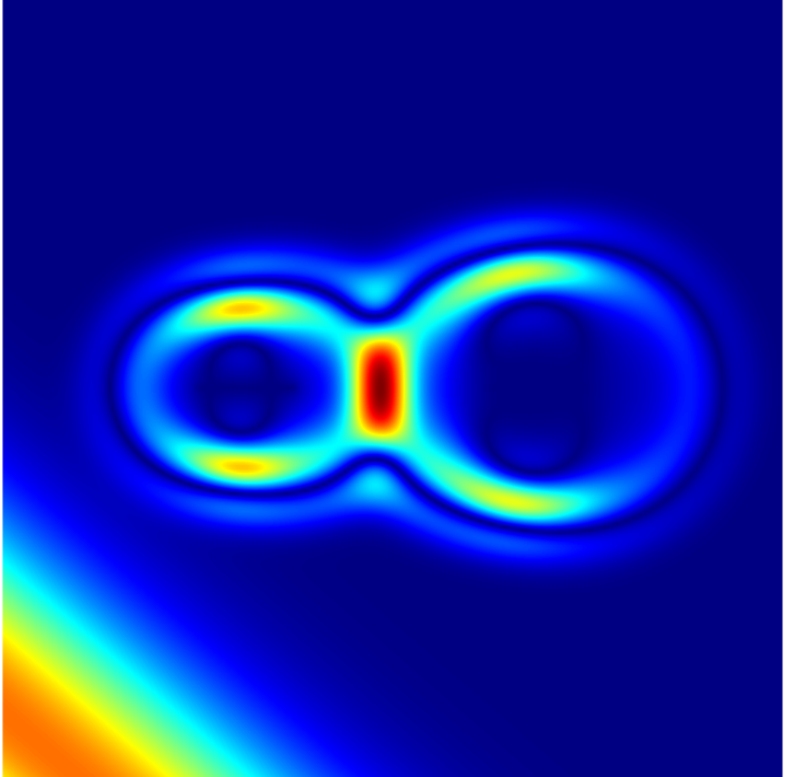}
\caption{Time traces of the total field near an annular obstacles with Dirichlet boundary conditions and incident field~\eqref{eq:3d_inc}. The row panels depict the fields on the $xz$-plane (left), $xy$-plane (middle) and $yz$-plane(right).}
\label{fig:aperture3Da}
\end{figure}

\begin{figure}
\centering
\includegraphics[height=50mm]{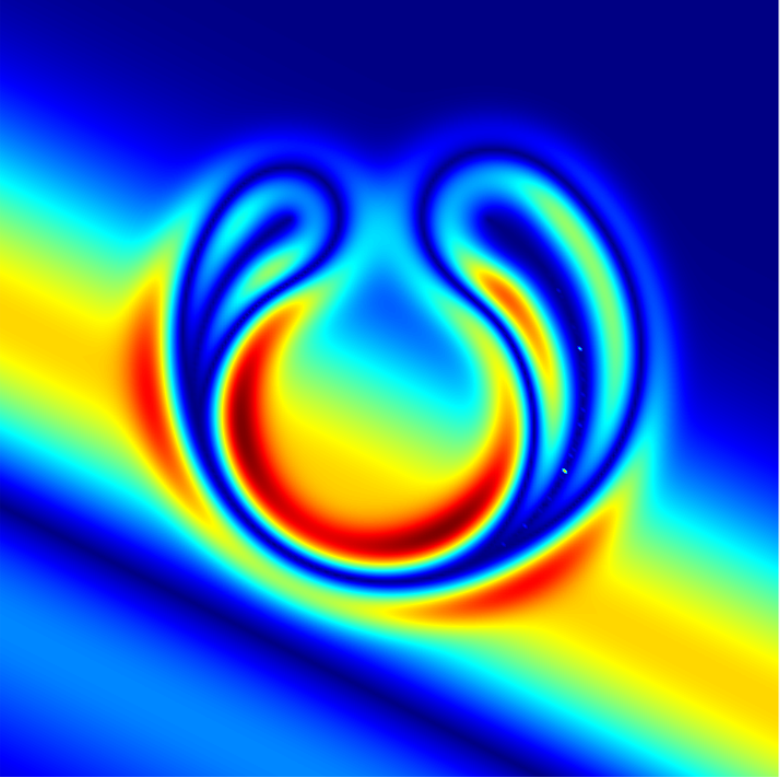}\includegraphics[height=50mm]{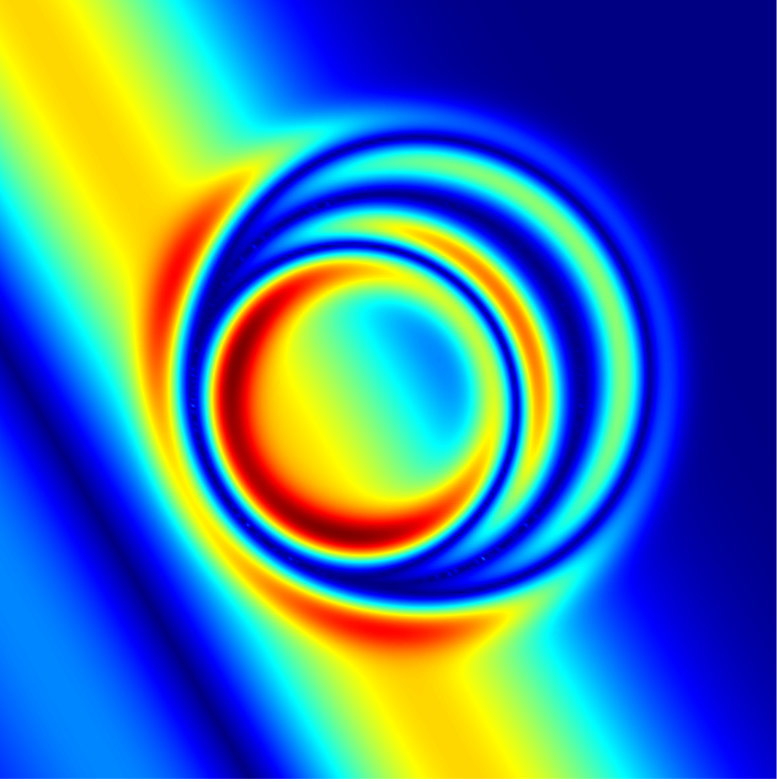}\\\
\includegraphics[height=50mm]{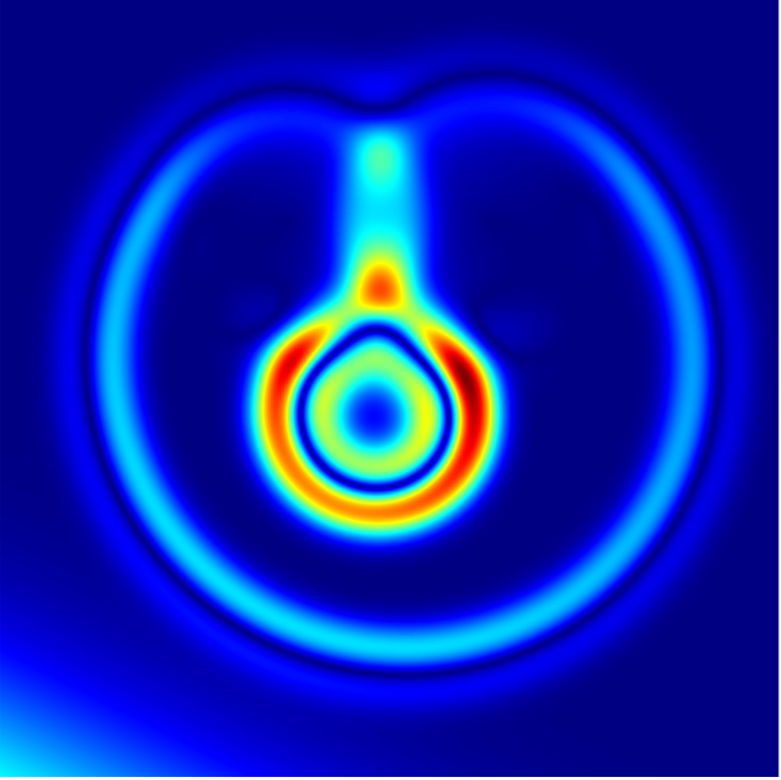}\includegraphics[height=50mm]{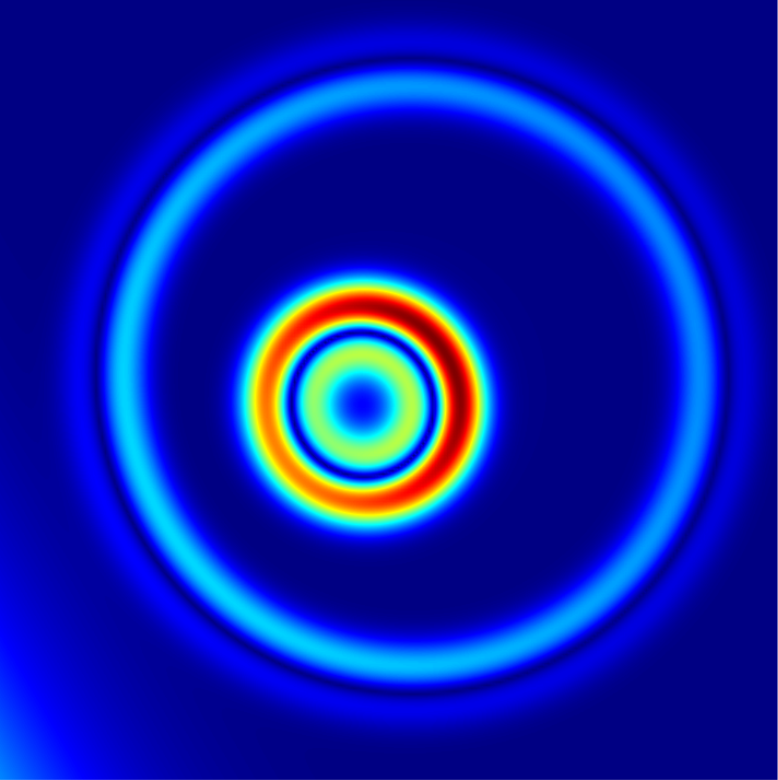}
\caption{Time traces of the total field near an annular obstacles with Dirichlet boundary conditions and incident field~\eqref{eq:3d_inc}. The row panels depict the fields on the $xz$-plane (left)and the $xy$-plane (right).}
\label{fig:aperture3Db}
\end{figure}

We present in Figure~\ref{fig:aperture3Da} and Figure~\ref{fig:aperture3Db} time traces of the total field around the annular obstacles and the spherical cavity shell.

\section{Summary and Conclusions}

We have shown that CQ methods used in conjunction with Nystr\"om discretizations for the BIE solution of their associated Laplace domain Helmholtz problems can deliver high-order numerical solutions of two dimensional wave equations in unbounded domains with non smooth open/closed boundaries.  In particular, the CQ-RK5 methods require the largest time steps to reach desired accuracy levels, and more importantly the gains in accuracy over the CQ-BDF2 and CQ-RK3 are not outweighed by the increased number of Laplace domain problems that ought to be solved in this CQ version. We focused on two classes of high-order Nystr\"om discretizations, one based on Alpert quadratures, and one on QBX methods. Our choice was motivated by the fact that both of these quadratures are rather straightforward to implement,  are applicable seamlessly to real as well as complex wavenumbers, are compatible with fast algorithms, and can be extended to three dimensional settings. The extension of the CQ Nystr\"om methodology presented in this paper to the solution of three dimensional wave equations, for scattering problems beyond the axisymmetric case presented herein, and the acceleration of the solvers via fast algorithms are currently underway.  

\section*{Acknowledgments}
 Catalin Turc gratefully acknowledges support from NSF through contract DMS-1908602.


\bibliographystyle{acm}
\bibliography{CQ_EPC} 

\begin{thebibliography}{10}

\bibitem{alpert1999hybrid}
Bradley~K Alpert.
\newblock Hybrid gauss-trapezoidal quadrature rules.
\newblock {\em SIAM Journal on Scientific Computing}, 20(5):1551--1584, 1999.

\bibitem{turc_corner_N}
A.~Anand, J.~S. Ovall, and C.~Turc.
\newblock Well-conditioned boundary integral equations for two-dimensional
  sound-hard scattering problems in domains with corners.
\newblock {\em J. Integral Equations Appl.}, 24(3):321--358, 2012.

\bibitem{anand2012well}
Akash Anand, Jeff Ovall, Catalin Turc, et~al.
\newblock Well-conditioned boundary integral equations for two-dimensional
  sound-hard scattering problems in domains with corners.
\newblock {\em Journal of Integral Equations and Applications}, 24(3):321--358,
  2012.

\bibitem{anderson2020high}
Thomas~G Anderson, Oscar~P Bruno, and Mark Lyon.
\newblock High-order, dispersionless fast-hybrid wave equation solver. part i:
  O (1) sampling cost via incident-field windowing and recentering.
\newblock {\em SIAM Journal on Scientific Computing}, 42(2):A1348--A1379, 2020.

\bibitem{banjai2010multistep}
Lehel Banjai.
\newblock Multistep and multistage convolution quadrature for the wave
  equation: algorithms and experiments.
\newblock {\em SIAM Journal on Scientific Computing}, 32(5):2964--2994, 2010.

\bibitem{banjai2011runge}
Lehel Banjai, Christian Lubich, and Jens~Markus Melenk.
\newblock Runge--kutta convolution quadrature for operators arising in wave
  propagation.
\newblock {\em Numerische Mathematik}, 119(1):1--20, 2011.

\bibitem{banjai2009rapid}
Lehel Banjai and Stefan Sauter.
\newblock Rapid solution of the wave equation in unbounded domains.
\newblock {\em SIAM Journal on Numerical Analysis}, 47(1):227--249, 2009.

\bibitem{barnett2020high}
Alex Barnett, Leslie Greengard, and Thomas Hagstrom.
\newblock High-order discretization of a stable time-domain integral equation
  for 3d acoustic scattering.
\newblock {\em Journal of Computational Physics}, 402:109047, 2020.

\bibitem{bendali2012extension}
Abderrahmane Bendali, Francis Collino, Mbarek Fares, and Bassam Steif.
\newblock Extension to nonconforming meshes of the combined current and charge
  integral equation.
\newblock {\em IEEE transactions on antennas and propagation},
  60(10):4732--4744, 2012.

\bibitem{betcke2017overresolving}
Timo Betcke, Nicolas Salles, and Wojciech Smigaj.
\newblock Overresolving in the laplace domain for convolution quadrature
  methods.
\newblock {\em SIAM Journal on Scientific Computing}, 39(1):A188--A213, 2017.

\bibitem{turc1}
Y.~Boubendir and C.~Turc.
\newblock Wave-number estimates for regularized combined field boundary
  integral operators in acoustic scattering problems with neumann boundary
  conditions.
\newblock {\em IMA Journal of Numerical Analysis}, 33(4):1176--1225, 2013.

\bibitem{bremer2012nystrom}
James Bremer.
\newblock On the nystr{\"o}m discretization of integral equations on planar
  curves with corners.
\newblock {\em Applied and Computational Harmonic Analysis}, 32(1):45--64,
  2012.

\bibitem{bruno2012second}
Oscar~P Bruno and St{\'e}phane~K Lintner.
\newblock Second-kind integral solvers for te and tm problems of diffraction by
  open arcs.
\newblock {\em Radio Science}, 47(06):1--13, 2012.

\bibitem{chapko2000numerical}
Roman Chapko, Rainer Kress, and Lars M{\"o}nch.
\newblock On the numerical solution of a hypersingular integral equation for
  elastic scattering from a planar crack.
\newblock {\em IMA journal of numerical analysis}, 20(4):601--619, 2000.

\bibitem{cohl1999compact}
Howard~S Cohl and Joel~E Tohline.
\newblock A compact cylindrical green’s function expansion for the solution
  of potential problems.
\newblock {\em The astrophysical journal}, 527(1):86, 1999.

\bibitem{ColtonKress}
D.~Colton and R.~Kress.
\newblock {\em Inverse acoustic and electromagnetic scattering theory},
  volume~93 of {\em Applied Mathematical Sciences}.
\newblock Springer-Verlag, Berlin, second edition, 1998.

\bibitem{conway2010exact}
John~T Conway and Howard~S Cohl.
\newblock Exact fourier expansion in cylindrical coordinates for the
  three-dimensional helmholtz green function.
\newblock {\em Zeitschrift f{\"u}r angewandte Mathematik und Physik},
  61:425--443, 2010.

\bibitem{dominguez2012fully}
Victor Dominguez, Sijiang~L Lu, and Francisco-Javier Sayas.
\newblock A fully discrete calder{\'o}n calculus for two dimensional time
  harmonic waves.
\newblock {\em arXiv preprint arXiv:1210.7017}, 2012.

\bibitem{dominguez2014nystrom}
V{\'\i}ctor Dom{\'\i}nguez, Sijiang~L Lu, and Francisco-Javier Sayas.
\newblock A nystr{\"o}m flavored calder{\'o}n calculus of order three for two
  dimensional waves, time-harmonic and transient.
\newblock {\em Computers \& Mathematics with Applications}, 67(1):217--236,
  2014.

\bibitem{dominguez2015well}
Victor Dominguez, Mark Lyon, and Catalin Turc.
\newblock Well-posed boundary integral equation formulations and
  nystr$\backslash$" om discretizations for the solution of helmholtz
  transmission problems in two-dimensional lipschitz domains.
\newblock {\em Journal of Integral Equations and Applications}, 28(3):395--440,
  2016.

\bibitem{dominguez2008dirac}
V{\'\i}ctor Dom{\'\i}nguez, M-L Rap{\'u}n, and F-J Sayas.
\newblock Dirac delta methods for helmholtz transmission problems.
\newblock {\em Advances in Computational Mathematics}, 28(2):119--139, 2008.

\bibitem{dominguez2015fully}
V{\'\i}ctor Dom{\'\i}nguez, Tonatiuh S{\'a}nchez-Vizuet, and Francisco-Javier
  Sayas.
\newblock A fully discrete calder{\'o}n calculus for the two-dimensional
  elastic wave equation.
\newblock {\em Computers \& Mathematics with Applications}, 69(7):620--635,
  2015.

\bibitem{epstein2010debye}
Charles~L Epstein and Leslie Greengard.
\newblock Debye sources and the numerical solution of the time harmonic maxwell
  equations.
\newblock {\em Communications on Pure and Applied Mathematics: A Journal Issued
  by the Courant Institute of Mathematical Sciences}, 63(4):413--463, 2010.

\bibitem{epstein2019high}
Charles~L Epstein, Leslie Greengard, and Michael O'Neil.
\newblock A high-order wideband direct solver for electromagnetic scattering
  from bodies of revolution.
\newblock {\em Journal of Computational Physics}, 387:205--229, 2019.

\bibitem{et1986formulation}
A~Bamberger et~T.~Ha~Duong and JC~Nedelec.
\newblock Formulation variationnelle espace-temps pour le calcul par potentiel
  retard{\'e} de la diffraction d'une onde acoustique (i).
\newblock {\em Mathematical methods in the applied sciences}, 8(1):405--435,
  1986.

\bibitem{ganesh2023high}
M~Ganesh and F~Le~Lou{\"e}r.
\newblock A high-order algorithm for time-domain scattering in three
  dimensions.
\newblock {\em Advances in Computational Mathematics}, 49(4):46, 2023.

\bibitem{hao2014high}
Sijia Hao, Alex~H Barnett, Per-Gunnar Martinsson, and P~Young.
\newblock High-order accurate methods for nystr{\"o}m discretization of
  integral equations on smooth curves in the plane.
\newblock {\em Advances in Computational Mathematics}, 40(1):245--272, 2014.

\bibitem{helsing2020extended}
Johan Helsing and Anders Karlsson.
\newblock An extended charge-current formulation of the electromagnetic
  transmission problem.
\newblock {\em SIAM Journal on Applied Mathematics}, 80(2):951--976, 2020.

\bibitem{helsing2008corner}
Johan Helsing and Rikard Ojala.
\newblock Corner singularities for elliptic problems: Integral equations,
  graded meshes, quadrature, and compressed inverse preconditioning.
\newblock {\em Journal of Computational Physics}, 227(20):8820--8840, 2008.

\bibitem{klockner2013quadrature}
Andreas Kl{\"o}ckner, Alexander Barnett, Leslie Greengard, and Michael O'Neil.
\newblock Quadrature by expansion: A new method for the evaluation of layer
  potentials.
\newblock {\em Journal of Computational Physics}, 252:332--349, 2013.

\bibitem{KressCorner}
R.~Kress.
\newblock A {N}ystr\"om method for boundary integral equations in domains with
  corners.
\newblock {\em Numer. Math.}, 58(2):145--161, 1990.

\bibitem{KressH}
R.~Kress.
\newblock On the numerical solution of a hypersingular integral equation in
  scattering theory.
\newblock {\em J. Comput. Appl. Math.}, 61(3):345--360, 1995.

\bibitem{kusmaul}
R.~Kussmaul.
\newblock Ein numerisches {V}erfahren zur {L}\"osung des {N}eumannschen
  {N}eumannschen {A}ussenraumproblems f\"ur die {H}elmholtzsche
  {S}chwingungsgleichung.
\newblock {\em Computing (Arch. Elektron. Rechnen)}, 4:246--273, 1969.

\bibitem{Labarca2019convolution}
Ignacio Labarca, Luiz~M Faria, and Carlos P{\'e}rez-Arancibia.
\newblock Convolution quadrature methods for time-domain scattering from
  unbounded penetrable interfaces.
\newblock {\em Proceedings of the Royal Society A}, 475(2227):20190029, 2019.

\bibitem{lu2014efficient}
Wangtao Lu and Ya~Yan Lu.
\newblock Efficient high order waveguide mode solvers based on boundary
  integral equations.
\newblock {\em Journal of Computational Physics}, 272:507--525, 2014.

\bibitem{martensen}
{E.} Martensen.
\newblock \"{U}ber eine {M}ethode zum r\"aumlichen {N}eumannschen {P}roblem mit
  einer {A}nwendung f\"ur torusartige {B}erandungen.
\newblock {\em Acta Math.}, 109:75--135, 1963.

\bibitem{sayas2013retarded}
Francisco-Javier Sayas.
\newblock Retarded potentials and time domain boundary integral equations.
\newblock 2013.

\bibitem{siegel2018local}
Michael Siegel and Anna-Karin Tornberg.
\newblock A local target specific quadrature by expansion method for evaluation
  of layer potentials in 3d.
\newblock {\em Journal of Computational Physics}, 364:365--392, 2018.

\bibitem{taskinen2006current}
Matti Taskinen and Pasi Yla-Oijala.
\newblock Current and charge integral equation formulation.
\newblock {\em IEEE Transactions on Antennas and Propagation}, 54(1):58--67,
  2006.

\bibitem{yilmaz2004time}
Ali~E Yilmaz, Jian-Ming Jin, and Eric Michielssen.
\newblock Time domain adaptive integral method for surface integral equations.
\newblock {\em IEEE Transactions on Antennas and Propagation},
  52(10):2692--2708, 2004.

\bibitem{young2012high}
P~Young, Sijia Hao, and Per-Gunnar Martinsson.
\newblock A high-order nystr{\"o}m discretization scheme for boundary integral
  equations defined on rotationally symmetric surfaces.
\newblock {\em Journal of Computational Physics}, 231(11):4142--4159, 2012.

\end{thebibliography}

\end{document}